\documentclass{article}

\usepackage[textheight=9in,
textwidth=6.5in,
top=1in,
headheight=14pt,
headsep=25pt,
footskip=30pt]{geometry}
\usepackage{hyperref}
\hypersetup{
	colorlinks=true,
	linkcolor=black,
	anchorcolor=black,
	citecolor=black,
	urlcolor=black,
	pdftitle={Multiple Bayesian Filtering as Message Passing},
	pdfauthor={Giorgio M. Vitetta, Pasquale Di Viesti, Emilio Sirignano and Francesco Montorsi},
	pdfkeywords={Hidden Markov Model, Factor Graph, Particle Filter, Sum-Product Algorithm, Marginalized Particle Filter, Kalman Filter, Multiple Particle Filtering}}

\usepackage{amssymb}
\usepackage{amsfonts}
\usepackage{amsmath}
\usepackage{eurosym}
\usepackage{amssymb}
\usepackage{amsfonts}
\usepackage{amsmath}
\usepackage{graphicx}
\graphicspath{{Figures/}}
\usepackage[ruled]{algorithm2e}
\usepackage{subcaption}
\usepackage{xcolor}
\usepackage{titlesec}
\usepackage[titletoc,toc,title]{appendix}
\usepackage[noadjust]{cite}

\DeclareRobustCommand{\cev}[1]{\reflectbox{\ensuremath{\vec{\reflectbox{\ensuremath{#1}}}}}}
\begin{document}

\title{Double Bayesian Smoothing as Message Passing}
\maketitle

\begin{abstract}
Recently, a novel method for developing filtering algorithms, based on the
interconnection of two Bayesian filters and called double Bayesian
filtering, has been proposed. In this manuscript we show that the same
conceptual approach can be exploited to devise a new smoothing method,
called double Bayesian smoothing. A double Bayesian smoother combines a
double Bayesian filter, employed in its forward pass, with the
interconnection of two backward information filters used in its backward
pass. As a specific application of our general method, a detailed derivation
of double Bayesian smoothing algorithms for conditionally linear Gaussian
systems is illustrated. Numerical results for two specific dynamic systems
evidence that these algorithms can achieve a better complexity-accuracy
tradeoff and tracking capability than other smoothing techniques recently
appeared in the literature.
\end{abstract}

\begin{center}
	\begin{tabular}{ccc}
	Pasquale Di Viesti$^\dagger$ & Giorgio M. Vitetta$^\dagger$ &Emilio Sirignano$^\dagger$  \\\vspace{.5cm}
	\href{mailto:pasquale.diviesti@unimore.it}{pasquale.diviesti@unimore.it} &
	\href{mailto:giorgio.vitetta@unimore.it}{giorgio.vitetta@unimore.it} &
	\href{mailto:emilio.sirignano@unimore.it}{emilio.sirignano@unimore.it}\\
	\end{tabular}
\bigskip

$^\dagger${Dept. of Engineering "Enzo Ferrari", University of Modena and Reggio Emilia.}
\end{center}

\bigskip

\textbf{Keywords:} Hidden Markov Model, Smoothing, Factor Graph, Particle Filter, Kalman
Filter, Sum-Product Algorithm.

\bigskip

\section{Introduction\label{sec:intro}}

The problem of \emph{Bayesian} \emph{smoothing} for a \emph{state space model%
} (SSM) concerns the development of recursive algorithms able to estimate
the \emph{probability density function} (pdf) of the model state on a given
observation interval, given a batch of noisy measurements acquired over it 
\cite{Anderson_1979}, \cite{Sar_2013}; the estimated pdf is known as a \emph{%
	smoothed} or \emph{smoothing} pdf. Two general methods are available in the
literature for recursively calculating smoothing densities;\ they are known
as the \emph{forward filtering-backward smoothing recursion} (e.g., see \cite%
{Doucet_2000} and \cite{Kitagawa_1987}) and the method based on the \emph{%
	two-filter smoothing formula} (e.g., see \cite{Kitagawa_1994} and \cite%
{Bresler_1986}). Both methods are based on the idea that the smoothing
densities can be computed by combining the predicted and/or filtered
densities generated by a Bayesian filtering method with the statistical
information produced in the backward pass by a different filtering method;
the latter method is \emph{paired with the first one} and, in the case of
the two-filter smoothing formula, is known as \emph{backward information
	filtering} (BIF). Unluckily, closed form solutions for Bayesian smoothing
are available for \emph{linear Gaussian} and \emph{linear Gaussian mixture}
models only \cite{Anderson_1979, Sar_2013, Vo_2012}. This has motivated the
development of various methods based on approximating smoothing densities in
different ways. For instance, the use of Gaussian approximations for the
smoothing densities and of sigma points techniques for solving moment
matching integrals has been investigated in \cite{Sar_2010,
	Kok_2016,Gar_2017}. Another class of methods (usually known as \emph{%
	particle smoothers}) is based on the exploitation of \emph{sequential Monte
	Carlo} techniques, i.e. on approximating smoothing densities through a set
of weighted particles (e.g., see \cite%
{Doucet_2000,Douc_2011,Kitagawa_1994,Kitagawa_1996,Godsill_2004,Lindsten_2013}
and references therein). Recently, substantial attention has been also paid
to the development of smoothing algorithms for the class of \emph{%
	conditionally linear Gaussian} SSMs \cite{Hot_2019,
	Briers_2010,Fong_2002,Lindsten_2016, Vitetta_2018}. In this case, the above
mentioned approximate methods can benefit from the so called \emph{%
	Rao-Blackwellization} technique, i.e. from the marginalisation of the linear
substructure of any conditionally linear Gaussian model; this can
significantly reduce the overall computational complexity of both
sigma-point based Gaussian smoothing \cite{Hot_2019} and particle smoothing 
\cite{Briers_2010,Fong_2002,Lindsten_2016, Vitetta_2018} (that is usually
known as \emph{Rao-Blackwellized particle smoothing}, RBPS, in this case).

In this manuscript, we propose a novel general method for the development of 
\emph{computationally efficient} particle smoothers. Our method exploits the
same conceptual approach illustrated in \cite{Vitetta_DiViesti_2019} in the
context of Bayesian filtering and dubbed \emph{multiple Bayesian filtering}.
That approach is based on the idea of developing new filtering algorithms
by: a) interconnecting multiple heterogeneous Bayesian filters; b)
representing the processing accomplished by each Bayesian filter and the
exchange of statistical information among distinct filters as a message
passing over a proper factor graph. In \cite{Vitetta_DiViesti_2019} the
exploitation of this approach has been investigated in detail for the case
in which two Bayesian filters are interconnected, i.e. \emph{dual Bayesian
	filtering} (DBF) is employed. Moreover, it has been shown that accurate and
computationally efficient DBF algorithms can be devised if the considered
SSM is conditionally linear Gaussian. In this manuscript, we show that, if
DBF is employed in the \emph{forward} pass of a smoothing method, a BIF
method, paired with DBF and based on the interconnection of two backward
information filters can be devised by following some simple rules. Similarly
as DBF, our derivation of such a BIF method, called \emph{double backward
	information filtering} (DBIF), is based on a \emph{graphical model}. Such a
graphical model allows us to show that: a) the pdfs computed in DBIF can be
represented as messages passed on it; b) all the expressions of the passed
messages can be derived by applying the same rule, namely the so called 
\emph{sum-product algorithm} (SPA) \cite{Loeliger_2007}, \cite%
{Kschischang_2001}, to it; c) iterative algorithms can be developed in a
natural fashion once the cycles it contains have been identified and the
order according to which messages are passed on them (i.e., the \emph{%
	message scheduling}) has been established; d) the statistical information
generated by a DBIF algorithm in the backward pass can be easily merged with
those produced by its paired DBF technique in the forward pass in order to
evaluate the required smoothed pdfs. To exemplify the usefulness of the
resulting smoothing method, based on the combination of DBF and DBIF, and
called \emph{double Bayesian smoothing} (DBS), the two DBF algorithms
proposed in \cite{Vitetta_DiViesti_2019} for the class of conditionally
linear Gaussian SSMs are taken into consideration, and the BITF algorithm
paired with each of them and a simplified version of it are derived. This
leads to the development of four new DBS algorithms, two generating an
estimate of the joint smoothing density over the whole observation interval,
the other two an estimate of the marginal smoothing densities over the same
interval. Our computer simulations for two specific conditionally linear
Gaussian\ SSMs evidence that, in the first case, the derived DBS algorithms
perform very closely to the RBPS technique proposed in \cite{Lindsten_2016}
and to the particle smoothers devised in \cite{Vitetta_2018}, but at lower
computational cost and time. In the second case, instead, two of the devised
DBS techniques represent the only technically useful options, thanks to
their good tracking capability. In fact, such techniques are able to operate
reliably even when their competitors diverge in the forward pass.

It is worth stressing that the\ technical contribution provided by this
manuscript represents a significant advancement with respect to the
application of factor graph theory to particle smoothing illustrated in \cite%
{Vitetta_2018}. In fact, in that manuscript, we also focus on conditionally
linear Gaussian models, but assume that the forward pass is accomplished by 
\emph{marginalized particle filtering} (MPF; also known as \emph{%
	Rao-Blackwellized particle filtering}); in other words, Bayesian filtering
is based on the interconnection of a particle filter with a \emph{bank} 
\emph{of Kalman filters}. In this manuscript, instead, the general method we
propose applies to a couple of arbitrary interconnected Bayesian filters.
Moreover, the specific smoothing algorithms we derive assume that the
forward pass is carried out by a filtering algorithm based on the
interconnection of a particle filter with a \emph{single} extended Kalman
filter.

The remaining part of this manuscript is organized as follows. In Section %
\ref{sec:Factorgraphs}, a general graphical model, on which the processing
accomplished in DBF, DBIF, and DBS is based, is illustrated. In Section \ref%
{Graph_mod_CLG}, a specific instance of the graphical model illustrated in
the previous section is developed under the assumptions that the filters
employed in the forward pass are an extended Kalman filter and a particle
filter, and that the considered SSM is conditionally linear Gaussian. Then,
the scheduling and the computation of the messages passed over this model
are analysed in detail and new DBS algorithms are devised. The differences
and similarities between these algorithms and other known smoothing
techniques are analysed in Section \ref{Comparison}. A comparison, in terms
of accuracy, computational complexity, and execution time, between the
proposed techniques and three smoothers recently appeared in the literature,
is provided in Section \ref{num_results} for two conditionally linear
Gaussian SSMs. Finally, some conclusions are offered in Section \ref%
{sec:conc}.

\emph{Notations}: The same notation as refs. \cite%
{Vitetta_2018,Vitetta_DiViesti_2019} and \cite{Vitetta_2019} is adopted.

\section{Graphical Model for a Couple of Interconnected Bayesian Information
	Filters and Message Passing on it\label{sec:Factorgraphs}}

In this manuscript, we consider a discrete-time SSM whose $D-$dimensional 
\emph{hidden state} in the $k-$th interval is denoted $\mathbf{x}%
_{k}\triangleq \lbrack x_{0,k},x_{1,k},...,$ $x_{D-1,k}]^{T}$, and whose 
\emph{state update} and \emph{measurement models} are expressed by 
\begin{equation}
	\mathbf{x}_{k+1}=\mathbf{f}_{k}\left( \mathbf{x}_{k}\right) +\mathbf{w}_{k}
	\label{eq:X_update}
\end{equation}%
and%
\begin{eqnarray}
	\mathbf{y}_{k} &\triangleq &[y_{0,k},y_{1,k},...,y_{P-1,k}]^{T}  \notag \\
	&=&\mathbf{h}_{k}\left( \mathbf{x}_{k}\right) +\mathbf{e}_{k},
	\label{meas_mod}
\end{eqnarray}%
respectively, with $k=1$, $2$, $...$, $T$. Here, $\mathbf{f}_{k}\left( 
\mathbf{x}_{k}\right) $ ($\mathbf{h}_{k}\left( \mathbf{x}_{k}\right) $) is a
time-varying $D-$dimensional ($P-$dimensional) real function, $T$ is the
duration of the observation interval and $\mathbf{w}_{k}$ ($\mathbf{e}_{k}$)
is the $k-$th element of the process (measurement) noise sequence $\left\{ 
\mathbf{w}_{k}\right\} $ ($\left\{ \mathbf{e}_{k}\right\} $); this sequence
consists of $D-$dimensional ($P-$dimensional) \emph{independent and
	identically distributed} (iid) Gaussian noise vectors, each characterized by
a zero mean and a covariance matrix $\mathbf{C}_{w}$ ($\mathbf{C}_{e}$).
Moreover, statistical independence between $\left\{ \mathbf{e}_{k}\right\} $
and $\{\mathbf{w}_{k}\}$ is assumed.

From a statistical viewpoint, a complete statistical description of the
considered SSM is provided by the pdf $f(\mathbf{x}_{1})$ of its initial
state, its \emph{Markov model} $f(\mathbf{x}_{k+1}|\mathbf{x}_{k})$ and its 
\emph{observation model} $f(\mathbf{y}_{k}|\mathbf{x}_{k})$ for any $k$; the
first pdf is assumed to be known, whereas the last two pdfs can be easily
derived from Eq. (\ref{eq:X_update}) and Eq. (\ref{meas_mod}), respectively.

In the following, we focus on the problem of developing novel smoothing
algorithms and, in particular, algorithms for the estimation of the \emph{%
	joint} \emph{smoothed pdf} $f(\mathbf{x}_{1:T}|\mathbf{y}_{1:T})$ (problem 
\textbf{P.1}) and the sequence of \emph{marginal} \emph{smoothed pdfs\ }$\{f(%
\mathbf{x}_{k}|\mathbf{y}_{1:T}),\,k=1,2,...,T\}$ (problem \textbf{P.2});
here, $\mathbf{y}_{1:T}\triangleq \left[ \mathbf{y}_{1}^{T},\mathbf{y}%
_{2}^{T},...,\mathbf{y}_{T}^{T}\right] ^{T}$ is a $P\cdot T-$dimensional
vector. Note that, in principle, once problem \textbf{P.1} is solved,
problem \textbf{P.2} can be easily tackled; in fact, if the joint pdf $f(%
\mathbf{x}_{1:T}|\mathbf{y}_{1:T})$ is known, all the posterior pdfs $\{f(%
\mathbf{x}_{k}|\mathbf{y}_{1:T})\}$ can be evaluated by \emph{marginalization%
}.

The development of our smoothing algorithms is mainly based on the graphical
approach illustrated in our previous manuscripts \cite[Sec. III]%
{Vitetta_2018}, \cite[Sec. II]{Vitetta_DiViesti_2019} and \cite[Sec. III]%
{Vitetta_2019} for Bayesian filtering and smoothing. This approach consists
in the following steps:

1. The state vector $\mathbf{x}_{k}$ is partitioned in two substates,
denoted $\mathbf{x}_{k}^{(1)}$ and $\mathbf{x}_{k}^{(2)}$ and having sizes $%
D_{1}$ and $D_{2}=D-D_{1}$, respectively. Note that, if $\mathbf{\bar{x}}%
_{k}^{(i)}$ represents the portion of $\mathbf{x}_{k}$ not included in $%
\mathbf{x}_{k}^{(i)}$ (with $i=1$ and $2$), our assumptions entail that $%
\mathbf{\bar{x}}_{k}^{(1)}=\mathbf{x}_{k}^{(2)}$ and $\mathbf{\bar{x}}%
_{k}^{(2)}=\mathbf{x}_{k}^{(1)}$.

2. A \emph{sub-graph} that allows to represent both Bayesian filtering and
BIF for the substate $\mathbf{x}_{k}^{(i)}$ (with $i=1$ and $2$) as message
passing algorithms on it is developed, under the assumption that the
complementary substate $\mathbf{\bar{x}}_{k}^{(i)}$ is statistically known.
This means that filtered and predicted densities of $\mathbf{x}_{k}^{(i)}$
are represented as messages passed on the edges of this sub-graph and the
rules for computing them result from the application of the SPA to it.

3. The two sub-graphs devised in the previous step (one referring to $%
\mathbf{x}_{k}^{(1)}$, the other one to $\mathbf{x}_{k}^{(2)}$) are \emph{%
	interconnected}, so that a single graphical model referring to the whole
state $\mathbf{x}_{k}$ is obtained.

4. Algorithms for Bayesian filtering and BIF for the whole state $\mathbf{x}%
_{k}$ are derived by applying the SPA to the graphical model obtained in the
previous step.

\begin{figure}[tbp]
	\centering
	\includegraphics[width=0.6\textwidth]{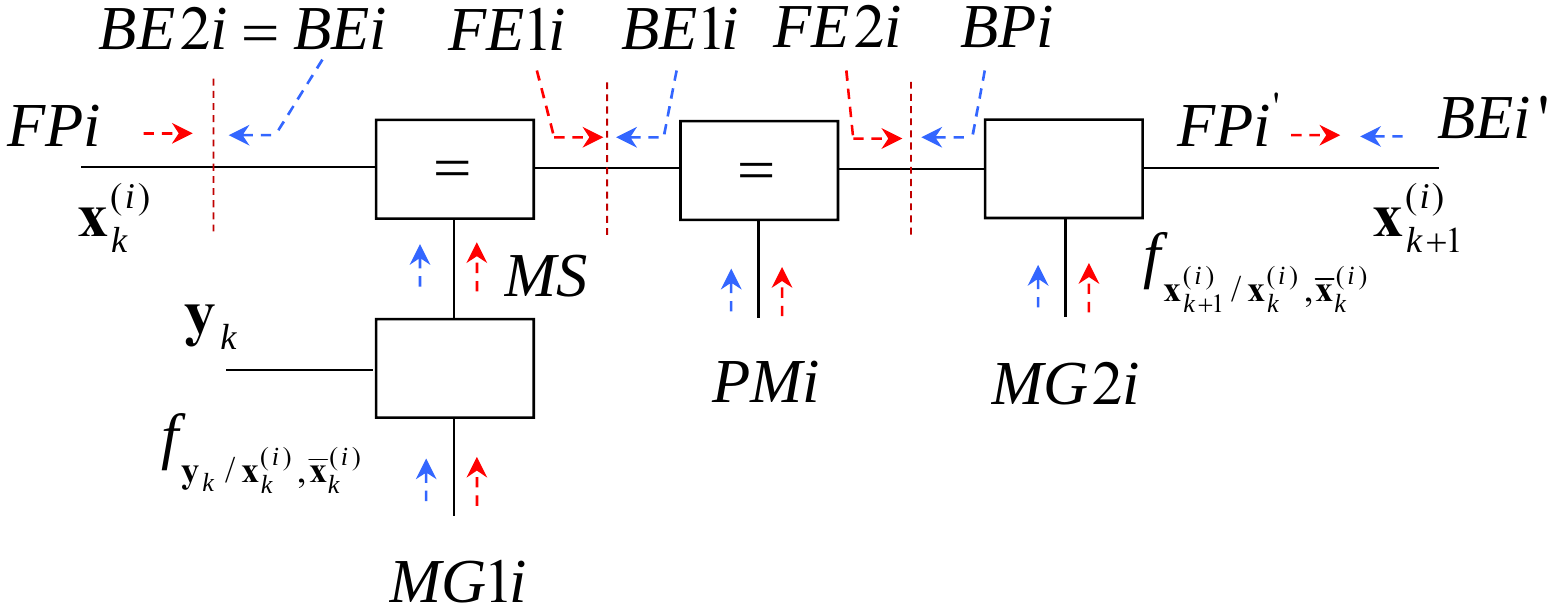}
	\caption{Factor graph involved in the $k-$th ($(T-k)-$th) recursion of
		Bayesian filtering (BIF) for the substate $\mathbf{x}_{k}^{(i)}$ and forward
		(backward) message passing on it. The flow of messages in the forward
		(backward) pass are indicated by red (blue) arrows, respectively; the brown
		vertical lines cutting each graph identify the partitioning associated with
		formulas (\protect\ref{factorisation3a}) (left cut), (\protect\ref%
		{factorisation3b}) (central cut) and (\protect\ref{factorisation3c}) (right
		cut). The messages $\vec{m}_{\mathrm{fp}}(\mathbf{x}_{k}^{(i)})$, $\cev{m}_{%
			\mathrm{bp}}(\mathbf{x}_{k}^{(i)})$, $\vec{m}_{\mathrm{fp}}(\mathbf{x}%
		_{k+1}^{(i)})$, $\cev{m}_{\mathrm{be}}(\mathbf{x}_{k+1}^{(i)})$, $m_{\mathrm{%
				ms}}(\mathbf{x}_{k}^{(i)})$, $m_{\mathrm{mg}l}(\mathbf{\bar{x}}_{k}^{(i)})$, 
		$m_{\mathrm{pm}}(\mathbf{x}_{k}^{(i)})$, $\vec{m}_{\mathrm{fe}l}(\mathbf{x}%
		_{k}^{(i)})$ and $\cev{m}_{\mathrm{be}l}(\mathbf{x}_{k}^{(i)})$ are denoted $FPi$, $BPi$, $FPi^{^{\prime }}$, $BEi^{^{\prime }}$, $MSi$, $MGli$, $PMi$, $FEli$ and $BEli$ respectively, to ease reading.}
	\label{Fig_1}
\end{figure}

Let us analyse now the steps 2.-4. in more detail. As far as \textbf{step 2.}
is concerned, the sub-graph devised for the substate $\mathbf{x}_{k}^{(i)}$
is based on the same principles illustrated in our manuscripts cited above
(in particular, ref. \cite{Vitetta_DiViesti_2019}) and is illustrated in
Fig. \ref{Fig_1}. The $k-$th recursion (with $k=1,2,...,T$) of Bayesian
filtering for the sub-state $\mathbf{x}_{k}^{(i)}$ is represented as a \emph{%
	forward} message passing on this factor graph, that involves the Markov
model $f(\mathbf{x}_{k+1}^{(i)}|\mathbf{x}_{k}^{(i)},\mathbf{\bar{x}}%
_{k}^{(i)})$ and the observation model $f(\mathbf{y}_{k}|\mathbf{x}%
_{k}^{(i)},\mathbf{\bar{x}}_{k}^{(i)})$. This allows to compute the messages 
$\vec{m}_{\mathrm{fe}1}(\mathbf{x}_{k}^{(i)})$, $\vec{m}_{\mathrm{fe}2}(%
\mathbf{x}_{k}^{(i)})$ and $\vec{m}_{\mathrm{fp}}(\mathbf{x}_{k+1}^{(i)})$,
that convey the \emph{first} \emph{filtered} pdf of $\mathbf{x}_{k}^{(i)}$,
the \emph{second} \emph{filtered} pdf of $\mathbf{x}_{k}^{(i)}$ and the 
\emph{predicted} pdf of $\mathbf{x}_{k+1}^{(i)}$, respectively, on the basis
of the messages $\vec{m}_{\mathrm{fp}}(\mathbf{x}_{k}^{(i)})$, $m_{\mathrm{ms%
}}(\mathbf{x}_{k}^{(i)})$ and $m_{\mathrm{pm}}(\mathbf{x}_{k}^{(i)})$; the
last three messages represent the predicted pdf of $\mathbf{x}_{k}^{(i)}$
evaluated in the previous (i.e., in the $(k-1)-$th) recursion of Bayesian
filtering, and the messages conveying the \emph{measurement} and the \emph{%
	pseudo-measurement} information, respectively, available in the $k-$%
recursion. The considered filtering algorithm requires the availability of
the messages $m_{\mathrm{pm}}(\mathbf{x}_{k}^{(i)})$, $m_{\mathrm{mg}1}(%
\mathbf{\bar{x}}_{k}^{(i)})$, $m_{\mathrm{mg}2}(\mathbf{\bar{x}}_{k}^{(i)})$%
, that are computed on the basis of external statistical information. The
presence of the messages $m_{\mathrm{mg}1}(\mathbf{\bar{x}}_{k}^{(i)})$ and $%
m_{\mathrm{mg}2}(\mathbf{\bar{x}}_{k}^{(i)})$ is due the fact that the
substate $\mathbf{\bar{x}}_{k}^{(i)}$ represents a \emph{nuisance state} for
the considered filtering algorithm; in fact, these messages convey filtered
(or predicted) pdfs of $\mathbf{\bar{x}}_{k}^{(i)}$ and are employed to
integrate out the dependence of the pdfs $f(\mathbf{y}_{k}|\mathbf{x}%
_{k}^{(i)},\mathbf{\bar{x}}_{k}^{(i)})$ and $f(\mathbf{x}_{k+1}^{(i)}|%
\mathbf{x}_{k}^{(i)},\mathbf{\bar{x}}_{k}^{(i)})$, respectively, on $\mathbf{%
	\bar{x}}_{k}^{(i)}$. Note also that these two messages are not necessarily
equal, since more refined information about $\mathbf{\bar{x}}_{k}^{(i)}$
could become available after that the message $m_{\mathrm{ms}}(\mathbf{x}%
_{k}^{(i)})$ has been computed. On the other hand, the message $m_{\mathrm{pm%
}}(\mathbf{x}_{k}^{(i)})$ conveys the statistical information provided by a
pseudo-measurement\footnote{%
	Generally speaking, a pseudo-measurement is a \emph{fictitious} measurement
	that is computed on the basis of statistical information provided by a
	filtering algorithm different from the one benefiting from it.} about $%
\mathbf{x}_{k}^{(i)}$. In Fig. \ref{Fig_1}, following \cite[Sec. II]%
{Vitetta_DiViesti_2019}, it is assumed that the pseudo-measurement $\mathbf{z%
}_{k}^{(i)}$ is available in the estimation of $\mathbf{x}_{k}^{(i)}$ and
that $m_{\mathrm{pm}}(\mathbf{x}_{k}^{(i)})$ represents the pdf of $\mathbf{z%
}_{k}^{(i)}$ conditioned on $\mathbf{x}_{k}^{(i)}$, that is%
\begin{equation}
	m_{\mathrm{pm}}\left( \mathbf{x}_{k}^{(i)}\right) \triangleq f\left( \mathbf{%
		z}_{k}^{(i)}\left\vert \mathbf{x}_{k}^{(i)}\right. \right) .
	\label{eq_pm_Fi}
\end{equation}

The computation of the messages $\vec{m}_{\mathrm{fe}1}(\mathbf{x}%
_{k}^{(i)}) $, $\vec{m}_{\mathrm{fe}2}(\mathbf{x}_{k}^{(i)})$ and $\vec{m}_{%
	\mathrm{fp}}(\mathbf{x}_{k+1}^{(i)})$ on the basis of the messages $\vec{m}_{%
	\mathrm{fp}}(\mathbf{x}_{k}^{(i)})$, $m_{\mathrm{ms}}(\mathbf{x}_{k}^{(i)})$%
, $m_{\mathrm{pm}}(\mathbf{x}_{k}^{(i)})$, $m_{\mathrm{mg}1}(\mathbf{\bar{x}}%
_{k}^{(i)})$ and $m_{\mathrm{mg}2}(\mathbf{\bar{x}}_{k}^{(i)})$ is based on
the two simple rules illustrated in \cite[Figs. 8-a) and 8-b), p. 1535]%
{Vitetta_2019} and can be summarized as follows. The first and second
filtered pdfs (i.e., the first and the second forward estimates) of $\mathbf{%
	x}_{k}^{(i)}$ are evaluated as 
\begin{equation}
	\vec{m}_{\mathrm{fe}1}\left( \mathbf{x}_{k}^{(i)}\right) =\vec{m}_{\mathrm{fp%
	}}\left( \mathbf{x}_{k}^{(i)}\right) m_{\mathrm{ms}}\left( \mathbf{x}%
	_{k}^{(i)}\right) ,  \label{eq_fe1_Fi}
\end{equation}%
and%
\begin{equation}
	\vec{m}_{\mathrm{fe}2}\left( \mathbf{x}_{k}^{(i)}\right) =\vec{m}_{\mathrm{fe%
		}1}\left( \mathbf{x}_{k}^{(i)}\right) m_{\mathrm{pm}}\left( \mathbf{x}%
	_{k}^{(i)}\right) ,  \label{eq_fe2_Fi}
\end{equation}%
respectively, where 
\begin{equation}
	m_{\mathrm{ms}}\left( \mathbf{x}_{k}^{(i)}\right) \triangleq \int f\left( 
	\mathbf{y}_{k}\left\vert \mathbf{x}_{k}^{(i)},\mathbf{\bar{x}}%
	_{k}^{(i)}\right. \right) \,m_{\mathrm{mg}1}\left( \mathbf{\bar{x}}%
	_{k}^{(i)}\right) \,d\mathbf{\bar{x}}_{k}^{(i)}  \label{eq_ms_Fi}
\end{equation}%
and $m_{\mathrm{pm}}(\mathbf{x}_{k}^{(i)})$ is defined in Eq. (\ref{eq_pm_Fi}%
). Equations (\ref{eq_fe1_Fi})-(\ref{eq_ms_Fi}) describe the processing
accomplished in the \emph{measurement update} of the considered recursion.
This is followed by the \emph{time update}, in which the new predicted pdf
(i.e., the new forward prediction) 
\begin{eqnarray}
	\vec{m}_{\mathrm{fp}}\left( \mathbf{x}_{k+1}^{(i)}\right) &=&\int \int
	f\left( \mathbf{x}_{k+1}^{(i)}\left\vert \mathbf{x}_{k}^{(i)},\mathbf{\bar{x}%
	}_{k}^{(i)}\right. \right) \vec{m}_{\mathrm{fe}2}\left( \mathbf{x}%
	_{k}^{(i)}\right)  \notag \\
	&&\cdot m_{\mathrm{mg}2}\left( \mathbf{\bar{x}}_{k}^{(i)}\right) d\mathbf{x}%
	_{k}\,d\mathbf{\bar{x}}_{k}^{(i)},  \label{eq_fpnew_Fi}
\end{eqnarray}%
is computed. The message passing procedure described above is initialised by
setting $\vec{m}_{\mathrm{fp}}(\mathbf{x}_{1}^{(i)})=f(\mathbf{x}_{1}^{(i)})$
(where $f(\mathbf{x}_{1}^{(i)})$ is the pdf resulting from the
marginalization of $f(\mathbf{x}_{1})$ with respect to $\mathbf{\bar{x}}%
_{1}^{(i)}$) in the first recursion and is run for $k=1,2,...,T$. Once this
procedure is over, BIF is executed for the substate $\mathbf{x}_{k}^{(i)}$;
its $(T-k)-$th recursion (with $k=T-1,T-2,...,1$) can be represented as a 
\emph{backward} message passing on the factor graph shown in Fig. \ref{Fig_1}%
. In this case, the messages $\cev{m}_{\mathrm{bp}}(\mathbf{x}_{k}^{(i)})$, $%
\cev{m}_{\mathrm{be}1}(\mathbf{x}_{k}^{(i)})$, $\cev{m}_{\mathrm{be}2}(%
\mathbf{x}_{k}^{(i)})=\cev{m}_{\mathrm{be}}(\mathbf{x}_{k}^{(i)})$, that
convey the \emph{backward predicted} pdf of $\mathbf{x}_{k}^{(i)}$, the 
\emph{first backward} \emph{filtered} pdf of $\mathbf{x}_{k}^{(i)}$ and the 
\emph{second backward} \emph{filtered} pdf of $\mathbf{x}_{k}^{(i)}$,
respectively, are evaluated on the basis of the messages $\cev{m}_{\mathrm{be%
}}(\mathbf{x}_{k+1}^{(i)})$, $m_{\mathrm{pm}}(\mathbf{x}_{k}^{(i)})$ and $m_{%
	\mathrm{ms}}(\mathbf{x}_{k}^{(i)})$, respectively; note that $\cev{m}_{%
	\mathrm{be}}(\mathbf{x}_{k+1}^{(i)})$ represents the backward filtered pdf
of $\mathbf{x}_{k}^{(i)}$ computed in the previous (i.e., in the $(T-(k+1))-$%
th) recursion of BIF. Moreover, the first and second backward filtered pdfs
of $\mathbf{x}_{k}^{(i)}$ are evaluated as (see Fig. \ref{Fig_1}) 
\begin{equation}
	\cev{m}_{\mathrm{be}1}\left( \mathbf{x}_{k}^{(i)}\right) =\cev{m}_{\mathrm{bp%
	}}\left( \mathbf{x}_{k}^{(i)}\right) \,m_{\mathrm{pm}}\left( \mathbf{x}%
	_{k}^{(i)}\right) ,  \label{eq_be1_Fi}
\end{equation}%
and%
\begin{equation}
	\cev{m}_{\mathrm{be}2}\left( \mathbf{x}_{k}^{(i)}\right) =\cev{m}_{\mathrm{be%
	}}\left( \mathbf{x}_{k}^{(i)}\right) =\cev{m}_{\mathrm{be}1}\left( \mathbf{x}%
	_{k}^{(i)}\right) \,m_{\mathrm{ms}}\left( \mathbf{x}_{k}^{(i)}\right) ,
	\label{eq_be2_Fi}
\end{equation}%
respectively, where $m_{\mathrm{pm}}(\mathbf{x}_{k}^{(i)})$ and $m_{\mathrm{%
		ms}}(\mathbf{x}_{k}^{(i)})$ are still expressed by Eq. (\ref{eq_pm_Fi}) and
Eq. (\ref{eq_ms_Fi}), respectively. The BIF message passing is initialised
by setting $\cev{m}_{\mathrm{be}}(\mathbf{x}_{T}^{(i)})=m_{\mathrm{fe}}(%
\mathbf{x}_{T}^{(i)})$ in its first recursion and is run for $%
k=T-1,T-2,...,1 $. Once the backward pass is over, a solution to problem 
\textbf{P.2} becomes available for the substate $\mathbf{x}_{k}^{(i)}$,
since the marginal smoothed pdf $f(\mathbf{x}_{k}^{(i)},\mathbf{y}_{1:T},%
\mathbf{z}_{1:T}^{(i)})$ (where $\mathbf{z}_{1:T}^{(i)}$ is the $P\cdot T-$%
dimensional vector resulting from the ordered concatenation of the all the
observed pseudo-measurements $\{\mathbf{z}_{k}^{(i)}\}$) can be evaluated as%
\footnote{%
	Note that, similarly as refs. \cite{Vitetta_2018} and \cite{Vitetta_2019}, a 
	\emph{joint} smoothed pdf is considered here in place of the corresponding 
	\emph{posterior} pdf.}%
\begin{eqnarray}
	f\left( \mathbf{x}_{k}^{(i)},\mathbf{y}_{1:T},\mathbf{z}_{1:T}^{(i)}\right)
	&=&\vec{m}_{\mathrm{fp}}\left( \mathbf{x}_{k}^{(i)}\right) \cev{m}_{\mathrm{%
			be}2}\left( \mathbf{x}_{k}^{(i)}\right)  \label{factorisation3a} \\
	&=&\vec{m}_{\mathrm{fe}1}\left( \mathbf{x}_{k}^{(i)}\right) \cev{m}_{\mathrm{%
			be}1}\left( \mathbf{x}_{k}^{(i)}\right)  \label{factorisation3b} \\
	&=&\vec{m}_{\mathrm{fe}2}\left( \mathbf{x}_{k}^{(i)}\right) \cev{m}_{\mathrm{%
			bp}}\left( \mathbf{x}_{k}^{(i)}\right) ,  \label{factorisation3c}
\end{eqnarray}%
with $k=1,2,...,T$. Note that, from a graphical viewpoint, formulas (\ref%
{factorisation3a})-(\ref{factorisation3c}) can be related with the three
different partitionings of the graph shown in Fig. \ref{Fig_1} (where a
specific partitioning is identified by a brown dashed vertical line cutting
the graph in two parts).

Given the graphical model represented in Fig. \ref{Fig_1}, \textbf{step 3.}
can be accomplished by adopting the same conceptual approach as \cite[Sec.
III]{Vitetta_2018} and \cite[Par. II-B]{Vitetta_DiViesti_2019}, where the
factor graphs on which smoothing and filtering, respectively, are based are
obtained by merging two sub-graphs, each referring to a distinct substate.
For this reason, in this case, the graphical model for the whole state $%
\mathbf{x}_{k}$ is obtained by interconnecting two distinct factor graphs,
each structured like the one shown in Fig. \ref{Fig_1}. In \cite[Par. II-B]%
{Vitetta_DiViesti_2019}, message passing on the resulting graph is described
in detail for the case of \emph{Bayesian filtering}. In this manuscript,
instead, our analysis of message passing concerns BIF\emph{\ and smoothing
	only}. The devised graph and the messages passed on it are shown in Fig. \ref%
{Fig_2}. Note that, in developing our graphical model, it has been assumed
that the \emph{smoothed} pdf referring to $\mathbf{x}_{k}^{(i)}$ (and
conveyed by the message $m_{\mathrm{sm}}(\mathbf{x}_{k}^{(i)})$) is computed
on the basis of Eq. (\ref{factorisation3a}), i.e. by merging the messages $%
\vec{m}_{\mathrm{fp}}(\mathbf{x}_{k}^{(i)})$ and $\cev{m}_{\mathrm{be}}(%
\mathbf{x}_{k}^{(i)})=\cev{m}_{\mathrm{be}2}(\mathbf{x}_{k}^{(i)})$.
Moreover, the following elements (identified by brown lines) have been added
to its $i-$th sub-graph (with $i=1$ and $2$): a) two equality nodes; b) the
block BIF$_{i}{\rightarrow }$BIF$_{j}$ for extracting useful information
from the messages computed on the $i-$th sub-graph and delivered to the $j-$%
th one. The former elements allow the $i-$th backward information filter to
generate copies of the messages $\cev{m}_{\mathrm{be}}(\mathbf{x}%
_{k+1}^{(i)})$ and $m_{\mathrm{sm}}(\mathbf{x}_{k}^{(i)})$, that are made
available to the other sub-graphs. In the latter element, instead, the
messages $m_{\mathrm{pm}}(\mathbf{x}_{k}^{(i)})$ (see Eq. (\ref{eq_pm_Fi}))
and $m_{\mathrm{mg}q}(\mathbf{\bar{x}}_{k}^{(i)})$ (with $q=1$ and $2$; see
Eqs. (\ref{eq_ms_Fi}) and (\ref{eq_fpnew_Fi})) are computed; note that this
block is connected to \emph{oriented} edges only, i.e. to edges on which the
flow of messages is unidirectional.

\begin{figure}[tbp]
	\centering
	\includegraphics[width=0.6\textwidth]{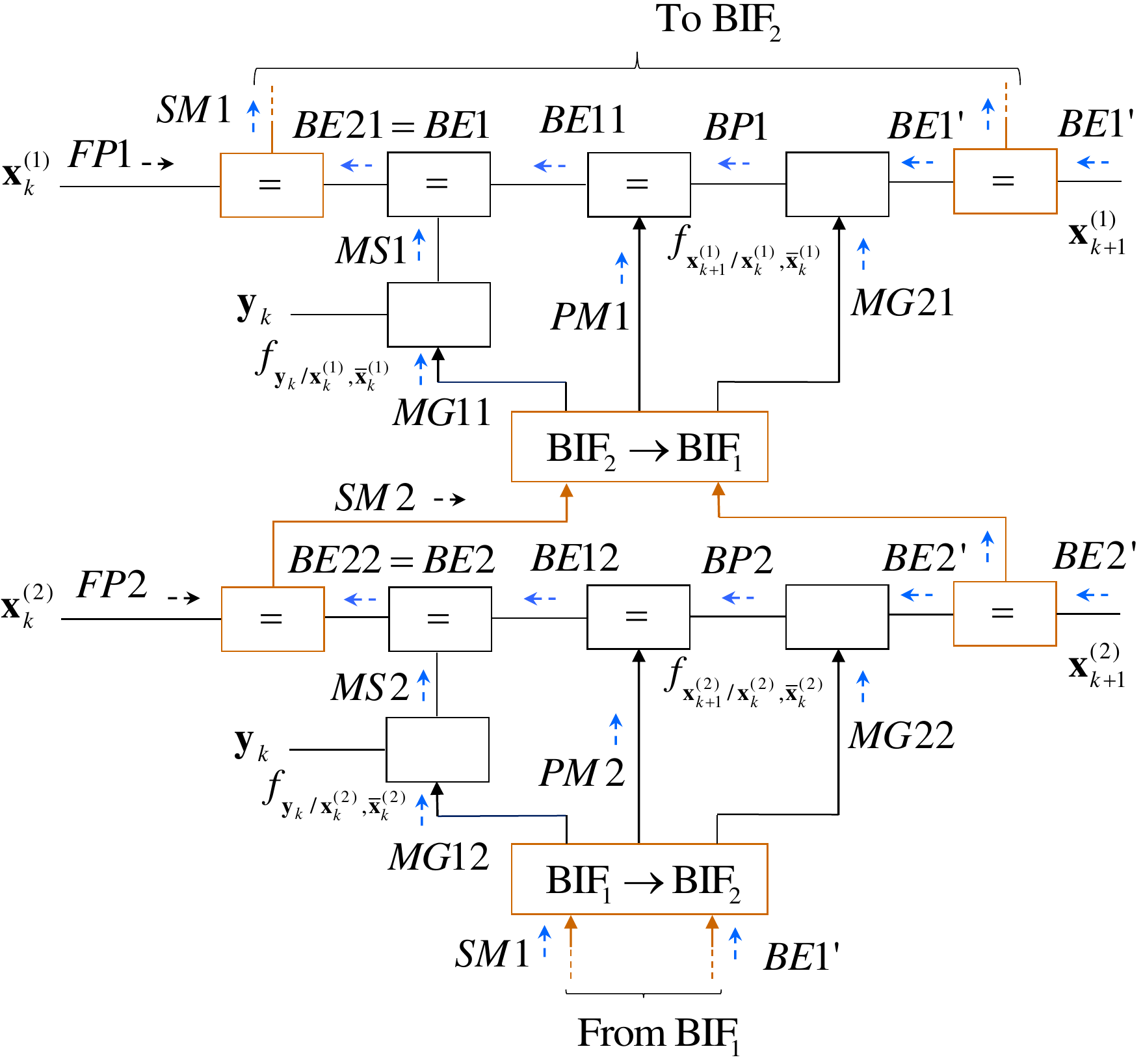}
	\caption{Graphical model based on the sub-graph shown in Fig. \protect\ref%
		{Fig_1} and referring to the interconnection of two backward information
		filters. The message computed in the backward (forward) pass are identified
		by blue (black) arrows. The message $m_{\mathrm{sm}}(\mathbf{x}_{k}^{(i)})$
		is denoted $SMi$ to ease reading.}
	\label{Fig_2}
\end{figure}

Given the graphical model represented in Fig. \ref{Fig_2}, \textbf{step 4}.
can be easily accomplished. In fact, recursive BIF and smoothing algorithms
can be derived by systematically applying the SPA to it after that a proper 
\emph{scheduling} has been established for message passing. In doing so, we
must always keep in mind that:\newline
1) Message passing on the $i-$th subgraph represents BIF/smoothing for the
substate $\mathbf{x}_{k}^{(i)}$; the exchange of messages between the
sub-graphs, instead, allows us to represent the interaction of two
interconnected BIF/smoothing algorithms in a effective and rigorous way.%
\newline
2) Different approximations can be used for the predicted/filtered/smoothed
pdfs computed in the message passing on each of the two sub-graphs and for
the involved Markov/observation models. For this reason, generally speaking,
the two interconnected filtering/BIF/smoothing algorithms are not required
to be of the same type.\newline
3) The $k-$th recursion of the overall BIF algorithm is fed by the backward
estimates $\cev{m}_{\mathrm{be}}(\mathbf{x}_{k+1}^{(1)})$ ($BE1^{\prime }$)
and $\cev{m}_{\mathrm{be}}(\mathbf{x}_{k+1}^{(2)})$ ($BE2^{\prime }$), and
generates the new backward predictions $\cev{m}_{\mathrm{bp}}(\mathbf{x}%
_{k}^{(1)})$ ($BP1$) and $\cev{m}_{\mathrm{bp}}(\mathbf{x}_{k}^{(2)})$ ($BP2$%
), and the two couples of filtered densities $\{(\cev{m}_{\mathrm{be}1}(%
\mathbf{x}_{k}^{(i)})$, $\cev{m}_{\mathrm{be}2}(\mathbf{x}_{k}^{(i)}))$, $i=1
$, $2\}$ ($\{BE1i$, $BE2i$, $i=1$, $2\}$). Moreover, merging the predicted
densities computed in the forward pass (i.e., the messages $\{FPi\}$) with
the second backward filtered densities (i.e., the messages $\{BE2i=BEi\}$)
allows us to generate the smoothed pdfs for each substate according to Eq. (%
\ref{factorisation3a}). However, a \emph{joint} filtered/smoothed density
for the whole state $\mathbf{x}_{k}$ is \emph{unavailable}.\newline
4) Specific algorithms are employed to compute the pseudo-measurement and
the nuisance substate pdfs in the BIF$_{i}{\rightarrow }$BIF$_{j}$ blocks
appearing in Fig. \ref{Fig_2}. These algorithms depend on the considered SSM
and on the selected message scheduling; for this reason, a general
description of their structure cannot be provided.\newline
5) The graphical model shown in Fig. \ref{Fig_2}, unlike the one illustrated
in Fig. \ref{Fig_1}, is \emph{not cycle free}. The presence of cycles raises
the problems of identifying all the messages that can be iteratively refined
and establishing the order according to which they are computed. Generally
speaking, iterative message passing on the devised graphical model involves
both the couple of measurement updates and the backward prediction
accomplished in each of the interconnected backward information filters. In
fact, this should allow each filter to progressively refine the nuisance
substate density employed in its second measurement update and backward
prediction, and improve the quality of the pseudo-measurements exploited in
its first measurement update. For this reason, if $n_{i}$ iterations are
run, the overall computational complexity of each recursion is multiplied by 
$n_{i}$.\newline

The final important issue about the graphical model devised for both
Bayesian filtering and BIF concerns the possible presence of \emph{redundancy%
}. In all the considerations illustrated above, \emph{disjoint} substates $%
\mathbf{x}_{k}^{(1)}$ and $\mathbf{x}_{k}^{(2)}$ have been assumed.
Actually, in ref. \cite{Vitetta_DiViesti_2019}, it has been shown that our
graphical approach can be also employed if the substates $\mathbf{x}%
_{k}^{(1)}$ and $\mathbf{x}_{k}^{(2)}$ cover $\mathbf{x}_{k}$, but do not
necessarily form a partition of it. In other words, some overlapping between
these two substates is admitted. When this occurs, the forward/backward
filtering algorithm run over the whole graphical model contains a form of 
\emph{redundancy}, since $N_{d}\triangleq D_{1}+D_{2}-D$ elements of the
state vector $\mathbf{x}_{k}$ are independently estimated by the
interconnected forward/backward filters. The parameter $N_{d}$ can be
considered as the \emph{degree of redundancy} characterizing the
filtering/smoothing algorithm. Moreover, in ref. \cite{Vitetta_DiViesti_2019}%
, it has been shown that the presence of redundancy in a Bayesian filtering
algorithm can significantly enhance its tracking capability (i.e., reduce
its probability of divergence); however, this result is obtained at the
price of an increased complexity with respect to the case in which the
interconnected filters are run over disjoint substates.

\section{Double Backward Information Filtering and Smoothing Algorithms for
	Conditionally Linear Gaussian State Space Models\label{Graph_mod_CLG}}

In this section we focus on the development of two new DBS algorithms for
conditionally linear Gaussian models. We first describe the graphical models
on which these algorithms are based; then, we provide a detailed description
of the computed messages and their scheduling in a specific case.

\subsection{Graphical Modelling}

In this paragraph, we focus on a specific instance of the graphical model
illustrated in Fig. \ref{Fig_2}, since we make the same specific choices as
ref. \cite{Vitetta_DiViesti_2019}\ for both the considered SSM and the two
Bayesian filters employed in the forward pass. For this reason, we assume
that: a) the SSM described by eqs. (\ref{eq:X_update})-(\ref{meas_mod}) is
conditionally linear Gaussian \cite{Lindsten_2016}, \cite{Vitetta_2019}, 
\cite{Schon_2005}, so that its state vector $\mathbf{x}_{k}$ can be
partitioned into its \emph{linear component} $\mathbf{x}_{k}^{(L)}$ and its 
\emph{nonlinear component} $\mathbf{x}_{k}^{(N)}$ (having sizes $D_{L}$ and $%
D_{N}$, respectively, with $D_{N}+D_{L}=D$); b) the dual Bayesian filter
employed in the forward pass results from the interconnection of an\emph{\
	extended Kalman filter} with a \emph{particle filter}\footnote{%
	In particular, a \emph{sequential importance resampling} filter is employed 
	\cite{Arulampalam_2002}.} (these filters are denoted F$_{1}$ and F$_{2}$,
respectively), as described in detail in ref. \cite{Vitetta_DiViesti_2019}.
As far as the last point is concerned, it is also worth mentioning that, on
the one hand, filter F$_{2}$ estimates the nonlinear state component only
(so that $\mathbf{x}_{k}^{(2)}=\mathbf{x}_{k}^{(N)}$ and $\mathbf{\bar{x}}%
_{k}^{(2)}=\mathbf{x}_{k}^{(L)}$) and approximates the predicted/filtered
densities of this component through a set of $N_{p}$ weighted particles. On
the other hand, filter F$_{1}$ employs a Gaussian approximation of all its
predicted/filtered densities, and works on the \emph{whole system state} or
on the \emph{linear state component}. In the first case (denoted \textbf{C.1}
in the following), we have that $\mathbf{x}_{k}^{(1)}=\mathbf{x}_{k}$ and $%
\mathbf{\bar{x}}_{k}^{(1)}$ is empty, so that both F$_{1}$ and F$_{2}$
estimate the nonlinear state component (for this reason, the corresponding
degree of redundancy in the developed smoothing algorithm is $N_{d}=D_{N}$);
in the second case (denoted \textbf{C.2} in the following), instead, $%
\mathbf{x}_{k}^{(1)}=\mathbf{x}_{k}^{(L)}$ and $\mathbf{\bar{x}}_{k}^{(1)}=%
\mathbf{x}_{k}^{(N)}$, so that filters F$_{1}$ and F$_{2}$ estimate \emph{%
	disjoint} substates (consequently, $N_{d}=0$).

\begin{figure}[tbp]
	\centering
	\includegraphics[width=0.6\textwidth]{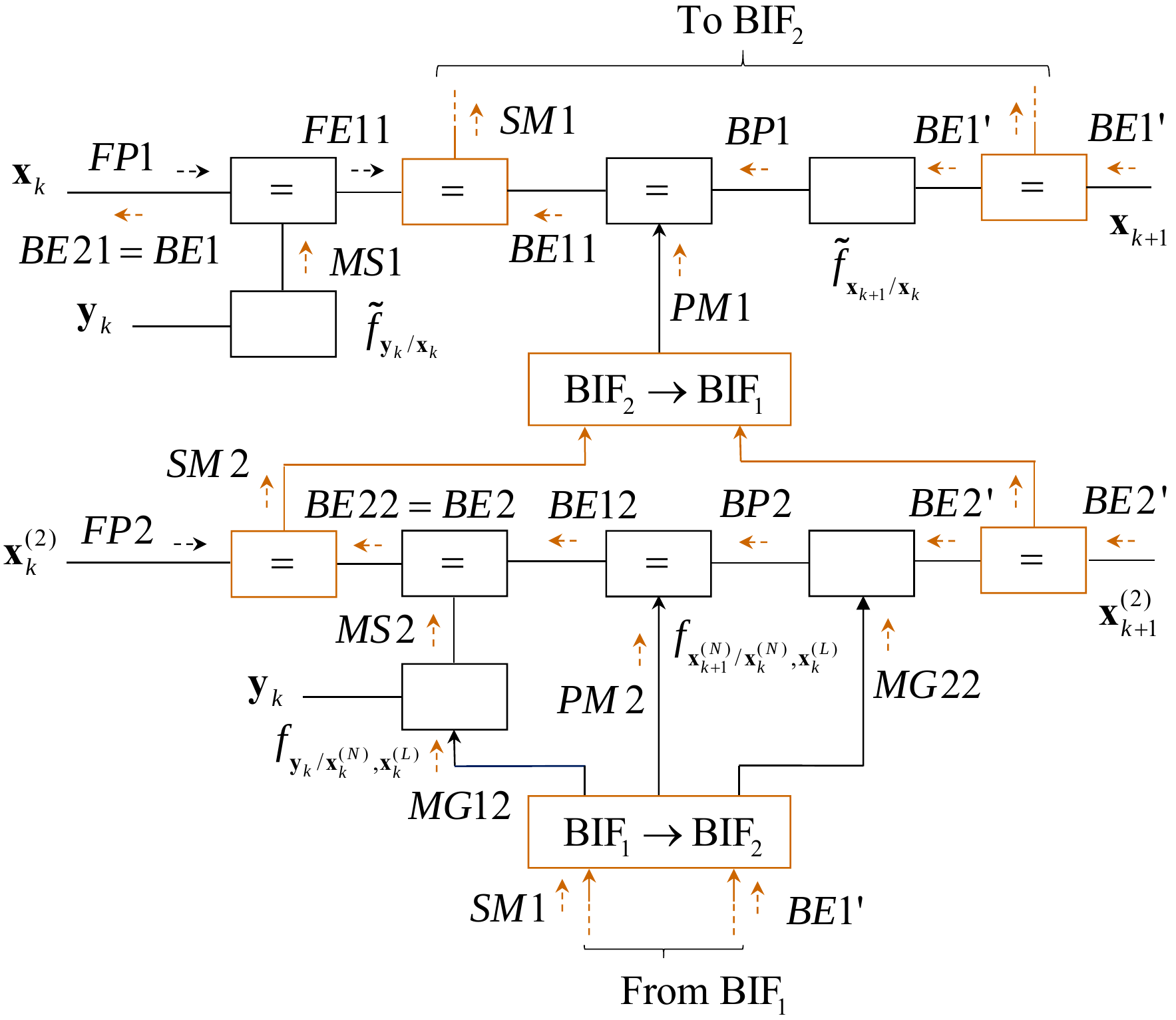}
	\caption{Graphical model referring to the interconnection of two backward
		information filters, one paired with an extended Kalman filter, the other
		one with a particle filter.}
	\label{Fig_3}
\end{figure}

Our selection of the forward filtering scheme has the following implications
on the developed DBIF scheme. The first backward information filter (denoted
BIF$_{1}$) is the backward filter associated with an extended Kalman filter
operating over on the \emph{whole system state} (case \textbf{C.1}) or on
the \emph{linear state component} (case \textbf{C.2}). The second backward
filter (denoted BIF$_{2}$), instead, is a backward filter associated with a
particle filter operating on the nonlinear state component only. In
practice, following \cite{Lindsten_2016,Vitetta_2018,Fong_2002}, BIF$_{2}$
is employed to update the weights of all the elements of the particle set
generated by filter F$_{2}$ in the forward pass. Then, based on the
graphical model shown in Fig. \ref{Fig_2}, the factor graph illustrated in
Fig. \ref{Fig_3} can be drawn for\ case \textbf{C.1}. It is important to
point out that:

1) The first backward information filter (BIF$_{1}$) is based on \emph{%
	linearised} (and, consequently, \emph{approximate}) Markov/measurement
models, whereas the second one (BIF$_{2}$) relies on \emph{exact} models, as
explained in more detail below. These models are the same as those employed
in ref. \cite{Vitetta_DiViesti_2019}.

2) Since the nuisance substate $\mathbf{\bar{x}}_{k}^{(1)}$ is empty, no
marginalization is required in BIF$_{1}$; for this reason, the messages $%
\{m_{\mathrm{mg}q}(\mathbf{\bar{x}}_{k}^{(1)})$; $q=1,2\}$ (i.e., $MG11$ and 
$MG21$) visible in Fig. \ref{Fig_2} do not appear in Fig. \ref{Fig_3}.
Moreover, the message $m_{\mathrm{sm}}(\mathbf{x}_{k}^{(1)})=m_{\mathrm{sm}}(%
\mathbf{x}_{k})$ is generated on the basis of Eq. (\ref{factorisation3b}),
instead of Eq. (\ref{factorisation3a}).

3) The backward filtered pdf $\cev{m}_{\mathrm{be}}(\mathbf{x}_{k+1}^{(2)})=%
\cev{m}_{\mathrm{be}}(\mathbf{x}_{k+1}^{(N)})$ and the smoothed pdf $m_{%
	\mathrm{sm}}(\mathbf{x}_{k}^{(2)})=m_{\mathrm{sm}}(\mathbf{x}_{k}^{(N)})$
(i.e., the messages $BE2^{^{\prime }}$ and $SM2$, respectively) feed the BIF$%
_{2}{\rightarrow }$BIF$_{1}$ block, where they are processed jointly to
generate the pseudo-measurement message $m_{\mathrm{pm}}(\mathbf{x}%
_{k}^{(1)})=m_{\mathrm{pm}}(\mathbf{x}_{k})$ ($PM1$) feeding filter F$_{1}$.
Similarly, the backward filtered pdf $\cev{m}_{\mathrm{be}}(\mathbf{x}%
_{k+1}^{(1)})=\cev{m}_{\mathrm{be}}(\mathbf{x}_{k+1})$ ($BE1\prime $) and
the smoothed pdf $m_{\mathrm{sm}}(\mathbf{x}_{k}^{(1)})=m_{\mathrm{sm}}(%
\mathbf{x}_{k})$ ($SM1$) feed the BIF$_{1}{\rightarrow }$BIF$_{2}$ block,
where the pseudo-measurement message $m_{\mathrm{pm}}(\mathbf{x}%
_{k}^{(2)})=m_{\mathrm{pm}}(\mathbf{x}_{k}^{(N)})$ ($PM2$) and the messages $%
\{m_{\mathrm{mg}q}(\mathbf{\bar{x}}_{k}^{(2)})=m_{\mathrm{mg}q}(\mathbf{x}%
_{k}^{(L)})$; $q=1,2\}$ (i.e., $MG12$ and $MG22$) are evaluated.

In the remaining part of this paragraph, we first provide various details
about the backward filters BIF$_{1}$ and BIF$_{2}$, and the way
pseudo-measurements are computed for each of them; then, we comment on how
the factor graph shown in Fig. \ref{Fig_3} should be modified if case 
\textbf{C.2} is considered.

BIF$_{1}$ - This backward filter is based on the \emph{linearized} versions
of Eqs. (\ref{eq:X_update}) and (\ref{meas_mod}), i.e. on the models (e.g.,
see \cite[pp. 194-195]{Anderson_1979} and \cite[Par. III-A]%
{Vitetta_DiViesti_2019})%
\begin{equation}
	\mathbf{x}_{k+1}=\mathbf{F}_{k}\,\mathbf{x}_{k}+\mathbf{u}_{k}+\mathbf{w}_{k}
	\label{state_up_approx}
\end{equation}%
and 
\begin{equation}
	\mathbf{y}_{k}=\mathbf{H}_{k}^{T}\,\mathbf{x}_{k}+\mathbf{v}_{k}+\mathbf{e}%
	_{k},  \label{meas_mod_approx}
\end{equation}%
respectively; here, $\mathbf{F}_{k}\triangleq \lbrack \partial \mathbf{f}%
_{k}\left( \mathbf{x}\right) /\partial \mathbf{x}]_{\mathbf{x=x}_{\mathrm{fe}%
		,k}}$, $\mathbf{u}_{k}\triangleq \mathbf{f}_{k}\left( \mathbf{x}_{\mathrm{fe}%
	,k}\right) -\mathbf{F}_{k}\,\mathbf{x}_{\mathrm{fe},k}$, $\mathbf{H}%
_{k}^{T}\triangleq \lbrack \partial \mathbf{h}_{k}\left( \mathbf{x}\right)
/\partial \mathbf{x}]_{\mathbf{x=x}_{\mathrm{fp},k}}$, $\mathbf{v}%
_{k}\triangleq \mathbf{h}_{k}\left( \mathbf{x}_{\mathrm{fp},k}\right) -%
\mathbf{H}_{k}^{T}\mathbf{x}_{\mathrm{fp},k}$ and $\mathbf{x}_{\mathrm{fp}%
	,k} $ ($\mathbf{x}_{\mathrm{fe},k}$) is the \emph{forward prediction }(\emph{%
	forward estimate}) of $\mathbf{x}_{k}$ computed by F$_{1}$ in its $(k-1)-$th
($k-$th) recursion. Consequently, the approximate models 
\begin{equation}
	\tilde{f}\left( \mathbf{x}_{k+1}\left\vert \mathbf{x}_{k}\right. \right) =%
	\mathcal{N}\left( \mathbf{x}_{k};\mathbf{F}_{k}\,\mathbf{x}_{k}+\mathbf{u}%
	_{k},\mathbf{C}_{w}\right)  \label{Markov_ekf_approx}
\end{equation}%
and 
\begin{equation}
	\tilde{f}\left( \mathbf{y}_{k}\left\vert \mathbf{x}_{k}\right. \right) =%
	\mathcal{N}\left( \mathbf{x}_{k};\mathbf{H}_{k}^{T}\,\mathbf{x}_{k}+\mathbf{v%
	}_{k},\mathbf{C}_{e}\right)  \label{Measurement_ekf_approx}
\end{equation}%
appear in the graphical model shown in Fig. \ref{Fig_3}.

BIF$_{2}$ - In developing this backward filter, the state vector $\mathbf{x}%
_{k}$ is represented as the ordered concatenation of its linear component $%
\mathbf{x}_{k}^{(L)}\triangleq \lbrack x_{0,k}^{(L)}$, $%
x_{1,k}^{(L)},...,x_{D_{L}-1,k}^{(L)}]^{T}$ and its nonlinear component $%
\mathbf{x}_{k}^{(N)}\triangleq \lbrack
x_{0,k}^{(N)},x_{1,k}^{(N)},...,x_{D_{N}-1,k}^{(N)}]^{T}$. Based on \cite[%
eq. (3)]{Vitetta_2019}, the Markov model%
\begin{equation}
	\mathbf{x}_{k+1}^{(N)}=\mathbf{A}_{k}^{(N)}\left( \mathbf{x}%
	_{k}^{(N)}\right) \mathbf{x}_{k}^{(L)}+\mathbf{f}_{k}^{(N)}\left( \mathbf{x}%
	_{k}^{(N)}\right) +\mathbf{w}_{k}^{(N)}  \label{eq:XN_update}
\end{equation}%
is adopted for the nonlinear state component (this model corresponds to the
last $D_{N}$ lines of Eq. (\ref{eq:X_update})); here, $\mathbf{f}_{k}^{(N)}(%
\mathbf{x}_{k}^{(N)})$ ($\mathbf{A}_{k}^{(N)}(\mathbf{x}_{k}^{(N)})$) is a
time-varying $D_{N}-$dimensional real function ($D_{N}\times D_{L}$ real
matrix) and $\mathbf{w}_{k}^{(N)}$ consists of the last $D_{N}$ elements of
the noise term $\mathbf{w}_{k}$ appearing in Eq. (\ref{eq:X_update}) (the
covariance matrix of $\mathbf{w}_{k}^{(N)}$ is denoted $\mathbf{C}_{w}^{(N)}$%
). Moreover, it is assumed that the observation model (\ref{meas_mod}) can
be put in the form (see \cite[eq. (31)]{Vitetta_DiViesti_2019} or \cite[eq.
(4)]{Vitetta_2019})%
\begin{equation}
	\mathbf{y}_{k}=\mathbf{g}_{k}\left( \mathbf{x}_{k}^{(N)}\right) +\mathbf{B}%
	_{k}\left( \mathbf{x}_{k}^{(N)}\right) \mathbf{x}_{k}^{(L)}+\mathbf{e}_{k},
	\label{eq:y_t}
\end{equation}%
where $\mathbf{g}_{k}(\mathbf{x}_{k}^{(N)})$ ($\mathbf{B}_{k}(\mathbf{x}%
_{k}^{(N)})$) is a time-varying $P-$dimensional real function ($P\times
D_{L} $ real matrix). Consequently, the considered backward filter is based
on the \emph{exact} pdfs 
\begin{eqnarray}
	&&f\left( \mathbf{x}_{k+1}^{(N)}\left\vert \mathbf{x}_{k}^{(N)},\mathbf{x}%
	_{k}^{(L)}\right. \right)  \notag \\
	&=&\mathcal{N}\left( \mathbf{x}_{k}^{(N)};\mathbf{A}_{k}^{(N)}\left( \mathbf{%
		x}_{k}^{(N)}\right) \mathbf{x}_{k}^{(L)}+\mathbf{f}_{k}^{(N)}\left( \mathbf{x%
	}_{k}^{(N)}\right) ,\mathbf{C}_{w}^{(N)}\right)  \notag \\
	&&  \label{f_x_N}
\end{eqnarray}%
and 
\begin{eqnarray}
	&&f\left( \mathbf{y}_{k}\left\vert \mathbf{x}_{k}^{(N)},\mathbf{x}%
	_{k}^{(L)}\right. \right)  \notag \\
	&=&\mathcal{N}\left( \mathbf{x}_{k};\mathbf{g}_{k}\left( \mathbf{x}%
	_{k}^{(N)}\right) +\mathbf{B}_{k}\left( \mathbf{x}_{k}^{(N)}\right) \mathbf{x%
	}_{k}^{(L)},\mathbf{C}_{e}\right) ,  \label{f_y_N}
\end{eqnarray}%
both appearing in the graphical model drawn in Fig. \ref{Fig_3}.

\emph{Computation of the pseudo-measurements for the first backward filter}
- Filter BIF$_{1}$ is fed by pseudo-measurement information about the \emph{%
	whole state} $\mathbf{x}_{k}$. The method for computing these information is
similar to the one illustrated in ref. \cite[Sects. III-IV]{Vitetta_2018}
and can be summarised as follows. The pseudo-measurements about the
nonlinear state component are represented by the $N_{p}$ particles conveyed
by the smoothed pdf $m_{\mathrm{sm}}(\mathbf{x}_{k}^{(N)})$ ($SM2$). On the
other hand, $N_{p}$ pseudo-measurements about the linear state component are
evaluated by means of the same method employed by \emph{marginalized
	particle filtering} (MPF) for this task. This method is based on the idea
that the random vector (see Eq. (\ref{eq:XN_update}))%
\begin{equation}
	\mathbf{z}_{k}^{(L)}\triangleq \mathbf{x}_{k+1}^{(N)}-\mathbf{f}%
	_{k}^{(N)}\left( \mathbf{x}_{k}^{(N)}\right) ,  \label{eq:z_L_l}
\end{equation}%
depending on the \emph{nonlinear state component} \emph{only}, must equal
the sum 
\begin{equation}
	\mathbf{A}_{k}^{(N)}\left( \mathbf{x}_{k}^{(N)}\right) \mathbf{x}_{k}^{(L)}+%
	\mathbf{w}_{k}^{(N)},  \label{eq:z_L_l_bis}
\end{equation}%
that depends on the \emph{linear state component}. For this reason, $N_{p}$
realizations of $\mathbf{z}_{k}^{(L)}$ (\ref{eq:z_L_l}) are computed in the
BIF$_{2}{\rightarrow }$BIF$_{1}$ block on the basis of the messages $\cev{m}%
_{\mathrm{be}}(\mathbf{x}_{k+1}^{(N)})$ ($BE2^{^{\prime }}$) and $m_{\mathrm{%
		sm}}(\mathbf{x}_{k}^{(N)})$, and are treated as measurements about $\mathbf{x%
}_{k}^{(L)}$.

\emph{Computation of the pseudo-measurements for the second backward filter}
- The messages $\cev{m}_{\mathrm{be}}(\mathbf{x}_{k+1})$ ($BE1^{\prime }$)
and $m_{\mathrm{sm}}(\mathbf{x}_{k})$ ($SM1$) feeding the BIF$_{1}{%
	\rightarrow }$BIF$_{2}$ block are employed for: a) generating the messages $%
\{m_{\mathrm{mg}q}(\mathbf{x}_{k}^{(L)})$; $q=1,2\}$ required to integrate
out the dependence of the state update and measurement models (i.e., of the
densities $f(\mathbf{x}_{k+1}^{(N)}|\mathbf{x}_{k}^{(N)}$, $\mathbf{x}%
_{k}^{(L)})$ (\ref{f_x_N}) and $f(\mathbf{y}_{k}|\mathbf{x}_{k}^{(N)},%
\mathbf{x}_{k}^{(L)})$ (\ref{f_y_N}), respectively) on the substate $\mathbf{%
	x}_{k}^{(L)}$; b) generating pseudo-measurement information about $\mathbf{x}%
_{k}^{(N)}$. As far as the last point is concerned, the approach we adopt is
the same as that developed for \emph{dual marginalized particle filtering }%
(dual MPF) in ref. \cite[Sec. V]{Vitetta_2019} and also adopted in particle
smoothing \cite[Sects. III-IV]{Vitetta_2018}. The approach relies on the
Markov model 
\begin{equation}
	\mathbf{x}_{k+1}^{(L)}=\mathbf{A}_{k}^{(L)}\left( \mathbf{x}%
	_{k}^{(N)}\right) \mathbf{x}_{k}^{(L)}+\mathbf{f}_{k}^{(L)}\left( \mathbf{x}%
	_{k}^{(N)}\right) +\mathbf{w}_{k}^{(L)},  \label{eq:XL1_update}
\end{equation}%
referring to the \emph{linear} state component (see \cite[eq. (1)]%
{Vitetta_2018} or \cite[eq. (3)]{Vitetta_2019}); in the last expression, $%
\mathbf{f}_{k}^{(L)}(\mathbf{x}_{k}^{(N)})$ ($\mathbf{A}_{k}^{(L)}(\mathbf{x}%
_{k}^{(N)})$) is a time-varying $D_{L}-$dimensional real function ($%
D_{L}\times D_{L}$ real matrix), and $\mathbf{w}_{k}^{(L)}$ consists of the
first $D_{L}$ elements of the noise term $\mathbf{w}_{k}$ appearing in (\ref%
{eq:X_update}) (the covariance matrix of $\mathbf{w}_{k}^{(L)}$ is denoted $%
\mathbf{C}_{w}^{(L)}$, and independence between $\{\mathbf{w}_{k}^{(L)}\}$
and $\{\mathbf{w}_{k}^{(N)}\}$ is assumed for simplicity). From Eq. (\ref%
{eq:XL1_update}) it is easily inferred that the random vector 
\begin{equation}
	\mathbf{z}_{k}^{(N)}\triangleq \mathbf{x}_{k+1}^{(L)}-\mathbf{A}%
	_{k}^{(L)}\left( \mathbf{x}_{k}^{(N)}\right) \,\mathbf{x}_{k}^{(L)}\text{,}
	\label{z_N_l}
\end{equation}%
must equal the sum 
\begin{equation}
	\mathbf{f}_{k}^{(L)}\left( \mathbf{x}_{k}^{(N)}\right) +\mathbf{w}_{k}^{(L)},
	\label{z_N_l_bis}
\end{equation}%
that depends on $\mathbf{x}_{k}^{(N)}$ \emph{only}; for this reason, $%
\mathbf{z}_{k}^{(N)}$ (\ref{z_N_l}) can be interpreted as a
pseudo-measurement about $\mathbf{x}_{k}^{(N)}$. In this case, the
pseudo-measurement information is conveyed by the message $m_{\mathrm{pm}}(%
\mathbf{x}_{k}^{(N)})$ ($PM2$) that expresses the\emph{\ correlation}
between the pdf of the random vector $\mathbf{z}_{k}^{(N)}$ (\ref{z_N_l})
(computed on the basis of the statistical information about the linear state
component made available by BIF$_{1}$) and the pdf obtained for $\mathbf{z}%
_{k}^{(N)}$ under the assumption that this vector is expressed by Eq. (\ref%
{z_N_l_bis}). The message $m_{\mathrm{pm}}(\mathbf{x}_{k}^{(N)})$ is
evaluated for each of the particles representing $\mathbf{x}_{k}^{(N)}$ in
BIF$_{2}$; this results in a set of $N_{p}$ particle weights employed in the
first measurement update of BIF$_{2}$ and different from those computed on
the basis of $\mathbf{y}_{k}$ (\ref{eq:y_t}) in its second measurement
update.

A graphical model similar to the one shown in Fig. \ref{Fig_3} can be easily
derived from the general model appearing in Fig. \ref{Fig_2} for\ case 
\textbf{C.2} too. The relevant differences with respect to case \textbf{C.1}
can be summarized as follows:

1) The backward filters BIF$_{1}$ and BIF$_{2}$ estimate $\mathbf{x}%
_{k}^{(1)}=\mathbf{x}_{k}^{(L)}$ and $\mathbf{x}_{k}^{(2)}=\mathbf{x}%
_{k}^{(N)}$, respectively; consequently, their nuisance substates are $%
\mathbf{\bar{x}}_{k}^{(1)}=\mathbf{x}_{k}^{(N)}$ and $\mathbf{\bar{x}}%
_{k}^{(2)}=\mathbf{x}_{k}^{(L)}$, respectively.

2) The BIF$_{2}{\rightarrow }$BIF$_{1}$ block is fed by the backward
predicted/smoothed pdfs computed by BIF$_{2}$; such pdfs are employed for:
a) generating the messages $m_{\mathrm{mg}1}(\mathbf{x}_{k}^{(N)})$ ($MG11$)
and $m_{\mathrm{mg}2}(\mathbf{x}_{k}^{(N)})$ ($MG21$) required to integrate
out the dependence of the Markov model (see Eq. (\ref{eq:XL1_update}))%
\begin{eqnarray}
	&&f\left( \mathbf{x}_{k+1}^{(L)}\left\vert \mathbf{x}_{k}^{(N)},\mathbf{x}%
	_{k}^{(L)}\right. \right)  \notag \\
	&=&\mathcal{N}\left( \mathbf{x}_{k}^{(L)};\mathbf{A}_{k}^{(L)}\left( \mathbf{%
		x}_{k}^{(N)}\right) \mathbf{x}_{k}^{(L)}+\mathbf{f}_{k}^{(L)}\left( \mathbf{x%
	}_{k}^{(N)}\right) ,\mathbf{C}_{w}^{(L)}\right)  \notag \\
	&&  \label{f_x_L}
\end{eqnarray}%
and of the measurement model $f(\mathbf{y}_{k}|\mathbf{x}_{k}^{(N)},\mathbf{x%
}_{k}^{(L)})$ (\ref{f_y_N}), respectively, on $\mathbf{x}_{k}^{(N)}$;\ b)
generating pseudo-measurement information about the substate $\mathbf{x}%
_{k}^{(L)}$ \emph{only}. As far as point a) is concerned, it is also
important to point out that the model $f(\mathbf{y}_{k}|\mathbf{x}_{k}^{(N)},%
\mathbf{x}_{k}^{(L)})$ ($f(\mathbf{x}_{k+1}^{(L)}|\mathbf{x}_{k}^{(N)},%
\mathbf{x}_{k}^{(L)})$) on which BIF$_{1}$ is based can be derived from Eq. (%
\ref{f_y_N}) (Eq. (\ref{f_x_L})) after setting $\mathbf{x}_{k}^{(N)}=\mathbf{%
	x}_{\mathrm{fp},k}^{(N)}$ ($\mathbf{x}_{k}^{(N)}=\mathbf{x}_{\mathrm{fe}%
	,k}^{(N)}$); here, $\mathbf{x}_{\mathrm{fp},k}^{(N)}$ ($\mathbf{x}_{\mathrm{%
		fe},k}^{(N)}$) denotes the prediction (the estimate) of $\mathbf{x}%
_{k}^{(N)} $ evaluated by the filter F$_{2}$ in the forward pass (further
details about this can be found in ref. \cite[Par. III-A]%
{Vitetta_DiViesti_2019})

The derivation of specific DBS algorithms based on the graphical model
illustrated in Fig. \ref{Fig_3} requires defining the scheduling of the
messages passed on it and deriving mathematical expressions for such
messages. These issues are investigated in detail in the following paragraph.

\subsection{Message Scheduling and Computation}

In this paragraph, the scheduling of a new \emph{recursive smoothing
	algorithm}, called \emph{double Bayesian smoothing algorithm} (DBSA) and
based on the graphical model illustrated in Fig. \ref{Fig_3}, and a
simplified version of it are described. Moreover, the expression of the
messages computed by the DBSA are illustrated.

\begin{figure}[tbp]
	\centering
	\includegraphics[width=0.6\textwidth]{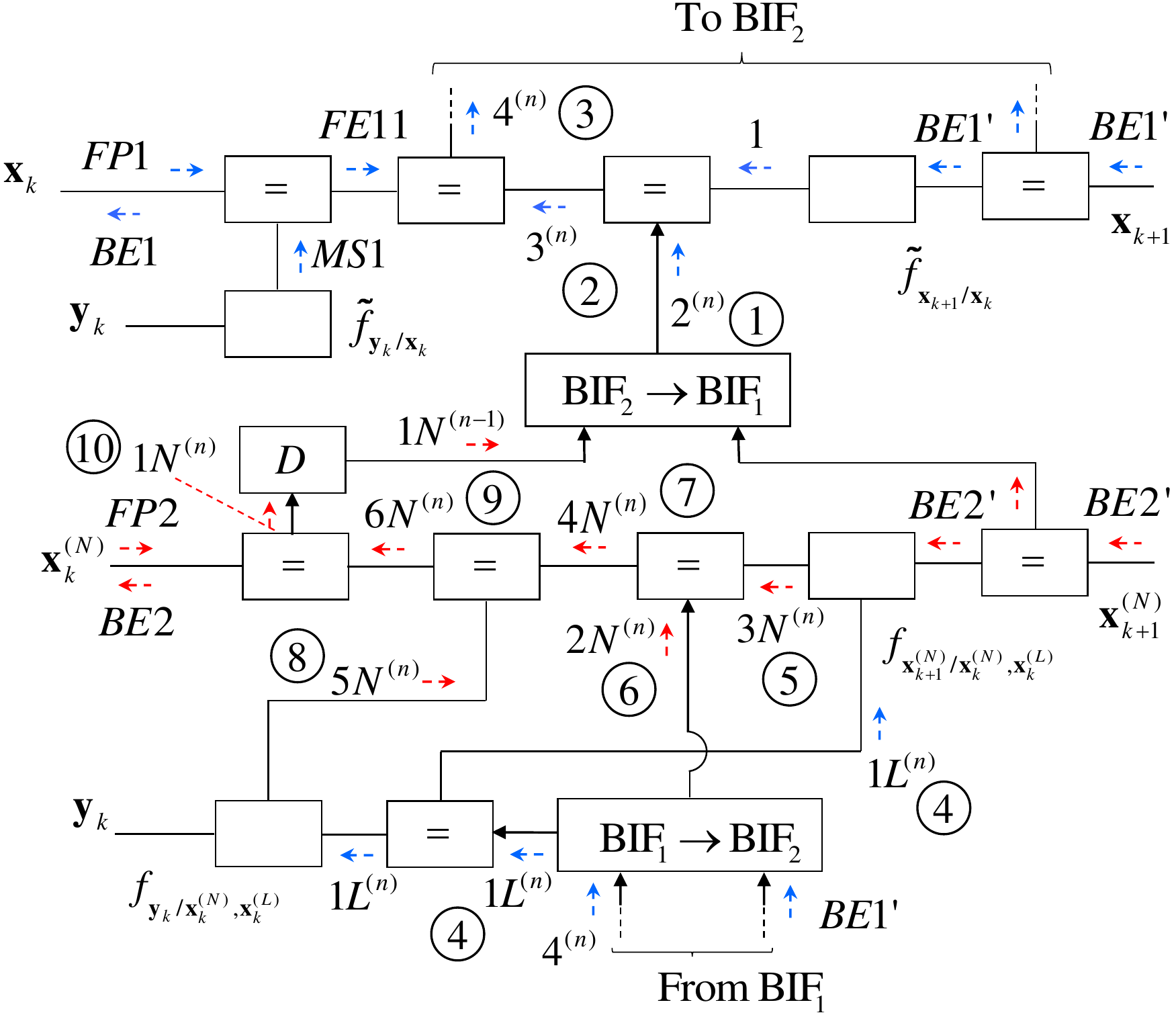}
	\caption{Representation of the message scheduling accomplished within the $%
		(T-k)-$th recursion of the backward pass of the DBSA; the circled integers $%
		1-10$ specify the order according to which the passed messages are computed
		in the $n-$th iteration embedded in the considered recursion. Blue and red
		arrows are employed to identify Gaussian messages and messages in particle
		form, respectively.}
	\label{Fig_4}
\end{figure}

The scheduling adopted in the DBSA mimics the one employed in ref. \cite%
{Vitetta_2018} (which, in turn, has been inspired by \cite{Fong_2002} and 
\cite{Lindsten_2016}). Moreover, in devising it, the presence of cycles in
the underlying graphical model has been accounted for by allowing multiple
passes of messages over the edges which such cycles consist of; this
explains why an iterative procedure is embedded in each recursion of the
DBSA. Our description of the devised scheduling is based on Fig. \ref{Fig_4}%
, that refers to the $(T-k)-$th recursion of the backward pass of the DBSA
(with $k=T-1$, $T-2$, $...$, $1$) and to the $n-$th iteration accomplished
within this recursion (with $n=1$, $2$, $...$, $n_{i}$, where $n_{i}$
represents the overall number of iterations). Note that, in this figure, a
simpler notation is adopted for most of the considered messages to ease
reading; in particular, the symbols $q$ , $q^{(n)}$, $qL$, $qL^{(n)}$, $qN$
and $qN^{(n)}$ are employed to represent the messages $m_{q}(\mathbf{x}_{k})$%
, $m_{q}^{(n)}(\mathbf{x}_{k})$, $m_{q}(\mathbf{x}_{k}^{(L)})$, $m_{q}^{(n)}(%
\mathbf{x}_{k}^{(L)})$, $m_{q}(\mathbf{x}_{k}^{(N)})$, $m_{q}^{(n)}(\mathbf{x%
}_{k}^{(N)})$, respectively (independently of the presence of an arrow and
of its direction in the considered message), and the presence of the
superscript $(n)$ in a given message means that such a message is computed
in the $n$-th iteration. Moreover, each of the passed messages conveys a
Gaussian pdf or a pdf in particle form. In the first case, the pdf of a
state/substate $\mathbf{x}$ is conveyed by the message%
\begin{equation}
	m_{G}\left( \mathbf{x}\right) =\mathcal{N(}\mathbf{x};\mathbf{\eta },\mathbf{%
		C}),  \label{Gaussian_message}
\end{equation}%
where $\mathbf{\eta }$ and $\mathbf{C}$ denote the mean and the covariance
of $\mathbf{x}$, respectively. In the second case, instead, its pdf is
conveyed by the message%
\begin{equation}
	m_{P}\left( \mathbf{x}\right) =\sum_{j=1}^{N_{p}}{m}_{P,j}\left( \mathbf{x}%
	\right) ,  \label{particle message}
\end{equation}%
where%
\begin{equation}
	m_{P,j}\left( \mathbf{x}\right) \triangleq w_{j}\,\delta \left( \mathbf{x}-%
	\mathbf{x}_{j}\right)  \label{particle message_j}
\end{equation}%
is its $j-$th \emph{component}; this represents the contribution of the $j-$%
th particle $\mathbf{x}_{j}$ and its weight $w_{j}$ to ${m}_{P}(\mathbf{x})$
(\ref{particle message}). The nature of each message can be easily inferred
from Fig. \ref{Fig_4}, since Gaussian messages and messages in particle form
are identified by blue and red arrows, respectively.

Before analysing the adopted scheduling, we need to define the \emph{input}
messages feeding the considered recursion of the DBSA and the outputs that
such a recursion produces. In the considered recursion, the DBSA \emph{input}
messages originate from:

1) The $k$-th recursion of the forward pass. These messages have been
generated by the DBF technique paired with the considered BIF scheme and, in
particular, by the DBF algorithm derived in ref. \cite[Par. III-B]%
{Vitetta_DiViesti_2019}, and have been stored (so that they are made
available to the backward pass).

2) The previous recursion (i.e., the $(T-k-1)-$th recursion) of the DBSA
itself.

As far as the input messages computed in the forward pass are concerned, BIF$%
_{1}$ is fed by the Gaussian messages (see Fig. \ref{Fig_4}) 
\begin{equation}
	\vec{m}_{\mathrm{fp}}\left( \mathbf{x}_{k}\right) \triangleq \mathcal{N}%
	\left( \mathbf{x}_{k};\mathbf{\eta }_{\mathrm{fp},k},\mathbf{C}_{\mathrm{fp}%
		,k}\right) .  \label{m_fp_L_MPF}
\end{equation}%
and%
\begin{equation}
	\vec{m}_{\mathrm{fe}1}\left( \mathbf{x}_{k}\right) \triangleq \mathcal{N}%
	\left( \mathbf{x}_{k};\mathbf{\eta }_{\mathrm{fe}1,k},\mathbf{C}_{\mathrm{fe}%
		1,k}\right) ,  \label{m_fe_L_MPF}
\end{equation}%
that convey the predicted pdf and the first filtered pdf, respectively,
computed by F$_{1}$ (in its $(k-1)-$th and in its $k-$th recursion,
respectively). The covariance matrix $\mathbf{C}_{\mathrm{fe}1,k}$ and the
mean vector $\mathbf{\eta }_{\mathrm{fe}1,k}$ are evaluated on the basis of
the associated precision matrix (see \cite[Eqs. (14)-(17)]{Vitetta_2018}) 
\begin{equation}
	\mathbf{W}_{\mathrm{fe}1,k}=\mathbf{H}_{k}\,\mathbf{W}_{e}\,\mathbf{H}%
	_{k}^{T}+\mathbf{W}_{\mathrm{fp},k}  \label{W_fe_L}
\end{equation}%
and of the associated transformed mean vector%
\begin{equation}
	\mathbf{w}_{\mathrm{fe}1,k}=\mathbf{H}_{k}\,\mathbf{W}_{e}\left( \mathbf{y}%
	_{k}-\mathbf{v}_{k}\right) +\mathbf{w}_{\mathrm{fp},k},  \label{w_fe_L}
\end{equation}%
respectively; here, $\mathbf{W}_{e}\triangleq \mathbf{C}_{e}^{-1}$, $\mathbf{%
	W}_{\mathrm{fp},k}\triangleq (\mathbf{C}_{\mathrm{fp},k})^{-1}$ and $\mathbf{%
	w}_{\mathrm{fp},k}\triangleq \mathbf{W}_{\mathrm{fp},k}\,\mathbf{\eta }_{%
	\mathrm{fp},k}$.

The backward filter BIF$_{2}$, instead, is fed by the forward messages $\vec{%
	m}_{\mathrm{fp}}(\mathbf{x}_{k}^{(N)})$ and $\vec{m}_{\mathrm{fe}1}(\mathbf{x%
}_{k}^{(N)})$, both in particle form (see Fig. \ref{Fig_4}); their $j-$th
components are represented by 
\begin{equation}
	\vec{m}_{\mathrm{fp},j}\left( \mathbf{x}_{k}^{(N)}\right) \triangleq
	w_{p}\,\delta \left( \mathbf{x}_{k}^{(N)}-\mathbf{x}_{k,j}^{(N)}\right) ,
	\label{m_fp_N_MPF}
\end{equation}%
and%
\begin{equation}
	\vec{m}_{\mathrm{fe}1,j}\left( \mathbf{x}_{k}^{(N)}\right) \triangleq w_{%
		\mathrm{fe},k,j}\,\delta \left( \mathbf{x}_{k}^{(N)}-\mathbf{x}%
	_{k,j}^{(N)}\right) ,  \label{m_fe_N_MPF}
\end{equation}%
respectively, with $j=1,2,...,N_{p}$; here, $\mathbf{x}_{k,j}^{(N)}$ is the $%
j-$th particle predicted by F$_{2}$ in the $(k-1)$-th recursion of the
forward pass (i.e., the $j-$th element of the particle set $S_{k}\triangleq
\{\mathbf{x}_{k,1}^{(N)}$, $\mathbf{x}_{k,2}^{(N)}$, $...$, $\mathbf{x}%
_{k,N_{p}}^{(N)}\}$), whereas $w_{p}\triangleq 1/N_{p}$ and $w_{\mathrm{fe}%
	,k,j}$ represent the (normalised) weights assigned to this particle in\ the
messages $\vec{m}_{\mathrm{fp}}(\mathbf{x}_{k}^{(N)})$ and $\vec{m}_{\mathrm{%
		fe}1}(\mathbf{x}_{k}^{(N)})$, respectively.

On the other hand, the input messages originating from the previous
recursion of the backward pass are the backward filtered Gaussian pdf 
\begin{equation}
	\cev{m}_{\mathrm{be}}\left( \mathbf{x}_{k+1}\right) \triangleq \mathcal{N}%
	\left( \mathbf{x}_{k+1};\mathbf{\eta }_{\mathrm{be},k+1},\mathbf{C}_{\mathrm{%
			be},k+1}\right)  \label{mess_be_l}
\end{equation}%
and the backward pdf 
\begin{equation}
	\cev{m}_{\mathrm{be}}\left( \mathbf{x}_{k+1}^{(N)}\right) \triangleq \delta
	\left( \mathbf{x}_{k+1}^{(N)}-\mathbf{x}_{\mathrm{be},k+1}^{(N)}\right) ,
	\label{mess_be_N_l}
\end{equation}%
that represents $\mathbf{x}_{k+1}^{(N)}$ through a single particle having
unit weight; these are computed by BIF$_{1}$ and BIF$_{2}$, respectively.
Consequently, in the considered recursion of the backward pass, all the
forward/backward input messages described above are processed to compute: 1)
the new backward pdfs $\cev{m}_{\mathrm{be}}(\mathbf{x}_{k})$ and $\cev{m}_{%
	\mathrm{be}}(\mathbf{x}_{k}^{(N)})$; 2) the smoothed statistical information
about $\mathbf{x}_{k}$ ($\mathbf{x}_{k}^{(N)}$) by properly merging forward
and backward messages generated by F$_{1}$ and BIF$_{1}$ (F$_{2}$ and BIF$%
_{2}$). As far as the last point is concerned, the evaluation of smoothed
information is based on the same conceptual approach as \cite%
{Fong_2002,Lindsten_2016,Vitetta_2018}. In fact, in our work, the \emph{joint%
} smoothing pdf $f(\mathbf{x}_{1:T}|\mathbf{y}_{1:T})$ is estimated by
providing multiple (say, $M$) \emph{realizations} of it. A single
realization (i.e., a single \emph{smoothed} state trajectory) is computed in
each backward pass; consequently, generating the whole output of the DBSA
requires running a single forward pass and $M$ distinct backward passes.
Moreover, the evaluation of the smoothed information is based on the
factorisation (\ref{factorisation3a}) or (\ref{factorisation3b}). In fact,
these formulas are exploited to merge the statistical information emerging
from the forward pass with that computed in any of the $M$ backward passes.

The message passing on which the DBSA is based can be divided in the three
consecutive phases listed below.

\textbf{Phase I} - In the first phase, the backward predicted pdf $\cev{m}%
_{1}(\mathbf{x}_{k})$ ($1$) is computed on the basis of the backward
filtered pdf $\cev{m}_{\mathrm{be}}(\mathbf{x}_{k+1})$ ($BE1^{\prime }$).

\textbf{Phase} \textbf{II} - In this phase, an iterative procedure for
computing and progressively refining the first backward filtered and the
smoothed pdfs of the whole state (BIF$_{1}$), and the second filtered and
the smoothed pdfs of the nonlinear state component (BIF$_{2}$) is carried
out. More specifically, in the $n$-th iteration of this procedure (with $n=1$%
, $2$, $...$, $n_{i}$), the ordered computation of the following messages is
accomplished in eight consecutive steps (see Fig. \ref{Fig_4}): 1) $%
m_{2}^{(n)}(\mathbf{x}_{k})$ ($2^{(n)}$; pdf conveying pseudo-measurement
information about $\mathbf{x}_{k}$); 2) $\cev{m}_{3}^{(n)}(\mathbf{x}_{k})$ (%
$3^{(n)}$; \emph{first backward filtered pdf of} $\mathbf{x}_{k}$); 3) $%
m_{4}^{(n)}(\mathbf{x}_{k})$ ($4^{(n)}$; \emph{smoothed pdf of} $\mathbf{x}%
_{k}$); 4) $m_{1}^{(n)}(\mathbf{x}_{k}^{(L)})$ ($1L^{(n)}$; pdf for
integrating out the dependence of $f(\mathbf{x}_{k+1}^{(N)}|\mathbf{x}%
_{k}^{(N)},\mathbf{x}_{k}^{(L)})$ and $f(\mathbf{y}_{k}|\mathbf{x}_{k}^{(N)},%
\mathbf{x}_{k}^{(L)})$ on $\mathbf{x}_{k}^{(L)}$), $\cev{m}_{3}^{(n)}(%
\mathbf{x}_{k}^{(N)})$ ($3N^{(n)}$; \emph{backward predicted pdf of} $%
\mathbf{x}_{k}^{(N)}$); 5) $m_{2}^{(n)}(\mathbf{x}_{k}^{(N)})$ ($2N^{(n)}$;
pdf conveying pseudo-measurement information about $\mathbf{x}_{k}^{(N)}$);
6) $m_{4}^{(n)}(\mathbf{x}_{k}^{(N)})$ ($4N^{(n)}$; \emph{first backward
	filtered pdf of} $\mathbf{x}_{k}^{(N)}$); 7) $m_{5}^{(n)}(\mathbf{x}%
_{k}^{(N)})$ ($5N^{(n)}$; message conveying measurement-based information
about $\mathbf{x}_{k}^{(N)}$); 8) $m_{6}^{(n)}(\mathbf{x}_{k}^{(N)})$ ($%
6N^{(n)}$; \emph{second backward filtered pdf of} $\mathbf{x}_{k}^{(N)}$), $%
m_{1}^{(n)}(\mathbf{x}_{k}^{(N)})$ ($1N^{(n)}$; \emph{smoothed pdf of }$%
\mathbf{x}_{k}^{(N)}$). Note that the message $m_{1}^{(n)}(\mathbf{x}%
_{k}^{(N)})$ computed in the last step of the $n$-th iteration is stored in
a memory cell (identified by the label `D'), so that it becomes available at
the beginning of the next iteration.

\textbf{Phase} \textbf{III} - In the third phase, the final smoothed pdf $%
m_{1}^{(n_{i})}(\mathbf{x}_{k}^{(N)})$ is exploited to compute: a) the final
backward pdf (i.e., the output message of BIF$_{2}$) $\cev{m}_{\mathrm{be}}(%
\mathbf{x}_{k}^{(N)})$; b) the new pseudo-measurement message $%
m_{2}^{(n_{i}+1)}(\mathbf{x}_{k})$, the final backward filtered pdf $\cev{m}%
_{3}^{(n_{i}+1)}(\mathbf{x}_{k})$, the final smoothed pdf $m_{4}^{(n_{i}+1)}(%
\mathbf{x}_{k})$ of $\mathbf{x}_{k} $ and, finally, the final backward pdf
(i.e., the output message of BIF$_{1}$) $m_{\mathrm{be}}(\mathbf{x}_{k})$.

\subsection{Message Computation}

The expressions of all the messages evaluated by the DBSA, with the
exception of the messages emerging from the BIF$_{1}{\rightarrow }$BIF$_{2}$
block and the BIF$_{2}{\rightarrow }$BIF$_{1}$ block, can be easily derived
by applying\ the few mathematical rules listed in Tables I-III of ref. \cite[%
App. A]{Vitetta_2019}; all such rules result from the application of the SPA
to equality nodes or nodes representing Gaussian functions. The derivation
of the algorithms for computing the pseudo-measurement messages $m_{2}^{(n)}(%
\mathbf{x}_{k})$ ($2^{(n)}$) and $m_{2}^{(n)}(\mathbf{x}_{k}^{(N)})$ ($%
2N^{(n)}$) emerging from the BIF$_{2}{\rightarrow }$BIF$_{1}$ block and the
BIF$_{1}{\rightarrow }$BIF$_{2}$ block, respectively, is based on the same
approach illustrated in refs. \cite[Par. IV-B]{Vitetta_2018} and \cite[%
Sects. IV-V]{Vitetta_2019}. On the other hand, the message $m_{1}^{(n)}(%
\mathbf{x}_{k}^{(L)})$ ($1L^{(n)}$)\ originating from the BIF$_{1}{%
	\rightarrow }$BIF$_{2}$ block results from marginalizing $m_{4}^{(n)}(%
\mathbf{x}_{k})$ ($4^{(n)}$) with respect to $\mathbf{x}_{k}^{(N)}$.

In the remaining part of this paragraph, the expressions of all the messages
computed in each of the three phases described above are provided for the ($%
T-k)-$th recursion of the backward pass; the derivation of these expressions
is sketched in Appendix \ref{app:A}.

\textbf{Phase I }- The computation of the \emph{backward predicted pdf} 
\begin{equation}
	\cev{m}_{1}\left( \mathbf{x}_{k}\right) \triangleq \mathcal{N}\left( \mathbf{%
		x}_{k};\mathbf{\eta }_{1,k},\mathbf{C}_{1,k}\right)  \label{m_bp}
\end{equation}%
of $\mathbf{x}_{k}$ involves $\cev{m}_{\mathrm{be}}(\mathbf{x}_{k+1})$ (\ref%
{mess_be_l}) and the pdf $\tilde{f}(\mathbf{x}_{k+1}|\mathbf{x}_{k})$ (\ref%
{Markov_ekf_approx}). Its parameters $\mathbf{\eta }_{1,k}$ and $\mathbf{C}%
_{1,k}$ are evaluated on the basis of the associated precision matrix%
\begin{equation}
	\mathbf{W}_{1,k}\triangleq \left( \mathbf{C}_{1,k}\right) ^{-1}=\mathbf{F}%
	_{k}^{T}\,\mathbf{P}_{k+1}\,\mathbf{W}_{\mathrm{be},k+1}\mathbf{F}_{k}
	\label{W_bp_x_l}
\end{equation}%
and of the associated transformed mean vector%
\begin{eqnarray}
	\mathbf{w}_{1,k} &\triangleq &\mathbf{W}_{1,k}\mathbf{\eta }_{1,k} =\mathbf{F}_{k}^{T}[\mathbf{P}_{k+1}\,\mathbf{w}_{\mathrm{be},k+1}-\mathbf{%
		W}_{\mathrm{be},k+1}\,\mathbf{Q}_{k+1}\,\mathbf{W}_{w}\,\mathbf{u}_{k}], 
	\notag \\
	&&  \label{w_bp_x_l}
\end{eqnarray}%
respectively; here, $\mathbf{W}_{\mathrm{be},k+1}\triangleq (\mathbf{C}_{%
	\mathrm{be},k+1})^{-1}$, $\mathbf{P}_{k+1}\triangleq \mathbf{\mathbf{I}}_{D}-%
\mathbf{W}_{\mathrm{be},k+1}\,\mathbf{Q}_{k+1}$, $\mathbf{Q}_{k+1}\triangleq
(\mathbf{W}_{w}+\mathbf{W}_{\mathrm{be},k+1})^{-1}$, $\mathbf{W}%
_{w}\triangleq (\mathbf{C}_{w})^{-1}$ and $\mathbf{w}_{\mathrm{be}%
	,k+1}\triangleq \mathbf{W}_{\mathrm{be},k+1}\,\mathbf{\eta }_{\mathrm{be}%
	,k+1}$.

\textbf{Phase} \textbf{II} - In the $n-$th iteration of this phase, the
eight consecutive steps listed below are carried out; for each step, all the
computed messages are described.

Step 1) - In this step, the message%
\begin{equation}
	m_{1}^{(n-1)}\left( \mathbf{x}_{k}^{(N)}\right)
	=\sum_{j=1}^{N_{p}}W_{1,k,j}^{(n-1)}\,\delta \left( \mathbf{x}_{k}^{(N)}-%
	\mathbf{x}_{k,j}^{(N)}\right) ,  \label{M_1_x_k_N}
\end{equation}%
computed in the previous iteration and conveying the smoothed pdf of $%
\mathbf{x}_{k}^{(N)}$ generated by F$_{2}$ and BIF$_{2}$ (see step 8)) is
processed jointly with $\cev{m}_{\mathrm{be}}(\mathbf{x}_{k+1}^{(N)})$ (\ref%
{mess_be_N_l}) in the BIF$_{2}{\rightarrow }$BIF$_{1}$ block to generate the
message%
\begin{equation}
	m_{2}^{(n)}\left( \mathbf{x}_{k}\right) =\mathcal{\mathcal{N}}\left( \mathbf{%
		x}_{k};\mathbf{\eta }_{2,k}^{(n)},\mathbf{C}_{2,k}^{(n)}\right) ,
	\label{m_pm_x_l}
\end{equation}%
that conveys the pseudo-measurement information provided to BIF$_{1}$. The
mean vector $\mathbf{\eta }_{2,k}^{(n)}$ and the covariance matrix $\mathbf{C%
}_{2,k}^{(n)}$ are evaluated as 
\begin{equation}
	\mathbf{\eta }_{2,k}^{(n)}=\left[ \left( \mathbf{\eta }_{L,k}^{(n)}\right)
	^{T},\left( \mathbf{\eta }_{N,k}^{(n)}\right) ^{T}\right] ^{T}
	\label{eta_pm_l_k}
\end{equation}%
and%
\begin{equation}
	\mathbf{C}_{2,k}^{(n)}=\left[ 
	\begin{array}{cc}
		\mathbf{C}_{LL,k}^{(n)} & \mathbf{C}_{LN,k}^{(n)} \\ 
		\left( \mathbf{C}_{LN,k}^{(n)}\right) ^{T} & \mathbf{C}_{NN,k}^{(n)}%
	\end{array}%
	\right] ,  \label{C_pm_l_k}
\end{equation}%
respectively, where%
\begin{equation}
	\mathbf{\eta }_{X,k}^{(n)}\triangleq \sum_{j=1}^{N_{p}}W_{1,k,j}^{(n-1)}\,%
	\mathbf{\eta }_{X,k,j}  \label{eta_pm_l_L_k}
\end{equation}%
is a $D_{X}$-dimensional mean vector (with $X=L$ and $N)$, 
\begin{equation}
	\mathbf{C}_{XY,k}^{(n)}\triangleq \sum_{j=1}^{N_{p}}W_{1,k,j}^{(n-1)}\,%
	\mathbf{r}_{XY,k,j}-\mathbf{\eta }_{X,k}\left( \mathbf{\eta }_{Y,k}\right)
	^{T}  \label{C_pm_l_L_k_bis}
\end{equation}%
is a $D_{X}\times D_{Y}$ covariance (or cross-covariance) matrix (with $%
XY=LL $, $NN$ and $LN)$, $\mathbf{\eta }_{L,k,j}=\mathbf{\tilde{\eta}}_{k,j}$%
, $\mathbf{\eta }_{N,k,j}=\mathbf{x}_{k,j}^{(N)}$, $\mathbf{r}%
_{LL,k,j}\triangleq \mathbf{\tilde{C}}_{k,j}+\mathbf{\tilde{\eta}}_{k,j}(%
\mathbf{\tilde{\eta}}_{k,j})^{T}$, $\mathbf{r}_{NN,k,j}\triangleq \mathbf{x}%
_{k,j}^{(N)}(\mathbf{x}_{k,j}^{(N)})^{T}$ and $\mathbf{r}_{LN,k,j}\triangleq 
\mathbf{\tilde{\eta}}_{k,j}(\mathbf{x}_{k,j}^{(N)})^{T}$. The covariance
matrix $\mathbf{\tilde{C}}_{k,j}$ and the mean vector $\mathbf{\tilde{\eta}}%
_{k,j}$ are computed on the basis of the associated precision matrix%
\begin{equation}
	\mathbf{\tilde{W}}_{k,j}\triangleq \left( \mathbf{\tilde{C}}_{k,j}\right)
	^{-1}=\left( \mathbf{A}_{k,j}^{(N)}\right) ^{T}\mathbf{W}_{w}^{(N)}\mathbf{A}%
	_{k,j}^{(N)}  \label{eq:W_pm_L_j}
\end{equation}%
and of the associated transformed mean vector%
\begin{equation}
	\mathbf{\tilde{w}}_{k,j}\triangleq \mathbf{\tilde{W}}_{k,j}\mathbf{\tilde{%
			\eta}}_{k,j}=\left( \mathbf{A}_{k,j}^{(N)}\right) ^{T}\mathbf{W}_{w}^{(N)}%
	\mathbf{z}_{k,j}^{(L)},  \label{eq:w_pm_L_j}
\end{equation}%
respectively; here, $\mathbf{A}_{k,j}^{(N)}\triangleq \mathbf{A}_{k}^{(N)}(%
\mathbf{x}_{k,j}^{(N)})$, 
\begin{equation}
	\mathbf{z}_{k,j}^{(L)}\triangleq \mathbf{x}_{\mathrm{be},k+1}^{(N)}-\mathbf{f%
	}_{k,j}^{(N)}  \label{PM_z_L}
\end{equation}%
is an iteration-independent pseudo-measurement (see Eq. (\ref{eq:z_L_l}))
and $\mathbf{f}_{k,j}^{(N)}\triangleq \mathbf{f}_{k}^{(N)}(\mathbf{x}%
_{k,j}^{(N)})$. Note that, in the first iteration, 
\begin{equation}
	W_{1,k,j}^{(n-1)}=W_{1,k,j}^{(0)}=w_{\mathrm{fe},k,j},
	\label{initial_weights}
\end{equation}%
for any $j$, i.e. $m_{1}^{(0)}(\mathbf{x}_{k}^{(N)})=\vec{m}_{\mathrm{fe}1}(%
\mathbf{x}_{k}^{(N)}$) (see Eq. (\ref{m_fe_N_MPF})) since the initial
information available about the particle set are those originating from the
forward pass. For this reason, the particles $\,\{\mathbf{x}_{k,j}^{(N)}\}$
and their weights $\{w_{\mathrm{fe},k,j}\}$ are stored in the memory cell at
the beginning of the first iteration.

Step 2) - In this step, the first backward filtered pdf $\cev{m}_{3}^{(n)}(%
\mathbf{x}_{k})$ of $\mathbf{x}_{k}$ is computed as (see Fig. \ref{Fig_4}) 
\begin{eqnarray}
	\cev{m}_{3}^{(n)}\left( \mathbf{x}_{k}\right) &=&\cev{m}_{1}(\mathbf{x}%
	_{k})\,m_{2}^{(n)}(\mathbf{x}_{k})  \label{m_be1_x_la} \\
	&=&\mathcal{\mathcal{N}}\left( \mathbf{x}_{k};\mathbf{\eta }_{3,k}^{(n)},%
	\mathbf{C}_{3,k}^{(n)}\right) ,  \label{m_be1_x_laa}
\end{eqnarray}%
where the messages $\cev{m}_{1}(\mathbf{x}_{k})$ and $m_{2}^{(n)}(\mathbf{x}%
_{k})$ are given by Eq. (\ref{m_bp}) and Eq. (\ref{m_pm_x_l}), respectively.
The covariance matrix $\mathbf{C}_{3,k}^{(n)}$ and the mean vector $\mathbf{%
	\eta }_{3,k}^{(n)}$ are computed on the basis of the associated precision
matrix 
\begin{equation}
	\mathbf{W}_{3,k}^{(n)}\triangleq (\mathbf{C}_{3,k}^{(n)})^{-1}=\mathbf{W}%
	_{1,k}+\mathbf{W}_{2,k}^{(n)}  \label{W_be1_l_kn}
\end{equation}%
and the associated transformed mean vector%
\begin{equation}
	\mathbf{w}_{3,k}^{(n)}\triangleq \mathbf{W}_{3,k}^{(n)}\,\mathbf{\eta }%
	_{3,k}^{(n)}\,=\mathbf{w}_{1,k}+\mathbf{w}_{2,k}^{(n)},  \label{w_be1_l_kn}
\end{equation}%
respectively; here, $\mathbf{W}_{2,k}^{(n)}\triangleq (\mathbf{C}%
_{2,k}^{(n)})^{-1}$, $\mathbf{w}_{2,k}^{(n)}\triangleq \mathbf{W}%
_{2,k}^{(n)}\,\mathbf{\eta }_{2,k}^{(n)}$, and $\mathbf{W}_{1,k}$ and $%
\mathbf{w}_{1,k}$ are given by Eqs. (\ref{W_bp_x_l}) and (\ref{w_bp_x_l}),
respectively. From Eqs. (\ref{W_be1_l_kn})-(\ref{w_be1_l_kn}) the
expressions 
\begin{equation}
	\mathbf{C}_{3,k}^{(n)}=\mathbf{W}_{k}^{(n)}\mathbf{C}_{2,k}^{(n)}
	\label{C_be1_l_ka}
\end{equation}%
and%
\begin{equation}
	\mathbf{\eta }_{3,k}^{(n)}=\mathbf{W}_{k}^{(n)}\left[ \mathbf{C}%
	_{2,k}^{(n)}\,\mathbf{w}_{1,k}+\mathbf{\eta }_{2,k}^{(n)}\right]
	\label{eta_be1_l_ka}
\end{equation}%
can be easily inferred; here, $\mathbf{W}_{k}^{(n)}\triangleq \lbrack 
\mathbf{C}_{2,k}^{(n)}\mathbf{W}_{1,k}+\mathbf{I}_{D}]^{-1}$.

Step 3) - In this step, the smoothed pdf $m_{4}^{(n)}(\mathbf{x}_{k})$ of $%
\mathbf{x}_{k}$ is evaluated as (see Fig. \ref{Fig_4}) 
\begin{eqnarray}
	m_{4}^{(n)}\left( \mathbf{x}_{k}\right) &=&\vec{m}_{\mathrm{fe}1}\left( 
	\mathbf{x}_{k}\right) \cev{m}_{3}^{(n)}\left( \mathbf{x}_{k}\right)
	\label{m_sm_x_la} \\
	&=&\mathcal{\mathcal{N}}\left( \mathbf{x}_{k};\mathbf{\eta }_{4,k}^{(n)},%
	\mathbf{C}_{4,k}^{(n)}\right) ,  \label{m_sm_x_l}
\end{eqnarray}%
where the messages $\vec{m}_{\mathrm{fe}1}\left( \mathbf{x}_{k}\right) $ and 
$\cev{m}_{3}^{(n)}\left( \mathbf{x}_{k}\right) $ are given by Eqs. (\ref%
{m_fe_L_MPF}) and (\ref{m_be1_x_laa}), respectively. The covariance matrix $%
\mathbf{C}_{4,k}^{(n)}$ and the mean vector $\mathbf{\eta }_{4,k}^{(n)}$ are
computed on the basis of the associated precision matrix%
\begin{equation}
	\mathbf{W}_{4,k}^{(n)}=\mathbf{W}_{\mathrm{fe}1,k}+\mathbf{W}_{3,k}^{(n)}
	\label{W_sm_l_k}
\end{equation}%
and of the associated transformed mean vector%
\begin{equation}
	\mathbf{w}_{4,k}^{(n)}=\mathbf{w}_{\mathrm{fe}1,k}+\mathbf{w}_{3,k}^{(n)},
	\label{w_sm_l_k}
\end{equation}%
respectively. Note that Eq. (\ref{m_sm_x_la}) represents an instance of Eq. (%
\ref{factorisation3b}), since $\vec{m}_{\mathrm{fe}1}\left( \mathbf{x}%
_{k}\right) $ and$\,\cev{m}_{3}^{(n)}\left( \mathbf{x}_{k}\right) $
correspond to $\vec{m}_{\mathrm{fe}1}(\mathbf{x}_{k}^{(i)})$ and $\cev{m}_{%
	\mathrm{be}1}(\mathbf{x}_{k}^{(i)})$, respectively ($\mathbf{x}_{k}^{(i)}=%
\mathbf{x}_{k}^{(N)}$ in this case).

Step 4) - In this step, the message 
\begin{equation}
	m_{1}^{(n)}\left( \mathbf{x}_{k}^{(L)}\right) \triangleq \int
	m_{4}^{(n)}\left( \mathbf{x}_{k}\right) \,d\mathbf{x}_{k}^{(N)}=\mathcal{N(}%
	\mathbf{x}_{k}^{(L)};\mathbf{\tilde{\eta}}_{1,k}^{(n)},\mathbf{\tilde{C}}%
	_{1,k}^{(n)}),  \label{m_fe_L_EKF_2}
\end{equation}%
is computed in the BIF$_{1}{\rightarrow }$BIF$_{2}$ block. In practice, the
mean $\mathbf{\tilde{\eta}}_{1,k}^{(n)}$ and the covariance matrix $\mathbf{%
	\tilde{C}}_{1,k}^{(n)}$ are extracted from the mean $\mathbf{\eta }%
_{4,k}^{(n)}$ and the covariance matrix $\mathbf{C}_{4,k}^{(n)}$ of $%
m_{4}^{(n)}(\mathbf{x}_{k})$ (\ref{m_sm_x_l}), respectively, since $\mathbf{x%
}_{k}^{(L)}$ consists of the first $D_{L}$ elements of $\mathbf{x}_{k}$.

Then, the backward predicted pdf $\cev{m}_{3}^{(n)}(\mathbf{x}_{k}^{(N)})$
is evaluated as (see Fig. \ref{Fig_4})%
\begin{eqnarray}
	\cev{m}_{3}^{(n)}\left( \mathbf{x}_{k}^{(N)}\right) &=&\int \int f(\mathbf{x}%
	_{k+1}^{(N)}|\mathbf{x}_{k}^{(N)},\mathbf{x}_{k}^{(L)})\,\cev{m}_{\mathrm{be}%
	}\left( \mathbf{x}_{k+1}^{(N)}\right)  \notag \\
	&&\cdot m_{1}^{(n)}\left( \mathbf{x}_{k}^{(L)}\right) \,d\mathbf{x}%
	_{k}^{(N)}d\mathbf{x}_{k+1}^{(N)}.  \label{weight_bpa}
\end{eqnarray}%
Actually, what is really needed in our computations is the value taken on by
this message (and also by messages $m_{2}^{(n)}(\mathbf{x}_{k}^{(N)})$ and $%
m_{5}^{(n)}(\mathbf{x}_{k}^{(N)})$ evaluated in step 5) and in step 7),
respectively) for $\mathbf{x}_{k}^{(N)}=\mathbf{x}_{k,j}^{(N)}$ (see Eqs. (%
\ref{m_be1_k}), (\ref{message_6}) and (\ref{m_1_prod})); such a value,
denoted $w_{3,k,j}^{(n)}$, is computed as 
\begin{equation}
	w_{3,k,j}^{(n)}=D_{3,k,j}^{(n)}\exp \left( -\frac{1}{2}Z_{3,k,j}^{(n)}%
	\right) ,  \label{weight_3}
\end{equation}%
where%
\begin{equation}
	D_{3,k,j}^{(n)}=(2\pi )^{-D_{N}/2}(\det (\mathbf{C}%
	_{3,k,j}^{(N)}[n]))^{-1/2},  \label{D_3}
\end{equation}%
\begin{equation}
	Z_{3,k,j}^{(n)}\triangleq \left\Vert \mathbf{x}_{\mathrm{be},k+1}^{(N)}-%
	\mathbf{\eta }_{3,k,j}^{(N)}[n]\right\Vert _{\mathbf{W}%
		_{3,k,j}^{(N)}[n]}^{2},  \label{Z_bp}
\end{equation}%
$\left\Vert \mathbf{x}\right\Vert _{\mathbf{W}}^{2}\triangleq \mathbf{x}^{T}%
\mathbf{Wx}$ denotes the square of the norm of the vector $\mathbf{x}$ with
respect to the positive definite matrix $\mathbf{W}$, 
\begin{equation}
	\mathbf{\eta }_{3,k,j}^{(N)}[n]\triangleq \mathbf{A}_{k,j}^{(N)}\,\mathbf{%
		\tilde{\eta}}_{1,k}^{(n)}+\mathbf{f}_{k,j}^{(N)},  \label{eta_bp}
\end{equation}%
$\mathbf{W}_{3,k,j}^{(N)}[n]\triangleq (\mathbf{C}_{3,k,j}^{(N)}[n])^{-1}$
and%
\begin{equation}
	\mathbf{C}_{3,k,j}^{(N)}[n]\triangleq \mathbf{A}_{k,j}^{(N)}\,\mathbf{\tilde{%
			C}}_{1,k}^{(n)}\,\left( \mathbf{A}_{k,j}^{(N)}\right) ^{T}+\mathbf{C}%
	_{w}^{(N)}.  \label{cov_bp}
\end{equation}

Step 5) - In this step, the message $m_{2}^{(n)}(\mathbf{x}_{k}^{(N)})$,
conveying pseudo-measurement information about the nonlinear state
component, is computed in the BIF$_{1}{\rightarrow }$BIF$_{2}$ block. The
value $w_{2,k,j}^{(n)}$ taken on by this message for $\mathbf{x}_{k}^{(N)}=%
\mathbf{x}_{k,j}^{(N)}$ is evaluated as 
\begin{equation}
	w_{2,k,j}^{(n)}=D_{2,k,j}^{(n)}\exp \left( -\frac{1}{2}Z_{2,k,j}^{(n)}\right)
	\label{m_pm_x_N_l_j}
\end{equation}%
for any $j$; here, 
\begin{equation}
	Z_{2,k,j}^{(n)}\triangleq \left\Vert \mathbf{\check{\eta}}%
	_{z,k,j}^{(n)}\right\Vert _{\mathbf{\check{W}}_{z,k,j}^{(n)}}^{2}-\left\Vert 
	\mathbf{f}_{k,j}^{(L)}\right\Vert _{\mathbf{W}_{w}^{(L)}}^{2}-\left\Vert 
	\mathbf{\check{\eta}}_{2,k,j}^{(n)}\right\Vert _{\mathbf{\check{W}}%
		_{2,k,j}^{(n)}}^{2},  \label{Z_pm}
\end{equation}%
$\mathbf{W}_{w}^{(L)}\triangleq \lbrack \mathbf{C}_{w}^{(L)}]^{-1}$, $%
\mathbf{f}_{k,j}^{(L)}\triangleq \mathbf{f}_{k}^{(L)}(\mathbf{x}%
_{k,j}^{(N)}) $, $\mathbf{\check{W}}_{z,k,j}^{(n)}\triangleq (\mathbf{\check{%
		C}}_{z,k,j}^{(n)})^{-1}$, $\mathbf{\check{w}}_{z,k,j}^{(n)}\triangleq 
\mathbf{\check{W}}_{z,k,j}^{(n)}\mathbf{\check{\eta}}_{z,k,j}^{(n)}$, 
\begin{equation}
	\mathbf{\check{\eta}}_{z,k,j}^{(n)}=\mathbf{\tilde{\eta}}_{\mathrm{be},k+1}-%
	\mathbf{A}_{k,j}^{(L)}\,\mathbf{\tilde{\eta}}_{1,k}^{(n)},
	\label{eta_mess_z_N}
\end{equation}%
\begin{equation}
	\mathbf{\check{C}}_{z,k,j}^{(n)}=\mathbf{\tilde{C}}_{\mathrm{be},k+1}-%
	\mathbf{A}_{k,j}^{(L)}\,\mathbf{\tilde{C}}_{1,k}^{(n)}\,\left( \mathbf{A}%
	_{k,j}^{(L)}\right) ^{T},  \label{C_mess_Z_N_bis}
\end{equation}%
\begin{equation}
	\mathbf{\check{W}}_{2,k,j}^{(n)}\triangleq \left( \mathbf{\check{C}}%
	_{2,k,j}^{(n)}\right) ^{-1}=\mathbf{\check{W}}_{z,k,j}^{(n)}+\mathbf{W}%
	_{w}^{(L)},  \label{W_pm_x_N_l_j}
\end{equation}%
\begin{equation}
	D_{2,k,j}^{(n)}\triangleq (2\pi )^{-D_{L}/2}[\det (\mathbf{\check{C}}%
	_{k,j}^{(n)})]^{-1/2}  \label{D_2}
\end{equation}%
$\mathbf{\check{C}}_{k,j}^{(n)}\triangleq \mathbf{\check{C}}_{z,k,j}^{(n)}+%
\mathbf{C}_{w}^{(L)}$, $\mathbf{A}_{k,j}^{(L)}\triangleq \mathbf{A}%
_{k}^{(L)}(\mathbf{x}_{k,j}^{(N)})$, $\mathbf{\check{\eta}}_{2,k,j}^{(n)}$
is evaluated on the basis of the\ associated transformed mean vector%
\begin{equation}
	\mathbf{\check{w}}_{2,k,j}^{(n)}\triangleq \mathbf{\check{W}}_{2,k,j}^{(n)}\,%
	\mathbf{\check{\eta}}_{2,k,j}^{(n)}=\mathbf{\check{w}}_{z,k,j}^{(n)}+\mathbf{%
		W}_{w}^{(L)}\,\mathbf{f}_{k,j}^{(L)},  \label{w_pm_x_N_l_j}
\end{equation}%
and the mean $\mathbf{\tilde{\eta}}_{\mathrm{be},k+1}$ and the covariance
matrix $\mathbf{\tilde{C}}_{\mathrm{be},k+1}$ are extracted from the mean $%
\mathbf{\eta }_{\mathrm{be},k+1}$ and the covariance matrix $\mathbf{C}_{%
	\mathrm{be},k+1}$ of $\cev{m}_{\mathrm{be}}(\mathbf{x}_{k+1})$ (\ref%
{mess_be_l}), since they refer to $\mathbf{x}_{k}^{(L)}$ only.

Step 6) - In this step, the message $m_{4}^{(n)}(\mathbf{x}_{k}^{(N)})$,
conveying the first backward filtered pdf of $\mathbf{x}_{k}^{(N)}$, is
computed as (see Fig. \ref{Fig_4}) 
\begin{equation}
	\cev{m}_{4}^{(n)}(\mathbf{x}_{k}^{(N)})=\cev{m}_{3}^{(n)}(\mathbf{x}%
	_{k}^{(N)})\,m_{2}^{(n)}(\mathbf{x}_{k}^{(N)}).  \label{m_be1_k}
\end{equation}%
The value~$w_{4,k,j}^{(n)}$ taken on by this message for $\mathbf{x}%
_{k}^{(N)}=\mathbf{x}_{k,j}^{(N)}$ is given by (see Eqs. (\ref{weight_3})
and (\ref{m_pm_x_N_l_j})) 
\begin{equation}
	w_{4,k,j}^{(n)}\triangleq w_{2,k,j}^{(n)}\,w_{3,k,j}^{(n)}\,
	\label{w_fe_2_x_N_l}
\end{equation}%
for any $j$.

Step 7) - In this step, the message conveying measurement-based information
about $\mathbf{x}_{k}^{(N)}$ is computed as (see Fig. \ref{Fig_4})%
\begin{eqnarray}
	m_{5}^{(n)}(\mathbf{x}_{k}^{(N)})\, &=&\int f(\mathbf{y}_{k}|\mathbf{x}%
	_{k}^{(N)},\,\mathbf{x}_{k}^{(L)})\,m_{1}^{(n)}(\mathbf{x}_{k}^{(L)})\,d%
	\mathbf{x}_{k}^{(L)}  \label{eq:mess_ms_PF} \\
	&=&\mathcal{N}\left( \mathbf{y}_{l};\mathbf{\bar{\eta}}_{5,k}^{(n)}\left( 
	\mathbf{x}_{k}^{(N)}\right) ,\mathbf{\bar{C}}_{5,k}^{(n)}\left( \mathbf{x}%
	_{k}^{(N)}\right) \right)  \notag \\
	&&  \label{eq:weight_before_resampling}
\end{eqnarray}%
where%
\begin{equation}
	\mathbf{\bar{\eta}}_{5,k}^{(n)}\left( \mathbf{x}_{k}^{(N)}\right) \triangleq 
	\mathbf{B}_{k}\left( \mathbf{x}_{k}^{(N)}\right) \,\mathbf{\tilde{\eta}}%
	_{1,k}^{(n)}+\mathbf{g}_{k}\left( \mathbf{x}_{k}^{(N)}\right)  \label{eta_5}
\end{equation}%
and%
\begin{equation}
	\mathbf{\bar{C}}_{5,k}^{(n)}\left( \mathbf{x}_{k}^{(N)}\right) \triangleq 
	\mathbf{B}_{k}\left( \mathbf{x}_{k}^{(N)}\right) \,\mathbf{\tilde{C}}%
	_{1,k}^{(n)}\,\,\mathbf{B}_{k}^{T}\left( \mathbf{x}_{k}^{(N)}\right) +%
	\mathbf{C}_{e}.  \label{C_5}
\end{equation}%
Consequently, the value taken on by $m_{5}^{(n)}(\mathbf{x}_{k}^{(N)})\,$\
for $\mathbf{x}_{k}^{(N)}=\mathbf{x}_{k,j}^{(N)}$ is 
\begin{eqnarray}
	w_{5,k,j}^{(n)}\, &=&\mathcal{N}\left( \mathbf{y}_{k};\mathbf{\bar{\eta}}%
	_{5,k,j}^{(n)},\mathbf{\bar{C}}_{5,k,j}^{(n)}\right)  \label{weight_5_a} \\
	&=&D_{5,k,j}^{(n)}\exp \left( -\frac{1}{2}Z_{5,k,j}^{(n)}\right) ,
	\label{weight_5_b}
\end{eqnarray}%
where%
\begin{equation}
	\mathbf{\bar{\eta}}_{5,k,j}^{(n)}\left( \mathbf{x}_{k}^{(N)}\right)
	\triangleq \mathbf{\bar{\eta}}_{5,k}^{(n)}\left( \mathbf{x}%
	_{k,j}^{(N)}\right) =\mathbf{B}_{k,j}\,\mathbf{\tilde{\eta}}_{1,k}^{(n)}+%
	\mathbf{g}_{k,j},  \label{eq:ms_fe_PF}
\end{equation}%
\begin{equation}
	\mathbf{\bar{C}}_{5,k,j}^{(n)}\triangleq \mathbf{\bar{C}}_{5,k}^{(n)}\left( 
	\mathbf{x}_{k,j}^{(N)}\right) =\mathbf{B}_{k,j}\,\mathbf{\tilde{C}}%
	_{1,k}^{(n)}\,\mathbf{B}_{k,j}^{T}+\mathbf{C}_{e},  \label{eq:C_sm_PF}
\end{equation}%
$\mathbf{B}_{k,j}\triangleq \mathbf{B}_{k}(\mathbf{x}_{k,j}^{(N)})$, $%
\mathbf{g}_{k,j}\triangleq \mathbf{g}_{l}(\mathbf{x}_{k,j}^{(N)})$, 
\begin{equation}
	D_{5,k,j}^{(n)}\triangleq (2\pi )^{-P/2}[\det (\mathbf{\bar{C}}%
	_{5,k,j}^{(n)})]^{-1/2}  \label{D_5}
\end{equation}%
\begin{equation}
	Z_{5,k,j}^{(n)}\triangleq \left\Vert \mathbf{y}_{k}-\mathbf{\bar{\eta}}%
	_{5,k,j}^{(n)}\right\Vert _{\mathbf{\bar{W}}_{5,k,j}^{(n)}}^{2}  \label{Z_5}
\end{equation}%
and $\mathbf{\bar{W}}_{5,k,j}^{(n)}\triangleq (\mathbf{\bar{C}}%
_{5,k,j}^{(n)})^{-1}$. Then, the message $\cev{m}_{6}^{(n)}(\mathbf{x}%
_{k}^{(N)})$ is evaluated as (see Fig. \ref{Fig_4})%
\begin{equation}
	\cev{m}_{6}^{(n)}\left( \mathbf{x}_{k}^{(N)}\right) =\cev{m}_{4}^{(n)}\left( 
	\mathbf{x}_{k}^{(N)}\right) \,m_{5}^{(n)}\left( \mathbf{x}_{k}^{(N)}\right)
	\,.  \label{message_6}
\end{equation}%
Its value for $\mathbf{x}_{k}^{(N)}=\mathbf{x}_{k,j}^{(N)}$ is given by (see
Eqs. (\ref{weight_3}), (\ref{m_pm_x_N_l_j}) and (\ref{weight_5_b}))%
\begin{eqnarray}
	w_{6,k,j}^{(n)}
	&=&w_{4,k,j}^{(n)}w_{5,k,j}^{(n)}\,=w_{2,k,j}^{(n)}%
	\,w_{3,k,j}^{(n)}w_{5,k,j}^{(n)}  \label{weight_6} \\
	&=&D_{6,k,j}^{(n)}\exp \left( -\frac{1}{2}Z_{6,k,j}^{(n)}\right)
	\label{weight_6b}
\end{eqnarray}%
where%
\begin{equation}
	D_{6,k,j}^{(n)}\triangleq D_{2,k,j}^{(n)}\,D_{3,k,j}^{(n)}\,D_{5,k,j}^{(n)}
	\label{D_6}
\end{equation}%
and 
\begin{equation}
	Z_{6,k,j}^{(n)}\triangleq Z_{2,k,j}^{(n)}+Z_{3,k,j}^{(n)}+Z_{5,k,j}^{(n)}.
	\label{Z_tot}
\end{equation}%
Note that the weight $w_{6,k,j}^{(n)}$ conveys the information provided by
the backward state transition ($w_{3,k,j}^{(n)}$), the pseudo-measurements ($%
w_{2,k,j}^{(n)}$) and the measurements ($w_{5,k,j}^{(n)}$).

Step 8) - In this step, the message $m_{1}^{(n)}(\mathbf{x}_{k}^{(N)})$,
conveying the smoothed pdf of $\mathbf{x}_{k}^{(N)}$ evaluated in the $n$-th
iteration, is computed as (see Fig. \ref{Fig_4}) 
\begin{equation}
	m_{1}^{(n)}\left( \mathbf{x}_{k}^{(N)}\right) =\vec{m}_{\mathrm{fp}}\left( 
	\mathbf{x}_{k}^{(N)}\right) \,\cev{m}_{6}^{(n)}\left( \mathbf{x}%
	_{k}^{(N)}\right) ;  \label{m_1_prod}
\end{equation}%
this formula represents an instance of Eq. (\ref{factorisation3a}), since $%
\vec{m}_{\mathrm{fp}1}(\mathbf{x}_{k}^{(N)})$ and$\,\cev{m}_{6}^{(n)}(%
\mathbf{x}_{k}^{(N)})$ correspond to $\vec{m}_{\mathrm{fp}}(\mathbf{x}%
_{k}^{(i)})$ and $\cev{m}_{\mathrm{be}2}(\mathbf{x}_{k}^{(i)})$,
respectively ($\mathbf{x}_{k}^{(i)}=\mathbf{x}_{k}^{(N)}$ in this case). \
The $j-$th component of $m_{1}^{(n)}(\mathbf{x}_{k}^{(N)})$ is evaluated as
(see Eqs. (\ref{m_fp_N_MPF}) and (\ref{weight_6}))%
\begin{eqnarray}
	m_{1,j}^{(n)}\left( \mathbf{x}_{k}^{(N)}\right) &=&\vec{m}_{\mathrm{fp}%
		,j}\left( \mathbf{x}_{k}^{(N)}\right) \,w_{6,k,j}^{(n)}  \label{m_sm_j_N1} \\
	&=&w_{1,k,j}^{(n)}\,\,\delta \left( \mathbf{x}_{k}^{(N)}-\mathbf{x}%
	_{k,j}^{(N)}\right) ,  \label{m_sm_j_N2}
\end{eqnarray}%
where 
\begin{equation}
	w_{1,k,j}^{(n)}\triangleq w_{p}\,w_{6,k,j}^{(n)}.  \label{w_sm_j_Na}
\end{equation}%
Then, the weights $\{w_{1,k,j}^{(n)}\}$ are normalized; the $j$-th
normalised weight is computed as%
\begin{equation}
	W_{1,k,j}^{(n)}\triangleq C_{k}^{(n)}\,w_{1,k,j}^{(n)}\,,
	\label{W_fe_2_x_N_l}
\end{equation}%
with $j=1,2,...,N_{p}$, where $C_{k}^{(n)}\triangleq
1/\sum\limits_{j=0}^{N_{p}-1}w_{1,k,j}^{(n)}$. Moreover, the weights $%
\{W_{1,k,j}^{(n)}\}$ are stored for the next iteration. This concludes the $%
n-$th iteration. Then, the index $n$ is increased by one, and a new
iteration is started by going back to step 1) if $n<n_{i}+1$; otherwise
(i.e., if $n=n_{i}+1$), we proceed with the next phase.

\textbf{Phase III} - In this phase, $\cev{m}_{\mathrm{be}}(\mathbf{x}%
_{k}^{(N)})$ (i.e., the BIF$_{2}$ output message) is computed first; then,
steps 1) and 2) of phase II are accomplished in order to compute all the
statistical information required for the evaluation of the backward estimate 
$\cev{m}_{\mathrm{be}}\left( \mathbf{x}_{k}\right) $ (i.e., the BIF$_{1}$
output message). More specifically, we first sample the set $S_{k}$ once on
the basis of the particle weights $\{W_{1,k,j}^{(n_{i})}\}$ computed in the
last iteration; if the $j_{k}$-th particle (i.e., $\mathbf{x}%
_{k,j_{k}}^{(N)} $) is selected, we set 
\begin{equation}
	\mathbf{x}_{\mathrm{be},k}^{(N)}=\mathbf{x}_{k,j_{k}}^{(N)},
	\label{selected_particle}
\end{equation}%
so that the message (see Eq. (\ref{mess_be_N_l}))%
\begin{equation}
	\cev{m}_{\mathrm{be}}\left( \mathbf{x}_{k}^{(N)}\right) \triangleq \delta
	\left( \mathbf{x}_{k}^{(N)}-\mathbf{x}_{\mathrm{be},k}^{(N)}\right) ,
	\label{m_be_N_new}
\end{equation}%
can be made available at the output of BIF$_{2}$. On the other hand, the
evaluation of the message $\cev{m}_{\mathrm{be}}\left( \mathbf{x}_{k}\right) 
$ is accomplished as follows. The messages $m_{2}^{(n_{i}+1)}(\mathbf{x}%
_{k}) $ and $\cev{m}_{3}^{(n_{i}+1)}\left( \mathbf{x}_{k}\right) $ are
computed first (see Eqs. (\ref{m_pm_x_l})-(\ref{PM_z_L}) and Eqs. (\ref%
{m_be1_x_la})-(\ref{w_be1_l_kn}), respectively). Then, the message $\cev{m}_{%
	\mathrm{be}}\left( \mathbf{x}_{k}\right) $ is computed as (see Fig. \ref%
{Fig_4}) 
\begin{eqnarray}
	\cev{m}_{\mathrm{be}}\left( \mathbf{x}_{k}\right) &=&\cev{m}_{\mathrm{be}%
		2}\left( \mathbf{x}_{k}\right) =\cev{m}_{\mathrm{be}1}^{(n_{i}+1)}\left( 
	\mathbf{x}_{k}\right) \,m_{\mathrm{ms}}\left( \mathbf{x}_{k}\right)
	\label{m_be2_x} \\
	&=&\mathcal{\mathcal{N}}\left( \mathbf{x}_{k};\mathbf{\eta }_{\mathrm{be}%
		2,k},\mathbf{C}_{\mathrm{be}2,k}\right) ,  \label{m_be2_xb}
\end{eqnarray}%
where%
\begin{equation}
	m_{\mathrm{ms}}\left( \mathbf{x}_{k}\right) =\mathcal{\mathcal{N}}\left( 
	\mathbf{x}_{k};\mathbf{\eta }_{\mathrm{ms},k},\mathbf{C}_{\mathrm{ms}%
		,k}\right)  \label{m_ms_x}
\end{equation}%
is the message conveying measurement information.

Moreover, the covariance matrices $\mathbf{C}_{\mathrm{ms},k}$ and $\mathbf{C%
}_{\mathrm{be}2,k}$, and the mean vectors $\mathbf{\eta }_{\mathrm{ms},k}$
and $\mathbf{\eta }_{\mathrm{be}2,k}$ are evaluated on the basis of the
associated precision matrices 
\begin{equation}
	\mathbf{W}_{\mathrm{ms},k}\triangleq (\mathbf{C}_{\mathrm{ms},k})^{-1}=%
	\mathbf{H}_{k}\mathbf{W}_{e}\,\mathbf{H}_{k}^{T},  \label{W_ms_x}
\end{equation}%
\begin{equation}
	\mathbf{W}_{\mathrm{be}2,k}\triangleq (\mathbf{C}_{\mathrm{be}2,k})^{-1}=%
	\mathbf{W}_{\mathrm{ms},k}+\mathbf{W}_{\mathrm{be}1,k}^{(n_{i}+1)},
	\label{W_be2_x}
\end{equation}%
and of the transformed mean vectors%
\begin{equation}
	\mathbf{w}_{\mathrm{ms},k}\triangleq \mathbf{W}_{\mathrm{ms},k}\,\mathbf{%
		\eta }_{\mathrm{ms},k}=\mathbf{H}_{k}\mathbf{W}_{e}\left( \mathbf{y}_{k}-%
	\mathbf{v}_{k}\right) \text{,}  \label{w_ms_x}
\end{equation}%
\begin{equation}
	\mathbf{w}_{\mathrm{be}2,k}\triangleq \mathbf{W}_{\mathrm{be}2,k}\,\mathbf{%
		\eta }_{\mathrm{be}2,k}=\mathbf{w}_{\mathrm{ms},k}+\mathbf{w}_{\mathrm{be}%
		1,k}^{(n_{i}+1)},  \label{w_be2_x}
\end{equation}%
respectively. The $k$-th recursion is now over.

In the DBSA, the \emph{first recursion} of the backward pass (corresponding
to $k=T-1$) requires the knowledge of the input messages $\cev{m}_{\mathrm{be%
}}(\mathbf{x}_{T})$ and $\cev{m}_{\mathrm{be}}(\mathbf{x}_{T}^{(N)})$.
Similarly as any BIF algorithm, the evaluation of these messages in DBIF is
based on the statistical information generated in the last recursion of the
forward pass. In particular, the above mentioned messages are still
expressed by Eqs. (\ref{mess_be_l}) and (\ref{mess_be_N_l}) (with $k=T-1$ in
both formulas), respectively. However, the vector $\mathbf{x}_{\mathrm{be}%
	,T}^{(N)}$ is generated by sampling the particle set $S_{T}$ on the basis of
the forward weights $\{w_{\mathrm{fe},T,j}\}$, since backward predictions
are unavailable at the final instant $k=T$. Therefore, if the $j_{T}$-th
particle of $S_{T}$ is selected, we set%
\begin{equation}
	\mathbf{x}_{\mathrm{be},T}^{(N)}=\mathbf{x}_{\mathrm{fe},T,j_{T}}^{(N)}
	\label{x_be_N_T}
\end{equation}%
in the message $\cev{m}_{\mathrm{be}}(\mathbf{x}_{T}^{(N)})$ entering the BIF%
$_{2}$ in the first recursion (see Eq. (\ref{mess_be_N_l})). As far as BIF$%
_{1}$ is concerned, following \cite{Vitetta_2018}, we choose%
\begin{equation}
	\mathbf{W}_{\mathrm{be},T}=\mathbf{W}_{\mathrm{fe}1,T}  \label{C_be_L_T}
\end{equation}%
and%
\begin{equation}
	\mathbf{w}_{\mathrm{be},T}=\mathbf{w}_{\mathrm{fe}1,T}  \label{eta_be_L_T}
\end{equation}%
for the message $\cev{m}_{\mathrm{be}}(\mathbf{x}_{T})$.

The DBSA is summarized in Algorithm 1. It
generates all the statistical information required to solve problems \textbf{%
	P.1} and \textbf{P.2}. Let us now discuss how this can be done in detail. As
far as problem\textbf{\ P.1} is concerned, it is useful to point out that
the DBSA produces a trajectory $\{\mathbf{x}_{\mathrm{be},k}^{(N)},k=1$, $2$%
, $...$, $T\}$ for the \emph{nonlinear} component (see Eq. (\ref%
{selected_particle})). Another trajectory, representing the time evolution
of the \emph{linear} state component only and denoted $\{\mathbf{x}_{\mathrm{%
		be},k}^{(L)},k=1$, $2$, $...$, $T\}$, can be evaluated by sampling the
message $m_{1}^{(n_{i})}(\mathbf{x}_{k}^{(L)})$ (see Eq. (\ref{m_fe_L_EKF_2}%
)) or by simply setting $\mathbf{x}_{\mathrm{be},k}^{(L)}=\mathbf{\tilde{\eta%
}}_{1,k}^{(n_{i})}$ (this task can be accomplished in phase III, after
sampling the particle set $S_{k}$; see also the task g- in phase III of
Algorithm 1.

\begin{algorithm}{
		\SetKw{a}{a-}
		\SetKw{b}{b-}
		\SetKw{c}{c-}%
		
		\SetKw{d}{d-}
		\SetKw{e}{e-}
		\SetKw{f}{f-}
		\SetKw{g}{g-}
		\SetKw{h}{h-}%

		\nl\textbf{Forward filtering}: For $k=1$ to $T$: Run the DBF,
		and store $\mathbf{W}_{\mathrm{fe}1,k}$ (\ref{W_fe_L}), $\mathbf{w}_{\mathrm{fe}1,k}$ (\ref{w_fe_L}), $S_{k}=\{\mathbf{x}%
		_{k,j}^{(N)}\}$ and $\{{w}_{\mathrm{fe},k,j}\}_{j=1}^{N_{p}}$.
		
		\nl\textbf{%
			Initialisation of backward filtering}: compute $\mathbf{x}_{\mathrm{be},T}^{(N)}$ (%
		\ref{x_be_N_T}), $\mathbf{W}%
		_{\mathrm{be},T}$ (\ref{C_be_L_T}) and $\mathbf{w}_{\mathrm{be},T}$ (\ref{eta_be_L_T}); then, compute $\mathbf{C}_{\mathrm{be},T}=(\mathbf{W}%
		_{\mathrm{be},T})^{-1}$, $\mathbf{\eta}_{\mathrm{be},T}=\mathbf{C}_{\mathrm{be},T}\mathbf{w}_{\mathrm{be},T}$.
		
		\nl\textbf{Backward filtering and smoothing}: \\
		\For {$k=T-1$ to $1$}{
			\a \textbf{Phase I}:
			
			%		- \emph{Marginalization}: extract $\mathbf{%
			%\tilde{\eta}}_{be,l+1}$ ($\mathbf{\tilde{C}}_{be,l+1}$) from $\mathbf{\eta}%
			%_{be,l+1}$ ($\mathbf{C}_{be,l+1}$).
			
			- \emph{Backward prediction in} BIF$_1$%
			: compute $\mathbf{W}_{1,k}$ (\ref{W_bp_x_l}) and
			$\mathbf{w}_{1,k}$ (%
			\ref{w_bp_x_l}).
			
			- \emph{Computation of iteration-independent information required in task b}:%
			For $j=1$ to $N_{p}$: compute $\mathbf{z}_{k,j}^{(L)}$ (\ref{PM_z_L}%
			), $\mathbf{\tilde{W}}_{k,j}$ (\ref{eq:W_pm_L_j}), $\mathbf{\tilde{w}}%
			_{k,j}$ (%
			\ref{eq:w_pm_L_j}), $\mathbf{\tilde{C}}_{k,j}=(\mathbf{%
				\tilde{W}}_{k,j})^{-1}$ and $\mathbf{\tilde{\eta}}_{k,j}=\mathbf{%
				\tilde{C}}_{k,j}\mathbf{\tilde{w}}_{k,j}$.
			
			- \emph{%
				Initialisation of particle weights}: Set  $W_{1,k,j}^{(0)}=w_{fe,k,j}$.\\
			
			\textbf{Phase II}:\\
			\For {$n=1$ to $n_{i}$}{
				
				\b Compute $\mathbf{\eta}%
				_{2,k}^{(n)}$ (\ref{eta_pm_l_k}) and $\mathbf{C}_{2,k}^{(n)}$ (\ref{C_pm_l_k}).
				
				% For $j=1$ to $N_{p}$: Compute $\tilde{w}_{3,l,j}^{(n)}$ (\ref{w__sm_j_N})
				% and $W_{sm,l,j}^{(n)}$ (\ref{W_fe_2_x_N_l}); 
				
				\c Compute $\mathbf{C}_{3,k}^{(n)}$ (\ref{C_be1_l_ka}), $\mathbf{\eta}_{3,k}^{(n)}$ (\ref{eta_be1_l_ka}), $\mathbf{W%
				}_{3,k}^{(n)}=(\mathbf{C}_{3,k}^{(n)})^{-1}$, $\mathbf{w}%
				_{3,k}^{(n)}=\mathbf{W}_{3,k}^{(n)}\mathbf{\eta }_{3,k}^{(n)}$, $%
				\mathbf{W}_{4,k}^{(n)}$ (\ref{W_sm_l_k}), $\mathbf{w}_{4,k}^{(n)}$ (\ref{w_sm_l_k}), $\mathbf{C}_{4,k}^{(n)}=(\mathbf{W}_{4,k}^{(n)})^{-1}$ and $%
				\mathbf{\eta}_{4,k}^{(n)}=\mathbf{C}_{4,k}^{(n)}\mathbf{w}_{4,k}^{(n)}$. Then, extract $\mathbf{\tilde{\eta}}_{1,k}^{(n)}$ ($\mathbf{\tilde{C}}%
				_{1,k}^{(n)}$) from $\mathbf{\eta}_{4,k}^{(n)}$ ($\mathbf{C}_{4,k}^{(n)}$%
				).
				
				\d For $j=1$ to $N_{p}$: compute $%
				\mathbf{\eta}_{3,k,j}^{(N)}[n]$ (\ref{eta_bp}) and $\mathbf{C}_{3,k,j}^{(N)}[n]$
				(\ref{cov_bp}). Then, compute $D_{3,k,j}^{(n)}$ (\ref{D_3}) and $Z_{3,k,j}^{(n)}$ (\ref{Z_bp}).
				
				\e For $j=1$ to $N_{p}$: compute $\mathbf{\check{\eta}}%
				_{z,k,j}^{(n)}$ (\ref{eta_mess_z_N}), $\mathbf{\check{C}}_{z,k,j}^{(n)}$ (\ref{C_mess_Z_N_bis}), $\mathbf{\check{W}}_{z,k,j}^{(n)}=(\mathbf{\check{C}}_{z,k,j}^{(n)})^{-1}$%
				, $\mathbf{\check{w}}_{z,k,j}^{(n)}=\mathbf{\check{W}}_{z,k,j}^{(n)}\mathbf{\check{\eta}}%
				_{z,k,j}^{(n)}$, $\mathbf{\check{W}}_{2,k,j}^{(n)}$ (\ref{W_pm_x_N_l_j}), $\mathbf{\check{w}%
				}_{2,k,j}^{(n)}$ (\ref{w_pm_x_N_l_j}). Then, compute $D_{2,k,j}^{(n)}$ (\ref{D_2}) and $Z_{2,k,j}^{(n)}$ (\ref{Z_pm}).

				\f For $j=1$ to $N_{p}$: Compute $%
				\mathbf{\bar{\eta}}_{5,k,j}^{(n)}$ (\ref{eq:ms_fe_PF}), $\mathbf{\bar{C}%
				}_{5,k,j}^{(n)}$ (\ref{eq:C_sm_PF}), $\mathbf{\bar{W}}_{5,k,j}^{(n)}=(\mathbf{\bar{C}}_{5,k,j}^{(n)})^{-1}$, $D_{5,k,j}^{(n)}$ (\ref{D_5}) and $Z_{5,k,j}^{(n)}$ (\ref{Z_5}).
				Then, compute $D_{6,k,j}^{(n)}$ (\ref{D_6}), $Z_{6,k,j}^{(n)}$ (\ref{Z_tot}), $w_{6,k,j}^{(n)}$ (\ref{weight_6b}), $w_{1,k,j}^{(n)}$ (\ref{w_sm_j_Na}) and $W_{1,k,j}^{(n)}$ (\ref{W_fe_2_x_N_l}). Store the weights $\{W_{1,k,j}^{(n)}\}$ for the next iteration.
				
				% close loop
				
			}
			
			\g \textbf{Phase III} - BIF$%
			_2$: Select the $j_{k}$-th
			particle 
			$\mathbf{x}_{k,j_{k}}^{(N)}$ by sampling the set $S_{k}$ on
			the basis of the weights $\{W_{1,k,j}^{(n_{i})}\}$, 
			set $\mathbf{x}_{\mathrm{be},k}^{(N)}=\mathbf{x}%
			_{k,j_{k}}^{(N)}$ and store $\mathbf{x}_{\mathrm{be},k}^{(N)}$ 
			for the next
			recursion.
			
			\h \textbf{Phase III} - BIF$_1$: 
			Compute $\mathbf{\eta}%
			_{2,k}^{(n_{i}+1)}$,  $\mathbf{C}_{2,k}^{(n_{i}+1)}$, $\mathbf{W}%
			_{3,k}^{(n_{i}+1)}$  and $\mathbf{w}_{3,k}^{(n_{i}+1)}$  (see steps
			1) and 2)). 
			Then, compute  $\mathbf{W}_{\mathrm{ms},k}$ (\ref{W_ms_x}), $\mathbf{w}_{\mathrm{ms},k}$
			(\ref{w_ms_x}), $\mathbf{W}_{\mathrm{be}2,k}$ (\ref{W_be2_x}), 
			$\mathbf{w}_{\mathrm{be}2,k}$
			(\ref{w_be2_x}), $\mathbf{C}_{\mathrm{be},k}=(\mathbf{W}_{\mathrm{be}2,k})^{-1}$ and
			$\mathbf{\eta}_{\mathrm{be},k}=\mathbf{C}_{\mathrm{be},k}\mathbf{w}_{\mathrm{be}2,k}$, and store $\mathbf{C}_{\mathrm{be},k}$ and $\mathbf{\eta}_{\mathrm{be},k}$ for the next recursion.
		}
		% closeloop
		\caption{Double Bayesian Smoothing}}
\end{algorithm}

Since the DBSA solves problem \textbf{P.1}, it also solves problem \textbf{%
	P.2}; in fact, once it has been run, an approximation of the marginal
smoothed pdf at any instant can be simply obtained by marginalization.
Unluckily, the last result is achieved at the price of a significant
computational cost, since $M$ backward passes are required. However, if we
are interested in solving problem \textbf{P.2} only, a simpler particle
smoother can be developed following the approach illustrated in ref. \cite%
{Vitetta_2018}, so that a single backward pass has to be run. In this pass,
the evaluation of the message $\cev{m}_{\mathrm{be}}(\mathbf{x}_{k}^{(N)})$
(i.e., of the particle $\mathbf{x}_{\mathrm{be},k}^{(N)}$) involves the
whole particle set $S_{k}$ and their weights $\{W_{1,k,j}^{(n_{i})}\}$ (see
Eq. (\ref{W_fe_2_x_N_l})) evaluated in the last phase of the $(T-k)-$th
recursion. More specifically, a new smoother is obtained by employing a
different method for evaluating $\mathbf{x}_{\mathrm{be},k}^{(N)}$ (see
phase III-BIF$_2$); it consists in computing the smoothed estimate%
\begin{equation}
	\mathbf{x}_{\mathrm{sm},k}^{(N)}=\sum%
	\limits_{j=1}^{N_{p}}W_{1,k,j}^{(n_{i})}\,\mathbf{x}_{k,j}^{(N)}
	\label{x_sm_N}
\end{equation}%
of $\mathbf{x}_{k}^{(N)}$ and, then, setting%
\begin{equation}
	\mathbf{x}_{\mathrm{be},k}^{(N)}=\mathbf{x}_{\mathrm{sm},k}^{(N)}.
	\label{x_be_N_new}
\end{equation}%
The resulting smoother is called \emph{simplified }DBSA (SDBSA) in the
following.

The computational complexity of the DBSA and the SDBSA can be reduced by
reusing the forward weights $\{w_{\mathrm{fe},k,j}\}$ in all the iterations
of phase II, so that step 7) can be skipped; this means that, for any $n$,
we set $w_{5,k,j}^{(n)}=w_{\mathrm{fe},k,j}$ in the evaluation of the $j-$th
particle weight $w_{6,k,j}^{(n)}$ according to Eq. (\ref{weight_6}) in step
8) of phase II. Our simulation results have evidenced that, at least for the
SSMs considered in Section \ref{num_results}, this modification does not
affect the estimation accuracy of the derived algorithms; for this reason,
it is always employed in our simulations.

The DBSA and the SDBSA refer to case \textbf{C.1}, i.e. to the case in which
the substates estimated by the interconnected forward/backward filters share
the substate $\mathbf{x}_{k}^{(N)}$. Let us focus now on case \textbf{C.2},
i.e. on the case on which the filters are run on disjoint substates. A
filtering technique, called \emph{simplified} DBF (SDBF), and based on the
interconnection of a particle filter (F$_{2}$) with a single Kalman filter (F%
$_{1}$), is developed for this case in ref. \cite[Par. III-B]%
{Vitetta_DiViesti_2019}. The BIF\ algorithm paired with it can be easily
derived following the approach illustrated above for the DBSA; the resulting
smoothing algorithm is dubbed \emph{disjoint} DBSA (DDBSA) in the following.
It is important to mention that, in deriving the DBSA, the following
relevant changes are made with respect to the DBSA (see Fig. \ref{Fig_2}):

1) The iterative procedure embedded in the $(T-k)-$th recursion of the
backward pass involves both the computation of the backward predicted pdf ($%
BP1$) and of the message $MS1$ in BIF$_{1}$; for this reason, it requires
marginalizing the pdfs $f(\mathbf{x}_{k+1}^{(N)}|\mathbf{x}_{k}^{(N)},%
\mathbf{x}_{k}^{(L)})$ and $f(\mathbf{y}_{k}|\mathbf{x}_{k}^{(N)},\mathbf{x}%
_{k}^{(L)})$, respectively, with respect to $\mathbf{x}_{k}^{(L)}$. This
result is achieved in the first iteration by setting $\mathbf{x}_{k}^{(N)}=%
\mathbf{x}_{\mathrm{fe},k}^{(N)}$ in both these pdfs, where $\mathbf{x}_{%
	\mathrm{fe},k}^{(N)}$ denotes the estimate of $\mathbf{x}_{k}^{(N)}$
computed by F$_{2}$ in the forward pass. In the following iterations, we set 
$\mathbf{x}_{k}^{(N)}=\mathbf{x}_{\mathrm{sm},k}^{(N)}$, where $\mathbf{x}_{%
	\mathrm{sm},k}^{(N)}$ represents the estimate of $\mathbf{x}_{k}^{(N)}$
evaluated on the basis of the statistical information provided by BIF$_{2}$
(through the message $SM2$).

2) The pseudo-measurement message $PM1$ (corresponding to $m_{2}^{(n)}(%
\mathbf{x}_{k})$ (\ref{m_pm_x_l}) in the DBSA) conveys information about $%
\mathbf{x}_{k}^{(L)}$only. Moreover, it is a Gaussian message, and its mean
and covariance matrix are given by $\mathbf{\eta }_{L,k}^{(n)}$ and $\mathbf{%
	C}_{LL,k}^{(n)}$ (see Eqs. (\ref{eta_pm_l_L_k}) and (\ref{C_pm_l_L_k_bis}),
respectively).

Finally, it is worth mentioning that a \emph{simplified} version of the
DDBSA (called SDDBSA) can be easily developed by making the same
modifications as those adopted in deriving the SDBSA from the DBSA.

\section{Comparison of the Developed Double Smoothing Algorithms with
	Related Techniques\label{Comparison}}

The DBSA and the DDBSA developed in the previous Section are conceptually
related to the Rao-Blackwellised particle smoothers proposed by Fong \emph{%
	et al. }\cite{Fong_2002} and by Lindsten \emph{et al.} \cite{Lindsten_2016}
(these algorithms are denoted Alg-F and Alg-L respectively, in the
following) and to the RBSS algorithm devised by Vitetta \emph{et al.} in
ref. \cite{Vitetta_2018}. In fact, all these techniques share with the DBSA
and the DDBSA the following important features: 1) all of them estimate the 
\emph{joint} smoothing density over the whole observation interval by
generating multiple \emph{realizations} from it; 2) they accomplish a single
forward pass and as many backward passes as the overall number of
realizations; 3)\ they combine Kalman filtering with particle filtering.
However, Alg-F, Alg-L and the RBSS\ algorithm employ, in both their forward
and backward passes, as many Kalman filters as the number of particles ($%
N_{p}$) to generate a particle-dependent estimate of the linear state
component only. On the contrary, the DBSA (DDBSA) employs a \emph{single}
extended Kalman filter (a \emph{single} Kalman filter), that estimates the
whole system state (the linear state component only); this substantially
reduces the memory requirements of particle smoothing and, consequently, the
overall number of memory accesses accomplished on the hardware platform on
which smoothing is run. As far as the last point is concerned, the memory
requirements\ of a smoothing algorithm can be roughly assessed by estimating
the overall number of real quantities that need to be stored in both its
forward pass and its backward pass. It can be shown that overall number of
real quantities to be stored by MPF, DBF and SDBF in the forward pass of the
considered smoothing algorithms is of order $\mathcal{O}(M_{MPF})$, $%
\mathcal{O}(M_{DBF})$, and $\mathcal{O}(M_{SDBF})$, respectively, with%
\footnote{%
	Note that the expressions (\ref{MPF_mem})-(\ref{SDBF_mem}) also account for
	the contributions due to measurement-based information (see Eqs. (\ref%
	{W_ms_x}) and (\ref{w_ms_x})).}%
\begin{equation}
	M_{MPF}=N_{p}T\,(2D_{L}^{2}+2D_{L}+D_{N}+1),  \label{MPF_mem}
\end{equation}%
\begin{equation}
	M_{DBF}=T\,(2D^{2}+2D+N_{p}\,D_{N}+N_{p})  \label{DBF_mem}
\end{equation}%
and 
\begin{equation}
	M_{SDBF}=T\,(2D_{L}^{2}+2D_{L}+N_{p}\,D_{N}+N_{p}\,).  \label{SDBF_mem}
\end{equation}%
Moreover, the overall number of real quantities to be stored by Alg-L, RBSS,
the DBSA and the DDBSA is approximately of order $\mathcal{O}(M_{Alg-L})$, $%
\mathcal{O}(M_{RBSS})$, $\mathcal{O}(M_{DBSA})$ and $\mathcal{O}(M_{DDBSA})$%
, respectively, with 
\begin{equation}
	M_{Alg-L}=M_{MPF}+D_{L}^{2}+D,  \label{Alg-L_mem}
\end{equation}%
\begin{equation}
	M_{RBSS}=M_{MPF}+D_{L}^{2}+D,  \label{RBSS_mem}
\end{equation}%
\begin{equation}
	M_{DBSA}=M_{DBF}+N_{p}+D^{2}+D+D_{N}  \label{DBSA_mem}
\end{equation}%
and%
\begin{equation}
	M_{DDBSA}=M_{SDBF}+N_{p}+D_{L}^{2}+D.  \label{DDBSA_mem}
\end{equation}%
The memory requirements of the SDBSA and the SDDBSA (the SPS algorithm) are
the same as those of the DBSA and the DDBSA (the RBSS algorithm),
respectively. Note also that the quantities $M_{DBSA}$ (\ref{DBSA_mem}) and $%
M_{DDBSA}$ (\ref{DDBSA_mem}) are smaller than $M_{Alg-L}$ (\ref{Alg-L_mem})
and $M_{RBSS}$ (\ref{RBSS_mem}), since $M_{MPF}$ is larger than $M_{DBF}$
and $M_{SDBF}$ because of its dependence on $N_{p}$.

The differences in the overall execution time measured for the simulated
smoothing algorithms are related not only to their requirements in terms of
memory resources, but also to their computational complexity. In our work,
the computational cost of the smoothing algorithms derived in the previous
section has been carefully assessed in terms of number of \emph{floating
	point operations} (flops) to be executed over the whole observation
interval. The general criteria adopted in estimating the computational cost
of an algorithm are the same as those illustrated in \cite[App. A, p. 5420]%
{Hoteit_2016} and are not repeated here for space limitations. A detailed
analysis of the cost required by each of the tasks accomplished by our
smoothing algorithms is provided in Appendix \ref{app:CDBFA}. Our
analysis leads to the conclusion that the overall computational cost of the
DBSA and of the DDBSA is approximately of order $\mathcal{O}(N_{DBSA})$ and $%
\mathcal{O}(N_{DDBSA})$, respectively, with 
\begin{eqnarray}
	N_{DBSA} &=&T\,\left\{ N_{DBF}+\,M\left[ 38D^{3}/3+20D_{N}^{3}/3+\right.
	\right.   \notag \\
	&&\left. \left. n_{i}N_{p}(2D_{L}^{2}D_{N}+2D_{L}D_{N}^{2} +D_{N}^{3}/3+5D_{L}^{3})+6n_{i}D^{3}\right] \right\} ,
	\label{DBSA_comp_compl}
\end{eqnarray}%
and%
\begin{eqnarray}
	N_{DDBSA} &=&T\,\left\{ N_{SDBF}+M\left[ 38D_{L}^{3}/3+20D_{N}^{3}/3\right.
	\right.   \notag \\
	&&\left. \left. +n_{i}N_{p}(2D_{L}^{2}D_{N}+2D_{L}D_{N}^{2} +D_{N}^{3}/3+5D_{L}^{3})+6n_{i}D_{L}^{3}\right] \right\} ;
	\label{DDBSA_comp_compl}
\end{eqnarray}%
here, $N_{DBF}$ and $N_{SDBF}$ represent the computational complexity of a
single recursion of the DBF and SDBF, respectively (see \cite[Eqs. (97) and
(98)]{Vitetta_DiViesti_2019}). Each of the expressions (\ref{DBSA_comp_compl}%
)-(\ref{DDBSA_comp_compl}) has been derived as follows. First, the costs of
all the tasks identified in Appendix \ref{app:CDBFA} have been summed;
then, the resulting expression has been simplified, keeping only the
dominant contributions due to matrix inversions, matrix products and
Cholesky decompositions, and discarding all the contributions that originate
from the evaluation of the matrices $\mathbf{A}_{k}^{(Z)}(\mathbf{x}%
_{k}^{(N)})$ (with $Z=L$ and $N$), $\mathbf{F}_{k}$, $\mathbf{H}_{k}$ and $%
\mathbf{B}_{k}$ and the functions $\mathbf{f}_{k}^{(Z)}(\mathbf{x}_{k}^{(N)})
$ (with $Z=L$ and $N$), $\mathbf{f}_{k}(\mathbf{x}_{k})$ and $\mathbf{g}_{k}(%
\mathbf{x}_{k}^{(N)})$. Moreover, the sampling of the particle set in each
recursion of the backward pass has been ignored.

From Eqs. (\ref{DBSA_comp_compl})-(\ref{DDBSA_comp_compl}) it is easily
inferred that the computational complexities of the DBSA and the DDBSA are
approximately of order $\mathcal{O}(n_{i}M\,N_{p}D_{L}^{3}T)$. A similar
approach can be followed for Alg-L and the RBSS algorithm; this leads to the
conclusion that their complexities are approximately of order $\mathcal{O}%
(M\,N_{p}D_{L}^{3}T)$, i.e. of the same order of the complexities of the
DBSA and of the DDBSA if $n_{i}=1$ is assumed.

On the other hand, the SDBSA and the SDDBSA are conceptually related to the
SPS algorithm devised by Vitetta \emph{et al.} in ref. \cite{Vitetta_2018}.
In fact, all these algorithms aim at solving problem \textbf{P.2} only
(consequently, they are unable to generate the joint smoothed pdf $f(\mathbf{%
	x}_{1:T}|\mathbf{y}_{1:T})$) and carry out a \emph{single backward pass}.
This property makes them much faster than Alg-L, the RBSS algorithm, the
DBSA and the DDBSA in the computation of marginal smoothed densities.
Finally, note that, similarly as the DBSA and the DDBSA techniques, the use
of the SDBSA and the SDDBSA requires a substantially smaller number of
memory accesses than the SPS algorithm, since the last algorithm employs MPF
in its forward pass. Moreover, the computational cost of the SDBSA and the
SDDBSA is approximately of order $\mathcal{O}(n_{i}N_{p}D_{L}^{3}T)$,
whereas that of the SPS algorithm is approximately of order $\mathcal{O}%
(N_{p}D_{L}^{3}T)$; consequently, they are all of the same order if $n_{i}=1$
is assumed.

\section{Numerical Results\label{num_results}}

In this section we first compare, in terms of accuracy and execution time,
the DBSA, the SDBSA, the DDBSA\ and the SDDBSA with Alg-L, the RBSS
algorithm, and the SPS algorithm for a specific conditionally linear
Gaussian SSM. The considered SSM is the same as the SSM\#2 defined in \cite%
{Vitetta_2018} and describes the bidimensional motion of an agent. Its state
vector in the $k$-th observation interval is defined as $\mathbf{x}%
_{k}\triangleq \lbrack \mathbf{v}_{k}^{T},\mathbf{p}_{k}^{T}]^{T}$, where $%
\mathbf{v}_{k}\triangleq \lbrack v_{x,k},v_{y,k}]^{T}$ and $\mathbf{p}%
_{k}\triangleq \lbrack p_{x,k},p_{y,k}]^{T}$ (corresponding to $\mathbf{x}%
_{k}^{(L)}$ and $\mathbf{x}_{k}^{(N)}$, respectively) represent the agent
velocity and position, respectively (their components are expressed in m/s
and in m, respectively). The state update equations are%
\begin{equation}
	\mathbf{v}_{k+1}=\rho \,\mathbf{v}_{k}+T_{s}\,\mathbf{a}_{k}(\mathbf{p}%
	_{k})+\left( 1-\rho \right) \,\mathbf{n}_{v,k}  \label{mod_1_v}
\end{equation}%
and%
\begin{equation}
	\mathbf{p}_{k+1}=\mathbf{p}_{k}+T_{s}\,\mathbf{v}_{k}+(T_{s}^{2}/2)\,\mathbf{%
		\ a}_{k}(\mathbf{p}_{k})+\mathbf{n}_{p,k},  \label{mod_1_p}
\end{equation}%
where $\rho $ is a forgetting factor (with $0<\rho <1$), $T_{s}$ is the
sampling interval, $\mathbf{n}_{v,k}$ is an \emph{additive Gaussian noise}
(AGN) vector characterized by the covariance matrix $\mathbf{I}_{2}$, 
\begin{equation}
	\mathbf{a}_{k}\left( \mathbf{p}_{k}\right) =-a_{0}\frac{\mathbf{p}_{k}}{
		\left\Vert \mathbf{p}_{k}\right\Vert }\frac{1}{1+\left( \left\Vert \mathbf{p}
		_{k}\right\Vert /d_{0}\right) ^{2}}  \label{acc_model}
\end{equation}%
is the acceleration due to a force applied to the agent (and pointing
towards the origin of our reference system), $a_{0}$ is a scale factor
(expressed in m/s$^{2}$), $d_{0}$ is a \emph{reference distance} (expressed
in m), and $\mathbf{n}_{p,k}$ is an AGN vector characterized by the
covariance matrix $\sigma _{p}^{2}\,\mathbf{I}_{2}$ and accounting for model
inaccuracy. The measurement vector available in the $k$-th interval for
state estimation is 
\begin{equation}
	\mathbf{y}_{k}=\mathbf{x}_{k}+\mathbf{e}_{k},  \label{mod_1_y}
\end{equation}%
where $\mathbf{e}_{k}\triangleq \lbrack \mathbf{e}_{v,k}^{T},\mathbf{e}%
_{p,k}^{T}]^{T}$ and $\mathbf{e}_{v,k}$ ($\mathbf{e}_{p,k}$) is an AGN
vector characterized by the covariance matrix $\sigma _{ev}^{2}\,\mathbf{I}%
_{2}$ ($\sigma _{ep}^{2}\,\mathbf{I}_{2}$).

In our computer simulations, following \cite{Vitetta_2018} and \cite%
{Vitetta_2019}, the estimation accuracy of the considered smoothing
techniques has been assessed by evaluating two \emph{root mean square errors}
(RMSEs), one for the linear state component, the other for the nonlinear
one, over an observation interval lasting $T=200$ $T_{s}$; these are denoted 
$RMSE_{L}($alg$)$ and $RMSE_{N}($alg$)$, respectively, where `alg' is the
acronym of the algorithm these parameters refer to. Our assessment of the
computational requirements is based, instead, on evaluating the average 
\emph{computation time} required for processing a single \emph{block} of
measurements (this quantity is denoted CTB$($alg$)$ in the following).
Moreover, the following values have been selected for the parameters of the
considered SSM: $\rho =0.995$, $T_{s}=0.01$ s, $\sigma _{p}$ $=5\cdot 10^{-3}
$ m, $\sigma _{e,p}=2\cdot 10^{-2}$ m, $\sigma _{e,v}=2\cdot 10^{-2}$ m/s, $%
a_{0}=0.5$ m/s$^{2}$, $d_{0}=5\cdot 10^{-3}$ m and $v_{0}=1$ m/s (the
initial position $\mathbf{p}_{0}\triangleq \lbrack p_{x,0},p_{y,0}]^{T}$ and
the initial velocity $\mathbf{v}_{0}\triangleq \lbrack v_{x,0},v_{y,0}]^{T}$
have been set to $[0.01$ m$,$ $0.01$ m$]^{T}$ and $[0.01$ m/s$,$ $0.01$ m/s$%
]^{T}$, respectively).

\begin{figure}[tbp]
	\centering
	\includegraphics[width=0.6\textwidth]{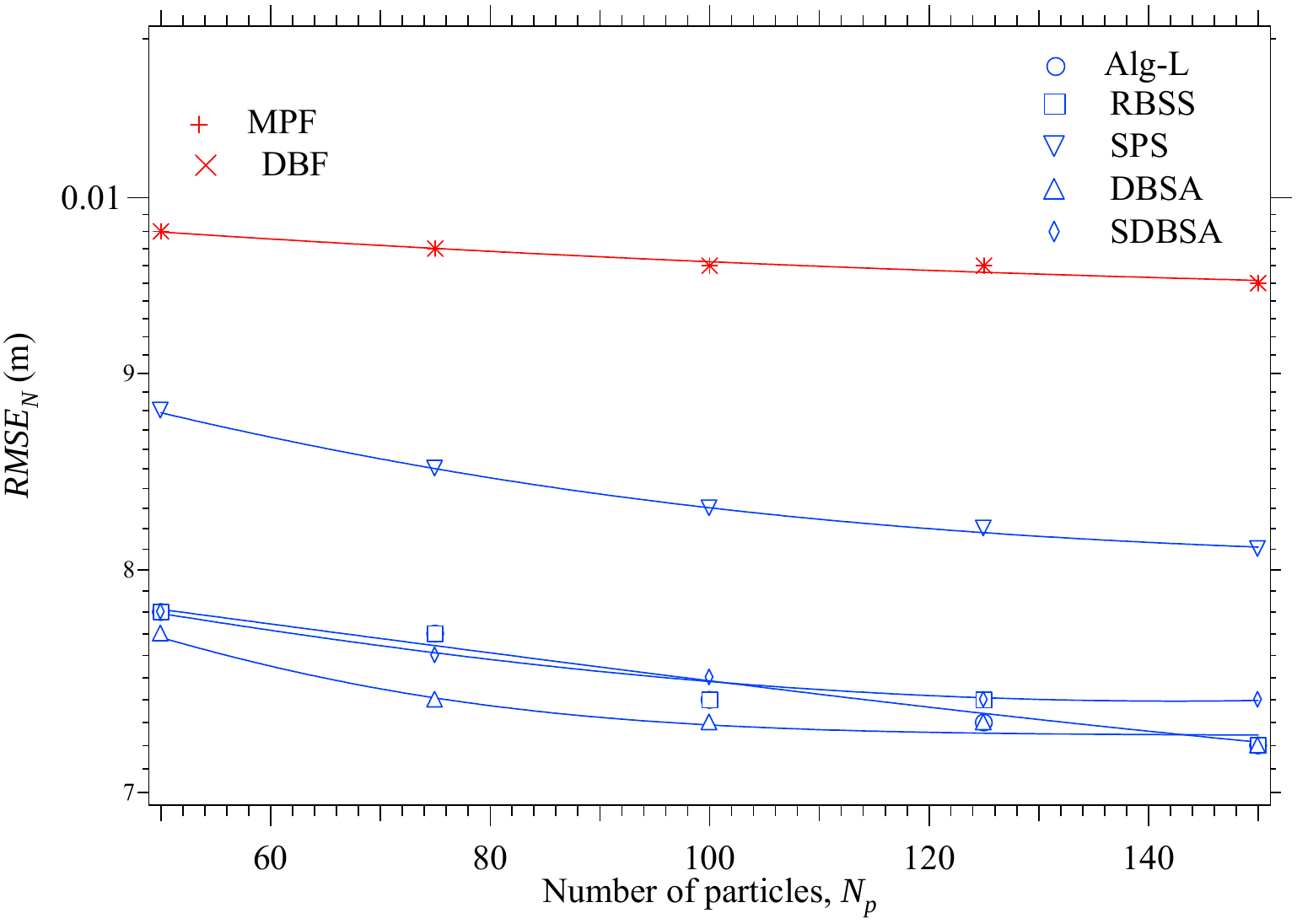}
	\caption{RMSE performance versus $N_{p}$ for the nonlinear component ($%
		RMSE_{N}$) of the state of SSM \#1; five smoothing algorithms (Alg-L, the
		DBSA, the SDBSA, the RBSS algorithm and the SPS algorithm) and two filtering
		techniques (MPF and DBF) are considered.}
	\label{Fig_rmsen}
\end{figure}

\begin{figure}[tbp]
	\centering
	\includegraphics[width=0.6\textwidth]{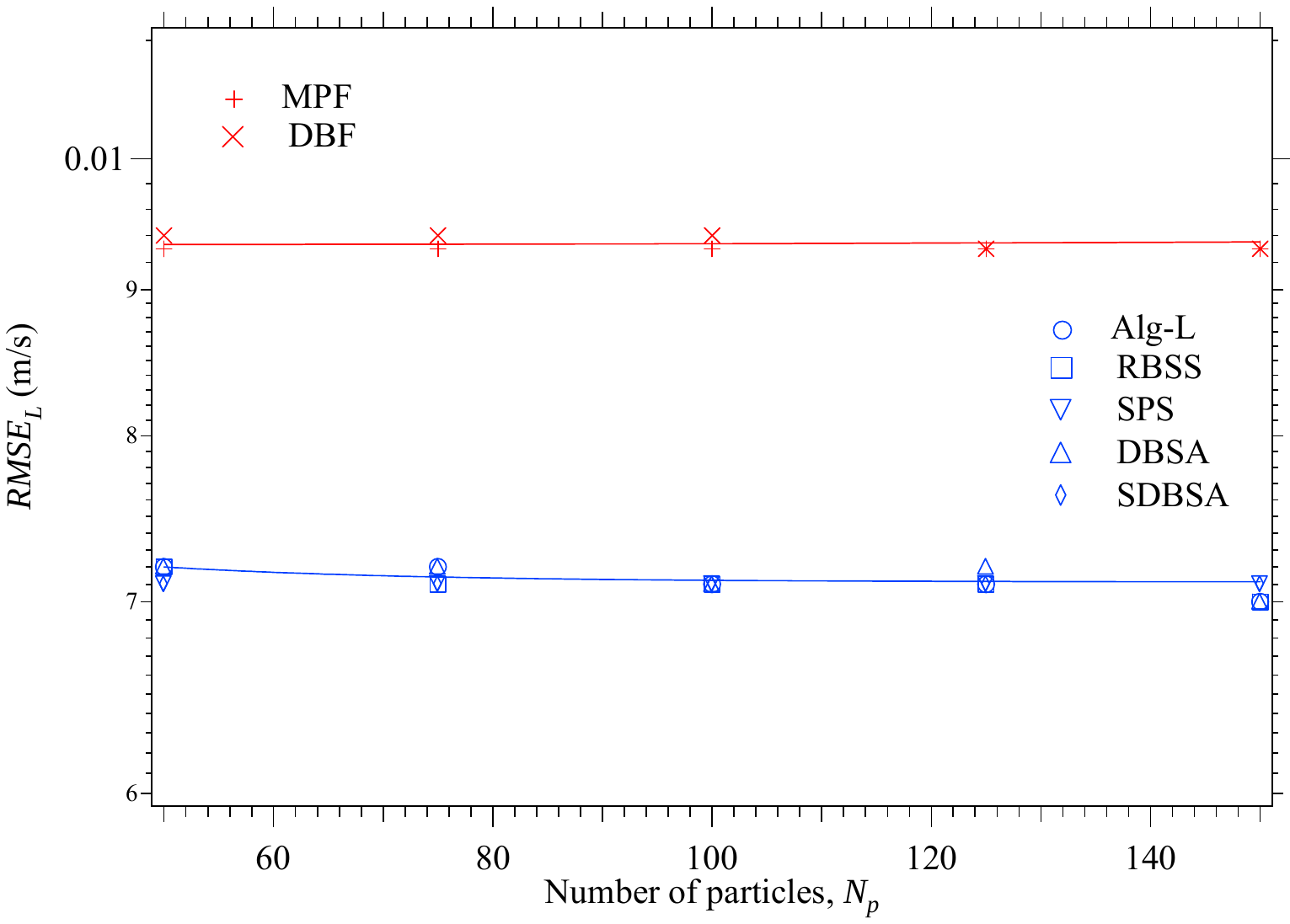}
	\caption{RMSE performance versus $N_{p}$ for the linear component ($RMSE_{L}$
		) of the state of SSM \#1; five smoothing algorithms (Alg-L, the DBSA, the
		SDBSA, the RBSS algorithm and the SPS algorithm) and two filtering
		techniques (MPF and DBF) are considered.}
	\label{Fig_rmsel}
\end{figure}

Some numerical results showing the dependence of $RMSE_{L}$ and $RMSE_{N}$
on the number of particles ($N_{p}$) for some of the considered smoothing
algorithms are illustrated in Figs. \ref{Fig_rmsen} and \ref{Fig_rmsel},
respectively (simulation results are indicated by markers, whereas
continuous lines are drawn to fit them, so facilitating the interpretation
of the available data). In this case, $n_{i}=1$ has been selected for all
the derived particle smoothers, $M=N_{p}$ has been chosen for all the
smoothing algorithms generating multiple trajectories and the range $[10,150]
$ has been considered for $N_{p}$ (since no real improvement is found for $%
N_{p}\gtrsim 150$). Moreover, $RMSE_{L}$ and $RMSE_{N}$ results are also
provided for MPF and DBF, since these filtering techniques are employed in
the forward pass of Alg-L, the RBSS algorithm and the SPS algorithm, and the
DBSA and the SDBSA, respectively; this allows us to assess the improvement
in estimation accuracy provided by the backward pass with respect to the
forward pass for each smoothing algorithm. These results show that:

1) The DBSA, the SDBSA, Alg-L and the RBSS algorithm achieve similar
accuracies in the estimation of both the linear and nonlinear state
components.

2) The SPS algorithm is slightly outperformed by the other smoothing
algorithms in terms of $RMSE_{N}$ only; for instance, $RMSE_{N}($SPS$)$ is
about $1.11$ times larger than $RMSE_{N}($SDBSA$)$ for $N_{p}=100$.

3) Even if the RBSS algorithm and the DBSA provide by far richer statistical
information than their simplified counterparts (i.e., than the SPS algorithm
and the SDBSA, respectively), they do not provide a significant improvement
in the accuracy of state estimation; for instance, $RMSE_{N}($SPS$)$ ($%
RMSE_{N}($SDBSA$)$) is about $1.12$ ($1.03$) time larger than $RMSE_{N}($RBSS%
$)$ ($RMSE_{N}($DBSA$)$) for $N_{p}=100$.

4) The accuracy improvement in terms of $RMSE_{L}$ ($RMSE_{N}$) provided by
all the smoothing algorithms except the SPS (by Alg-L, the RBSS algorithm,
the DBSA and the SDBSA) is about $24\%$ (about $23\%$) with respect to MPF
and DBF, for $N_{p}=100$. Moreover, the accuracy improvement in terms of $%
RMSE_{L}$ ($RMSE_{N}$) achieved by the SPS algorithm is about $24\%$ (about $%
14\%$) with respect to the MPF for $N_{p}=100$.

5) In the considered scenario, DBF is slightly outperformed by (perform
similarly as) MPF in the estimation of the linear (nonlinear) state
component; a similar result is reported in \cite{Vitetta_DiViesti_2019} for
a different SSM.

Our simulations have also evidenced that the DBSA and the SDBSA perform
similarly as the DDBSA and the SDDBSA; for this reason, RSME results
referring to the last two algorithms are not shown in Figs. \ref{Fig_rmsen}
and \ref{Fig_rmsel}. This leads to the conclusion that, in the considered
scenario, the presence of redundancy in double Bayesian smoothing does not
provide any improvement with respect to the case in which the two
interconnected filters operate on disjoint substates in the forward and in
the backward passes. Note that the same conclusion had been reached in ref. 
\cite[Sec. IV]{Vitetta_DiViesti_2019} for DBF only.

\begin{figure}[tbp]
	\centering
	\includegraphics[width=0.6\textwidth]{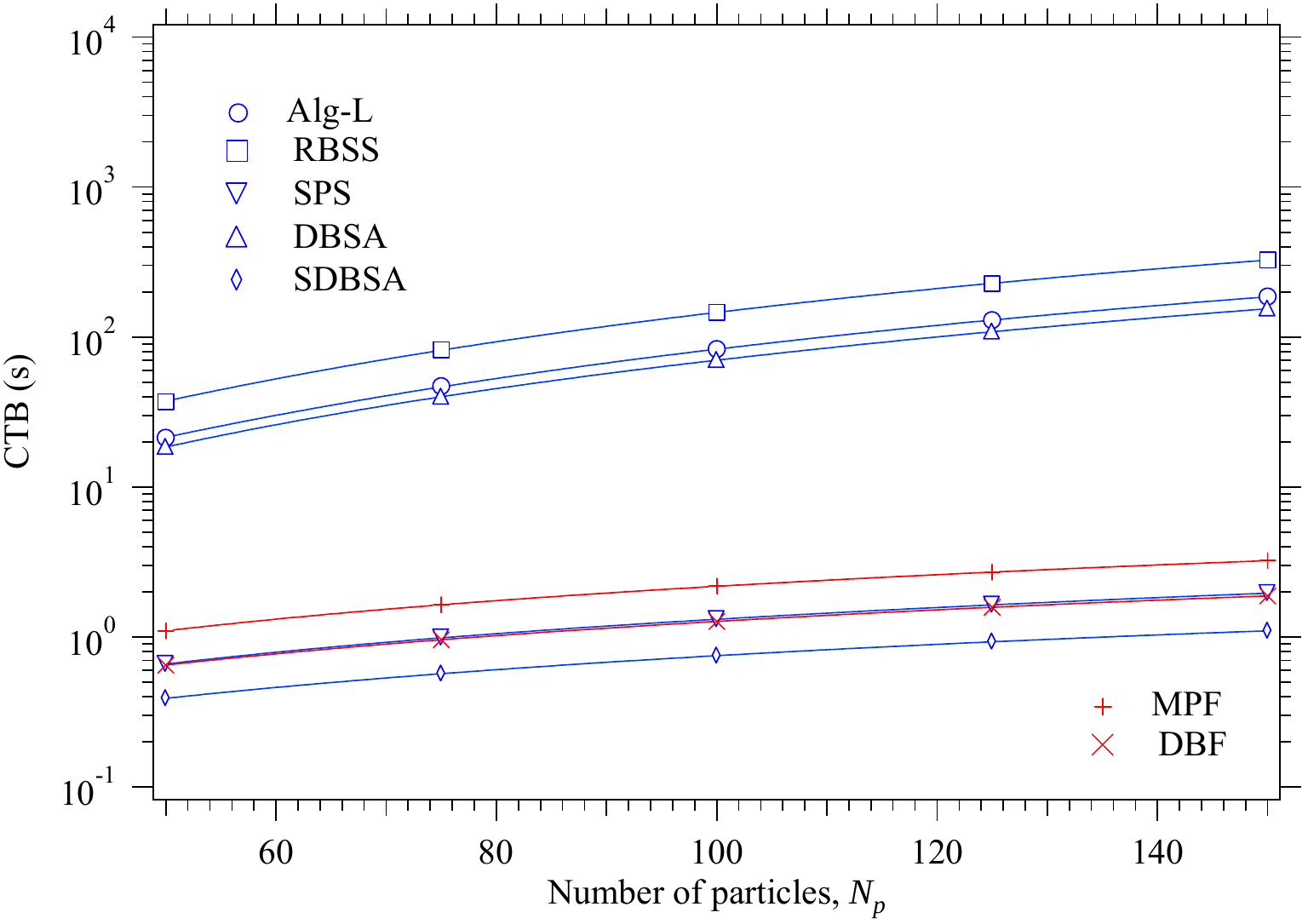}
	\caption{CTB versus $N_{p}$ for five smoothing algorithms (Alg-L, DBSA,
		SDBSA, the RBSS algorithm and the SPS algorithm) and two filtering
		techniques (MPF and DBF); SSM \#1 is considered.}
	\label{Fig_CTB}
\end{figure}

Despite their similar accuracies, the considered smoothing algorithms
require different computational efforts; this is easily inferred from the
numerical results appearing in Fig. \ref{Fig_CTB} and illustrating the
dependence of the CTB on $N_{p}$ for all the above mentioned filtering and
smoothing algorithms. In fact, these results show that CTB$($DBSA$)$\ is
approximately $0.85$ ($0.48$) times smaller than CTB$($Alg-L$)$ (CTB$($RBSS$)
$); this is in agreement with the mathematical results illustrated in
Section \ref{Comparison} about the complexity of Alg-L, the RBSS algorithms
and the DBSA, i.e. with the fact the complexities of all these smoothers are
approximately of order $\mathcal{O}(M\,N_{p}\,D_{L}^{3}T)$ (provided that $%
n_{i}=1$ is selected for the DBSA). Moreover, we have found that a $5.5\%$
reduction in CTB is obtained if the DDBSA is employed in place of the DBSA
(i.e., if double Bayesian smoothing is not redundant). Similar
considerations hold for the SDBSA, the SDDBSA and the SPS algorithm. In
fact, CTB$($SDBSA$)$ is approximately $0.57$\ times smaller than CTB$($SPS$)$%
; moreover, the CTB is reduced by $6.8\%$ if the SDDBSA is employed in place
of the SDBSA. It is also interesting to note that CTB$($DBF$)$ is
approximately $0.55$ times smaller than CTB$($MPF$)$ for the same value of $%
N_{p}$; once again, this result is in agreement with the results shown in 
\cite{Vitetta_DiViesti_2019} for a different SSM.

All the numerical results illustrated above lead to the conclusion that, in
the considered scenario, the DDBSA and the SDDBSA achieve the best
accuracy-complexity tradeoff in their categories of smoothing techniques.

\begin{figure}[tbp]
	\centering
	\includegraphics[width=0.6\textwidth]{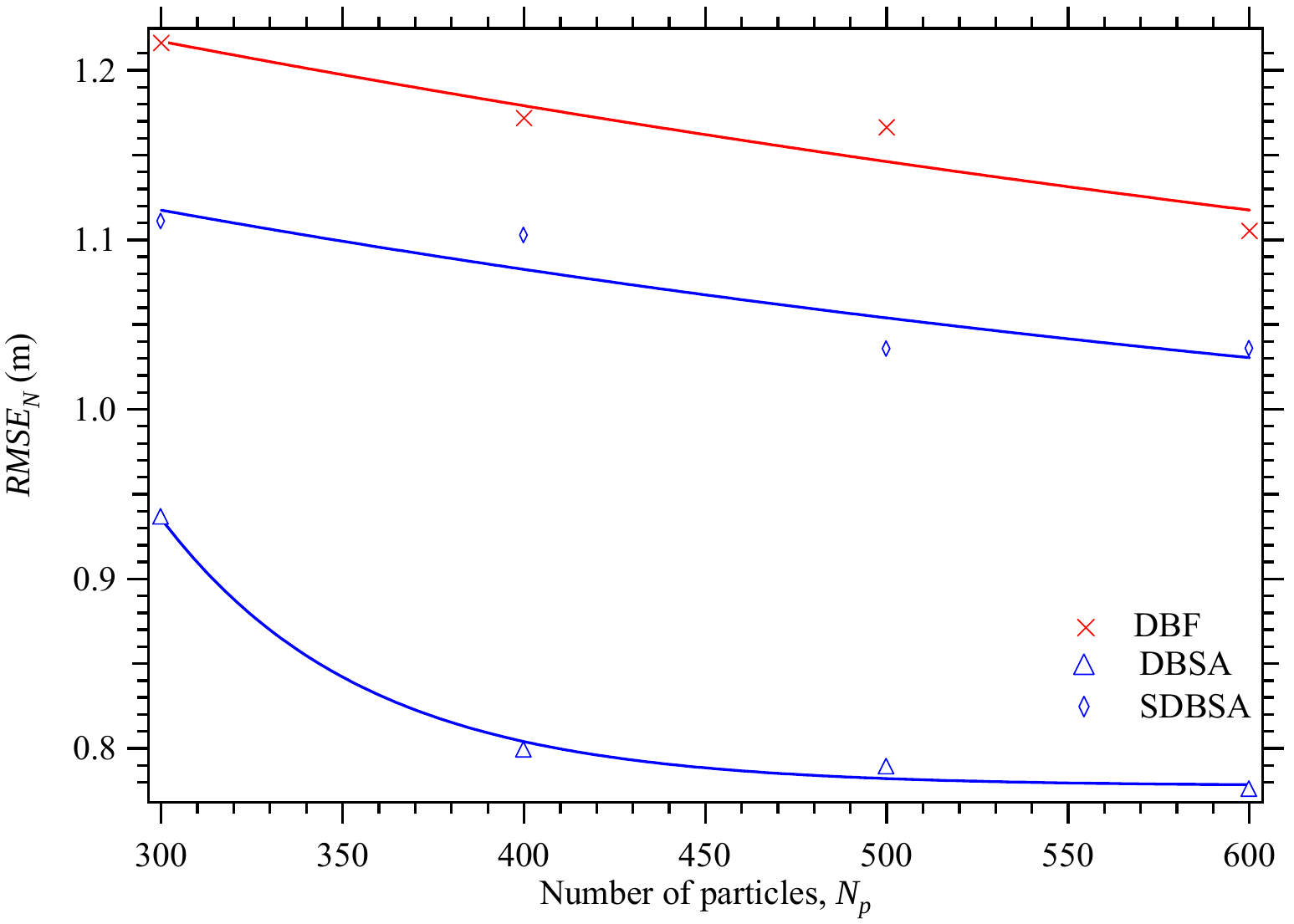}
	\caption{RMSE performance versus $N_{p}$ for the nonlinear component ($%
		RMSE_{N}$) of the state of SSM \#2; two smoothing algorithms (the DBSA and
		the SDBSA) and one filtering technique (DBF) are considered.}
	\label{Fig_rmsen2}
\end{figure}

The second SSM we considered is the same as the second SSM illustrated in 
\cite[Sec. IV]{Vitetta_DiViesti_2019} and refers to a sensor network
employing $P$ sensors placed on the vertices of a square grid (partitioning
a square area having side equal to $l$ m); these sensors receive the
reference signals radiated, at the same power level and at the same
frequency, by $N$ independent targets moving on a plane. Each target evolves
according to the motion model described by Eqs. (\ref{mod_1_v})-(\ref%
{mod_1_p}) with $\mathbf{a}_{k}(\mathbf{p}_{k})=\mathbf{0}$ for any $k$. In
this case, the considered SSM\ (denoted SSM\#2 in the following) refers to
the whole set of targets and its state vector $\mathbf{x}_{k}$ results from
the ordered concatenation of the vectors $\{\mathbf{x}_{k}^{(i)}$; $i=1$, $2$%
, $...$, $N\}$, where $\mathbf{x}_{k}^{(i)}\triangleq \lbrack (\mathbf{v}%
_{k}^{(i)})^{T},(\mathbf{p}_{k}^{(i)})^{T}]^{T}$, and \ $\mathbf{v}_{k}^{(i)}
$ and $\mathbf{p}_{k}^{(i)}$ represent the $i-$th target \emph{velocity} and
the \emph{position}, respectively. Moreover, the following additional
assumptions have been made about this SSM: 1) the process noises $\mathbf{n}%
_{p,k}^{(i)}$ and $\mathbf{n}_{v,k}^{(i)}$, affecting the $i-$th target
position and velocity, respectively, are given by $\mathbf{n}%
_{p,k}^{(i)}=(T_{s}^{2}/2)\,\mathbf{n}_{a,k}^{(i)}$ and $\mathbf{n}%
_{v,k}^{(i)}=T_{s}\mathbf{\,n}_{a,k}^{(i)}$, where $\{\mathbf{n}%
_{a,k}^{(i)}\}$ is two-dimensional AWGN, representing a random acceleration
and\ having covariance matrix $\sigma _{a}^{2}\,\mathbf{I}_{2}$ (with $i=1$, 
$2$, $...$, $N$); 2) the measurement acquired by the $q-$th sensor (with $q=1
$, $2$, $...$, $P$) in the $k$-th observation interval is given by 
\begin{equation}
	y_{q,k}=10\,\mathrm{log}_{10}\left( {\Psi }\sum_{i=1}^{N}\frac{\,d_{0}^{2}}{%
		\left\vert \left\vert \mathbf{s}_{q}-\mathbf{p}_{k}^{(i)}\right\vert
		\right\vert ^{2}}\right) +e_{k},  \label{mod_2_y}
\end{equation}%
where the measurement noise $\{e_{k}\}$ is AWGN with variance $\sigma
_{e}^{2}$, ${\Psi }$ denotes the normalised power received by each sensor
from any target at a distance $d_{0}$ from the sensor itself and $\mathbf{s}%
_{q}$ is the position of the considered sensor; 3) the overall measurement
vector $\mathbf{y}_{k}$ results from the ordered concatenation of the
measurements $\{y_{q,k}$; $q=1$, $2$, $...$, $P\}$ and, consequently,
provides information about the position only; 4) the initial position $%
\mathbf{p}_{0}^{(i)}\triangleq \lbrack p_{x,0}^{(i)},p_{y,0}^{(i)}]^{T}$ and
the initial velocity $\mathbf{v}_{0}^{(i)}\triangleq \lbrack
v_{x,0}^{(i)},v_{y,0}^{(i)}]^{T}$ of the $i-$th target are randomly selected
(with $i=1$, $2$, $...$, $N$). As far as the last point is concerned, it is
important to mention that, in our computer simulations, distinct targets are
placed in different squares of the partitioned area in a random fashion;
moreover, the initial velocity of each target is randomly selected within
the interval $(v_{\mathrm{min}},v_{\mathrm{max}})$ in order to ensure that
the trajectories of distinct targets do not cross each other in the
observation interval. The following values have been selected for the
parameters of SSM\#2: $P=25$, $l=10^{3}$ m, $T_{s}=1$ s, $\rho =1$, $\sigma
_{a}^{2}=0.1$ m/s$^{2}$, $\sigma _{e}^{2}=-35$ dB, ${\Psi =1}$, $d_{0}=1$ m, 
$v_{\mathrm{min}}=0$ m/s and $v_{\mathrm{min}}=0.1$ m/s. Moreover, $N=3$
targets have been observed over a time interval lasting $T=60$ $T_{s}$ s.
Our computer simulations have aimed at evaluating the accuracy achieved by
the considered smoothing algorithms in tracking the position of all the
targets. In practice, such an accuracy has been assessed by estimating the
average RMSE referring to the estimates of the whole set $\{\mathbf{p}%
_{k}^{(i)}$; $i=1$, $2$, $3\}$; note that, if the $i-$th target is
considered, its position $\mathbf{p}_{k}^{(i)}$ represents the \emph{\
	nonlinear} component of the associated substate $\mathbf{x}_{k}^{(i)}$,
because of the nonlinear dependence of $\mathbf{y}_{k}$ on it (see Eq. (\ref%
{mod_2_y})). Our computer simulations have evidenced that, in the considered
scenario, the MPF and the SDBF techniques diverge frequently in the
observation interval (some numerical results about the probability of
divergence area available in \cite[Sec. IV]{Vitetta_DiViesti_2019});
unluckily, when this occurs, all the smoothing algorithms that employ these
techniques in their forward pass (namely, Alg-L, the RBSS algorithm, the SPS
algorithm, the DDBSA and the SDDBSA) are unable to recover from this event
and, consequently, are useless. The DBF technique, instead, thanks to its
inner redundancy, is still able to track all the targets. Moreover, the two
smoothing algorithms employing this technique in their forward pass (namely,
the DBSA and the SDBSA), are able to improve the accuracy of position
estimates in their backward pass; this is evidenced by Fig. \ref{Fig_rmsen2}%
, that shows the dependence of $RMSE_{N}$ on the overall number of particles
($N_{p}$) for the DBF technique, the DBSA and the SDBSA (the range $[300,600]
$ is considered for $N_{p}$). Note that the SDBSA is outperformed by the
DBSA in terms of $RMSE_{N}$; for instance, $RMSE_{N}($SDBSA$)$ is about $1.31
$ times larger than $RMSE_{N}($DBSA$)$ for $N_{p}=500$. However, this result
is achieved at the price of a significantly higher complexity; in fact, CTB$(
$SDBSA$)$ is approximately equal to $2\cdot 10^{-3}\cdot $CTB$($DBSA$)$. 

\section{Conclusions\label{sec:conc}}

In this manuscript, factor graph methods have been exploited to develop new
smoothing algorithms based on the interconnection of two Bayesian filters in
the forward pass and of two backward information filters in the backward
pass. This has allowed us to develop a new approximate method for Bayesian
smoothing, called \emph{double Bayesian smoothing}. Four double Bayesian
smoothers have been derived for the class of conditionally linear Gaussian
systems and have been compared, in terms of both accuracy and execution
time, with other smoothing algorithms for two specific dynamic models. Our
simulation results lead to the conclusion that the devised algorithms can
achieve a better complexity-accuracy tradeoff and a better tracking
capability than other smoothing techniques recently appeared in the
literature.

\section*{Acknowledgment}

We would like to thank the anonymous Reviewers for their constructive
comments, that really helped us to improve the quality of this manuscript.

\appendices

\section{\label{app:A}}

In this Appendix, the derivation of the expressions of various messages
evaluated in each of the three phases the DBSA consists of is sketched.

\textbf{Phase I} - Formulas (\ref{W_bp_x_l}) and (\ref{w_bp_x_l}), referring
to the message $\cev{m}_{1}(\mathbf{x}_{k})$ (\ref{m_bp}), can be easily
computed by applying eqs. (IV.6)-(IV.8) of ref. \cite[Table 4, p.1304]%
{Loeliger_2007} in their backward form (with $A{\rightarrow }\mathbf{I}_{D}$%
, $X{\rightarrow }\mathbf{F}_{k}\mathbf{x}_{k}$, $Z{\rightarrow }\mathbf{x}%
_{k+1}$ and $Y{\rightarrow }\mathbf{u}_{k}+\mathbf{w}_{k}$) and, then, eqs.
(III.5)-(III.6) of \cite[Table 3, p.1304]{Loeliger_2007} (with $A{\
	\rightarrow }\mathbf{F}_{k}$, $X{\rightarrow }\mathbf{x}_{k}$ and $Y{\
	\rightarrow }\mathbf{F}_{k}\mathbf{x}_{k}$).

\textbf{Phase II} -Step 1) The message $m_{2}^{(n)}(\mathbf{x}_{k})$ (\ref%
{m_pm_x_l}) results from merging, in the BIF$_{2}{\rightarrow }$BIF$_{1}$
block, the statistical information about the \emph{nonlinear} state
component conveyed by the message $m_{1}^{(n-1)}(\mathbf{x}_{k}^{(N)})$ (\ref%
{M_1_x_k_N}) (and, consequently, by its $N_{p}$ components $%
\{m_{1,j}^{(n-1)}(\mathbf{x}_{k}^{(N)})=W_{1,k,j}^{(n-1)}\,\delta (\mathbf{x}%
_{k}^{(N)}-\mathbf{x}_{k,j}^{(N)})\}$) with those provided by the
pseudo-measurement $\mathbf{z}_{k}^{(L)}$ (\ref{eq:z_L_l}) about the \emph{\
	linear} state component. The method employed for processing this
pseudo-measurement is the same as that developed for MPF and can be
summarised as follows (additional mathematical details can be found in \cite[%
Sec. IV, p. 1527]{Vitetta_2019}):

1) The particles $\mathbf{x}_{k,j}^{(N)}$ and $\mathbf{x}_{\mathrm{be}
	,k+1}^{(N)}$, conveyed by the messages $m_{1}^{(n-1)}(\mathbf{x}_{k}^{(N)})$
(\ref{M_1_x_k_N}) and $\cev{m}_{\mathrm{be}}(\mathbf{x}_{k+1}^{(N)})$ (\ref%
{mess_be_N_l}), respectively, are employed to compute the $j-$th realization 
$\mathbf{z}_{k,j}^{(L)}$ (\ref{PM_z_L}) of $\mathbf{z}_{k}^{(L)}$ for $j=1$, 
$2$, $...$, $N_{p}$.

2) The pseudo-measurement $\mathbf{z}_{k,j}^{(L)}$ (\ref{PM_z_L}) is
exploited to generate the (particle-dependent) pdf%
\begin{equation}
	f_{j}^{(n)}\left( \mathbf{x}_{k}^{(L)}\right) =\mathcal{\mathcal{N}}\left( 
	\mathbf{x}_{k}^{(L)};\mathbf{\tilde{\eta}}_{k,j},\mathbf{\tilde{C}}%
	_{k,j}\right) ,  \label{eq:message_pm_L_j_tris}
\end{equation}%
that conveys pseudo-measurement information about $\mathbf{x}_{k}^{(L)}$ for
any $j$; the covariance matrix $\mathbf{\tilde{C}}_{k,j}$ and the mean
vector $\mathbf{\tilde{\eta}}_{k,j}$ of this message are computed on the
basis of the precision matrix $\mathbf{\tilde{W}}_{k,j}$ (\ref{eq:W_pm_L_j}%
)\ and the transformed mean vector $\mathbf{\tilde{w}}_{k,j}$ (\ref%
{eq:w_pm_L_j}), respectively.

3) The messages $\{m_{1,j}^{(n-1)}(\mathbf{x}_{k}^{(N)})\}$ are merged with
the pdfs $\{f_{j}^{(n)}(\mathbf{x}_{k}^{(L)})\}$ to generate the message $%
m_{2}^{(n)}(\mathbf{x}_{k})$ (\ref{m_pm_x_l}). The approach\ we adopt to
achieve this result is based on the fact that the message $m_{1,j}^{(n-1)}(%
\mathbf{x}_{k}^{(N)})$ and the pdf $f_{j}^{(n)}(\mathbf{x}_{k}^{(L)})$ refer
to the same particle (i.e., to the $j-$th particle $\mathbf{x}_{k,j}^{(N)}$,
but provide \emph{complementary} information (since they refer to the two
different components of the overall state $\mathbf{x}_{k}$). This allows us
to condense the statistical information conveyed by the sets $%
\{m_{1,j}^{(n-1)}(\mathbf{x}_{k}^{(N)})\}$ and $\{f_{j}^{(n)}(\mathbf{x}%
_{k}^{(L)})\}$ in the joint pdf%
\begin{equation}
	f^{(n)}(\mathbf{x}_{k}^{(L)},\mathbf{x}_{k}^{(N)})\triangleq
	w_{p}\sum\limits_{j=1}^{N_{p}}m_{1,j}^{(n-1)}\left( \mathbf{x}%
	_{k}^{(N)}\right) \,f_{j}^{(n)}\left( \mathbf{x}_{k}^{(L)}\right) .
	\label{joint_pdfbis}
\end{equation}%
referring to the whole state $\mathbf{x}_{k}$. Then, the message $%
m_{2}^{(n)}(\mathbf{x}_{k})$ (\ref{m_pm_x_l}) is computed by projecting the
pdf $f^{(k)}(\mathbf{x}_{k}^{(L)},\mathbf{x}_{k}^{(N)})$ (\ref{joint_pdfbis}%
) onto a single Gaussian pdf having the same \emph{mean }and \emph{covariance%
}.

Steps 2 and 3) The expression (\ref{m_be1_x_laa}) of $\cev{m}%
_{3}^{(n)}\left( \mathbf{x}_{k}\right) $ represents a straightforward
application of formula no. 2 of ref. \cite[App. A, TABLE I]{Vitetta_2019}
(with $\mathbf{W}_{1}{\rightarrow }\mathbf{W}_{1,k}$, $\mathbf{W}_{2}{%
	\rightarrow }\mathbf{W}_{2,k}^{(n)}$, $\mathbf{w}_{1}{\rightarrow }\mathbf{w}%
_{1,k}$ and $\mathbf{w}_{2}{\rightarrow }\mathbf{w}_{2,k}^{(n)}$). The same
considerations apply to the derivation of the expression (\ref{m_sm_x_l}) of 
$m_{4}^{(n)}\left( \mathbf{x}_{k}\right) $.

Step 4) The expression (\ref{weight_3}) of the weight $w_{3,k,j}^{(n)}$ is
derived as follows. First, we substitute the expression (\ref{f_x_N}) of $f(%
\mathbf{x}_{k+1}^{(N)}|\mathbf{x}_{k}^{(N)},\mathbf{x}_{k}^{(L)})$, and the
expressions of the messages $\cev{m}_{be}(\mathbf{x}_{k+1}^{(N)})$ (\ref%
{mess_be_N_l}) and $m_{1}^{(n)}(\mathbf{x}_{k}^{(L)})$ (\ref{m_fe_L_EKF_2})
in the \emph{right-hand side} (RHS)\ of Eq. (\ref{weight_bpa}). Then, the
resulting integral is solved by applying formula no. 1 of \cite[App. A,
TABLE\ II]{Vitetta_2019} in the integration with respect to $\mathbf{x}%
_{k}^{(L)}$ and the sifting property of the Dirac delta function in the
integration with respect to $\mathbf{x}_{k+1}^{(N)}$.

Step 5) - The derivation of the expression (\ref{m_pm_x_N_l_j}) of the
weight $w_{2,k,j}^{(n)}$ is similar to that illustrated for the particle
weights originating from the pseudo-measurements in dual MPF and can be
summarised as follows (additional mathematical details can be found in ref. 
\cite[Sec. V, pp. 1528-1529]{Vitetta_2019}). Two different Gaussian
densities are derived for the random vector $\mathbf{z}_{k}^{(N)}$ (\ref%
{z_N_l}),\emph{\ conditioned on} $\mathbf{x}_{k}^{(N)}$. The expression of
the first density originates from the definition (\ref{z_N_l}) and from the
knowledge of the \emph{joint} pdf of $\mathbf{x}_{k}^{(L)}$ and $\mathbf{x}%
_{k+1}^{(L)}$; this joint density is obtained from: a) the statistical
information provided by the message $m_{1}^{(n)}(\mathbf{x}_{k}^{(L)})$ (\ref%
{m_fe_L_EKF_2}) and the pdf $\mathcal{\mathcal{N}}(\mathbf{x}_{k+1}^{(L)},%
\mathbf{\tilde{\eta}}_{ \mathrm{be},k+1}\mathbf{\tilde{C}}_{\mathrm{be}%
	,k+1}) $ (resulting from integrating out the dependence of $\cev{m}_{\mathrm{%
		be}}(\mathbf{x}_{k+1})$ (\ref{mess_be_l}) on $\mathbf{x}_{k}^{(N)}$); b) the
Markov model $f(\mathbf{x}_{k+1}^{(L)}|\mathbf{x}_{k}^{(N)},\mathbf{x}%
_{k}^{(L)})$ (\ref{f_x_L}). This leads to the pdf 
\begin{equation}
	f_{1}^{(n)}\left( \mathbf{z}_{k}^{(N)}\left\vert \mathbf{x}_{k}^{(N)}\right.
	\right) =\mathcal{N}\left( \mathbf{z}_{k}^{(N)};\mathbf{\check{\eta}}%
	_{z,k}^{(n)}\left( \mathbf{x}_{k}^{(N)}\right) ,\mathbf{\check{C}}%
	_{z,k}^{(n)}\left( \mathbf{x}_{k}^{(N)}\right) \right) ,  \label{f_z_N}
\end{equation}%
where 
\begin{equation}
	\mathbf{\check{\eta}}_{z,k}^{(n)}\left( \mathbf{x}_{k}^{(N)}\right) =\mathbf{%
		\ \tilde{\eta}}_{\mathrm{be},k+1}-\mathbf{A}_{k}^{(L)}\left( \mathbf{x}%
	_{k}^{(N)}\right) \,\mathbf{\tilde{\eta}}_{1,k}^{(n)}  \label{eta_mess_z_Na}
\end{equation}%
and%
\begin{equation}
	\mathbf{\check{C}}_{z,k}^{(n)}\left( \mathbf{x}_{k}^{(N)}\right) =\mathbf{\ 
		\tilde{C}}_{\mathrm{be},k+1}-\mathbf{A}_{k}^{(L)}\left( \mathbf{x}%
	_{k}^{(N)}\right) \mathbf{\tilde{C}}_{1,k}^{(n)}\,\left( \mathbf{A}%
	_{k}^{(L)}\left( \mathbf{x}_{k}^{(N)}\right) \right) ^{T}.
	\label{C_mess_Z_N_a}
\end{equation}%
The second pdf of $\mathbf{z}_{k}^{(N)}$, instead, results from the fact
that this vector $\mathbf{z}_{k}^{(N)}$ must equal the sum (\ref{z_N_l_bis}%
); consequently, it is given by%
\begin{equation}
	f_{2}\left( \mathbf{z}_{k}^{(N)}\left\vert \mathbf{x}_{k}^{(N)}\right.
	\right) =\mathcal{N}\left( \mathbf{z}_{k}^{(N)};\mathbf{f}_{k}^{(L)}\left( 
	\mathbf{x}_{k}^{(N)}\right) ,\mathbf{C}_{w}^{(N)}\right) .  \label{fZ_N_bis}
\end{equation}%
Given the pdfs (\ref{f_z_N}) and (\ref{fZ_N_bis}), the message $\vec{m}%
_{3}^{(n)}(\mathbf{x}_{k}^{(N)})$ is expressed by their \emph{correlation},
i.e. it is computed as 
\begin{eqnarray}
	\vec{m}_{3}^{(n)}(\mathbf{x}_{k}^{(N)}) &=&\int f_{1}^{(n)}\left( \mathbf{z}%
	_{k}^{(N)}\left\vert \mathbf{x}_{k}^{(N)}\right. \right) \cdot f_{2}\left( \mathbf{z}_{k}^{(N)}\left\vert \mathbf{x}%
	_{k}^{(N)}\right. \right) d\mathbf{z}_{k}^{(N)}.
\end{eqnarray}%
Substituting Eqs. (\ref{f_z_N}) and (\ref{fZ_N_bis}) in the RHS of the last
expression, setting $\mathbf{x}_{k}^{(N)}=\mathbf{x}_{k,j}^{(N)}$ and
applying formula no. 4 of ref. \cite[Table II]{Vitetta_2019} to the
evaluation of the resulting integral yields Eq. (\ref{m_pm_x_N_l_j}); note
that $\mathbf{\ \check{\eta}}_{z,k,j}^{(n)}$ (\ref{eta_mess_z_N}) and $%
\mathbf{\check{C}}_{z,k,j}^{(n)}$ (\ref{C_mess_Z_N_bis}) represent the
values taken on by $\mathbf{\check{\eta}}_{z,k}^{(n)}(\mathbf{x}_{k}^{(N)})$
(\ref{eta_mess_z_Na}) and $\mathbf{\check{C}}_{z,k}^{(n)}(\mathbf{x}%
_{k}^{(N)})$ (\ref{C_mess_Z_N_a}), respectively, for $\mathbf{x}_{k}^{(N)}=%
\mathbf{x}_{k,j}^{(N)}$.

Step 7) The expression (\ref{weight_5_b}) of the weight $w_{5,k,j}^{(n)}$ is
derived as follows. First, we substitute the expressions (\ref{m_fe_L_EKF_2}%
) and (\ref{f_y_N}) of $m_{1}^{(n)}(\mathbf{x}_{k}^{(L)})$ and $f(\mathbf{y}%
_{k}|\mathbf{x}_{k}^{(N)},\,\mathbf{x}_{k}^{(L)})$, respectively, in the RHS
of Eq. (\ref{eq:mess_ms_PF}). Then, solving the resulting integral (see
formula no. 1 of ref. \cite[App. A, TABLE\ II]{Vitetta_2019}) produces Eq. (%
\ref{eq:weight_before_resampling}). Finally, setting $\mathbf{x}_{k}^{(N)}=%
\mathbf{x}_{k,j}^{(N)}$ in Eq. (\ref{eq:weight_before_resampling}) yields
Eq. (\ref{weight_5_b}).

\textbf{Phase III} - The expression (\ref{m_be2_xb}) of the message $\cev{m}%
_{\mathrm{be}}\left( \mathbf{x}_{k}\right) $ results from the application of
formula no. 2 of ref. \cite[App. A, TABLE I]{Vitetta_2019} to Eq. (\ref%
{m_be2_x}).

\section{Computational complexity of the devised double\\ Bayesian smoothers\label{app:CDBFA}}

In this appendix, the computational complexity of the tasks accomplished 
\emph{in a single recursion} of backward filtering and smoothing of the DBSA
is assessed in terms of flops. Moreover, we comment on how the illustrated
results can be also exploited to assess the computational complexity of a
single recursion of the DDBSA. In the following, $\mathcal{C}_{\mathbf{H}}$, 
$\mathcal{C}_{\mathbf{B}}$, $\mathcal{C}_{\mathbf{F}}$, $\mathcal{C}_{%
	\mathbf{A}^{(L)}}$ and $\mathcal{C}_{\mathbf{A}^{(N)}}$, and $\mathcal{C}_{%
	\mathbf{g}}$, $\mathcal{C}_{\mathbf{f}^{(L)}}$, $\mathcal{C}_{\mathbf{f}%
	^{(N)}}$ and $\mathcal{C}_{\mathbf{f}_{k}}$denote the cost due to the
evaluation of the matrices $\mathbf{H}_{k}$, $\mathbf{B}_{k}$, $\mathbf{F}%
_{k}$, $\mathbf{A}_{k}^{(L)}(\mathbf{x}_{k}^{(N)})$ and $\mathbf{A}%
_{k}^{(N)}(\mathbf{x}_{k}^{(N)})$, and of the functions $\mathbf{g}_{k}(%
\mathbf{x}_{k}^{(N)})$, $\mathbf{f}_{k}^{(L)}(\mathbf{x}_{k}^{(N)})$, $%
\mathbf{f}_{k}^{(N)}(\mathbf{x}_{k}^{(N)})$ and $\mathbf{f}_{k}(\mathbf{x}%
_{k})$, respectively. Moreover, similarly as \cite{Hoteit_2016}, it is
assumed that the computation of the inverse of any covariance matrix
involves a Cholesky decomposition of the matrix itself and the inversion of
a lower or upper triangular matrix. Finally, it is assumed that the
computation of the determinant of any matrix involves a Cholesky
decomposition of the matrix itself and the product of the diagonal entries
of a triangular matrix.

\textbf{Phase I} - The overall computational cost of this task is evaluated
as (see Eqs. (\ref{W_bp_x_l})-(\ref{w_bp_x_l}) and (\ref{eq:W_pm_L_j})-(\ref%
{PM_z_L})) 
\begin{eqnarray}
	\mathcal{C}_{1} &=&\mathcal{C}_{\mathbf{W}_{1,k}}+\mathcal{C}_{\mathbf{w}%
		_{1,k}}+N_{p}\left( \mathcal{C}_{\mathbf{z}_{k,j}^{(L)}}+\mathcal{C}_{%
		\mathbf{\tilde{W}}_{k,j}}\right.  \notag \\
	&&\left. +\mathcal{C}_{\mathbf{\tilde{w}}_{k,j}}+\mathcal{C}_{\mathbf{\tilde{%
				C}}_{k,j}}+\mathcal{C}_{\tilde{\eta}_{k,j}}\right) \triangleq \mathcal{C}_{%
		\mathrm{bp}}^{(1)}.  \label{C_F1_BP}
\end{eqnarray}%
Moreover, we have that: 1) the cost $\mathcal{C}_{\mathbf{W}_{1,k}}$ is
equal to $\mathcal{C}_{\mathbf{F}}+26D^{3}/3-D^{2}/2+5D/6$ flops; 2) the
cost $\mathcal{C}_{\mathbf{w}_{1,k}}$ is equal to $4D^{3}+4D^{2}-2D$ flops
(the cost for computing $\mathcal{C}_{\mathbf{F}}$ has been already
accounted for at point 1)); 3) the cost $\mathcal{C}_{\mathbf{z}%
	_{k,j}^{(L)}} $ is equal to $\mathcal{C}_{\mathbf{f}^{(N)}}+D_{N}$ flops; 4)
the cost $\mathcal{C}_{\mathbf{\tilde{W}}_{k,j}}$ is equal to $\mathcal{C}_{%
	\mathbf{A}^{(N)}}+4D_{N}^{3}-2D_{N}^{2}$ flops; 5) the cost $\mathcal{C}_{%
	\mathbf{\tilde{w}}_{k,j}}$ is equal to $2D_{N}^{3}+D_{N}^{2}-D_{N}$ flops
(the cost for computing $\mathcal{C}_{\mathbf{A}^{(N)}}$ has been already
accounted for at point 4)); 6) the cost $\mathcal{C}_{\mathbf{\tilde{C}}%
	_{k,j}}$ is equal to $2D_{N}^{3}/3+3D_{N}^{2}/2+5D_{N}/6$ flops; 7) the cost 
$\mathcal{C}_{\tilde{\eta}_{k,j}}$ is equal to $2D_{N}^{2}-D_{N}$ flops. The
expressions listed at points 1)-2) can be exploited for the DDBSA too; in
the last case, however, $D_{N}=0$ and $D=D_{L}$ must be assumed.

\textbf{Phase II} - The overall computational cost of this task is evaluated
as 
\begin{eqnarray}
	\mathcal{C}_{2} &=&n_{i}\left( \mathcal{C}_{\mathrm{pm}^{(1)}}+\mathcal{C}_{%
		\mathrm{be}1^{(1)}}+\mathcal{C}_{\mathrm{sm}^{(1)}}+\mathcal{C}_{\mathrm{bp}%
		^{(2)}}+\right.   \notag \\
	&&\left. \mathcal{C}_{\mathrm{pm}^{(2)}}+\mathcal{C}_{\mathrm{ms}^{(2)}}+%
	\mathcal{C}_{\mathrm{be}2^{(2)}}+\mathcal{C}_{\mathrm{sm}^{(2)}}\right) .
	\label{C_2}
\end{eqnarray}%
The terms appearing in the RHS of the last equation can be computed as
follows. First of all, we have that 
\begin{equation}
	\mathcal{C}_{\mathrm{pm}^{(1)}}=\mathcal{C}_{\eta _{2,k}^{(n)}}+\mathcal{C}_{%
		\mathbf{C}_{2,k}^{(n)}},  \label{C_PM1}
\end{equation}%
where (see Eqs. (\ref{eta_pm_l_k})-(\ref{C_pm_l_k})): 1) the cost $\mathcal{C%
}_{\eta _{2,k}^{(n)}}$ is equal to $2N_{p}D-D$ flops; 2) the cost $\mathcal{C%
}_{\mathbf{C}_{2,k}^{(n)}}$ is equal to $%
5N_{p}D_{L}^{2}+4N_{p}D_{N}^{2}+4N_{p}D_{L}D_{N}+D_{L}^{2}+D_{N}^{2}+D_{L}D_{N}
$ flops. The expressions listed at points 1)-2) can be exploited for the
DDBSA too; in the last case, however, $D_{N}=0$ and $D=D_{L}$ must be
assumed.

The second term appearing in the RHS of Eq. (\ref{C_2}) is evaluated as 
\begin{equation}
	\mathcal{C}_{\mathrm{be}1^{(1)}}=\mathcal{C}_{\mathbf{C}_{3,k}^{(n)}}+%
	\mathcal{C}_{\eta _{3,k}^{(n)}}+\mathcal{C}_{\mathbf{W}_{3,k}^{(n)}}+%
	\mathcal{C}_{\mathbf{w}_{3,k}^{(n)}},  \label{C_BE11}
\end{equation}%
where (see Eqs. (\ref{C_be1_l_ka})-(\ref{eta_be1_l_ka})): 1) the cost $%
\mathcal{C}_{\mathbf{C}_{3,k}^{(n)}}$ is equal to $14D^{3}/3+D^{2}/2+5D/6$
flops; 2) the cost $\mathcal{C}_{\eta _{3,k}^{(n)}}$ is equal to $4D^{2}-D$
flops (the cost for computing $\mathcal{C}_{\mathbf{W}_{k}^{(n)}}$ has been
already accounted for at point 1)); 3) the cost $\mathcal{C}_{\mathbf{W}%
	_{3,k}^{(n)}}$ is equal to $2D^{3}/3+3D^{2}/2+5D/6$ flops; 4) the cost $%
\mathcal{C}_{\mathbf{w}_{3,k}^{(n)}}$ is equal to $2D^{2}-D$ flops. The
expressions listed at points 1)-4) can be exploited for the DDBSA too; in
the last case, however, $D_{N}=0$ and $D=D_{L}$ must be assumed.

The third term appearing in the RHS of Eq. (\ref{C_2}) is computed as 
\begin{equation}
	\mathcal{C}_{\mathrm{sm}^{(1)}}=\mathcal{C}_{\mathbf{W}_{4,k}^{(n)}}+%
	\mathcal{C}_{\mathbf{w}_{4,k}^{(n)}}+\mathcal{C}_{\mathbf{C}_{4,k}^{(n)}}+%
	\mathcal{C}_{\eta _{4,k}^{(n)}},  \label{C_SM1}
\end{equation}%
where (see Eqs. (\ref{W_sm_l_k})-(\ref{w_sm_l_k})): 1) the cost $\mathcal{C}%
_{\mathbf{W}_{4,k}^{(n)}}$ is equal to $D^{2}$ flops; 2) the cost $\mathcal{C%
}_{\mathbf{w}_{4,k}^{(n)}}$ is equal to $D$ flops; 3) the cost $\mathcal{C}_{%
	\mathbf{C}_{4,k}^{(n)}}$ is equal to $2D^{3}/3+3D^{2}/2+5D/6$ flops; 4) the
cost $\mathcal{C}_{\eta _{4,k}^{(n)}}$ is equal to $2D^{2}-D$ flops. The
expressions listed at points 1)-4) can be exploited for the DDBSA too; in
the last case, however, $D_{N}=0$ and $D=D_{L}$ must be assumed.

The fourth term appearing in the RHS of Eq. (\ref{C_2}) is given by 
\begin{equation}
	\mathcal{C}_{\mathrm{bp}^{(2)}}=N_{p}\left( \mathcal{C}_{\eta
		_{3,k,j}^{(N)}}+\mathcal{C}_{\mathbf{C}_{3,k,j}^{(N)}}+\mathcal{C}%
	_{D_{3,k,j}^{(n)}}+\mathcal{C}_{Z_{3,k,j}^{(n)}}\right) ,  \label{C_BP2}
\end{equation}%
where (see Eqs. (\ref{eta_bp})-(\ref{cov_bp}) and (\ref{D_3})-(\ref{Z_bp})):
1) the cost $\mathcal{C}_{\eta _{3,k,j}^{(N)}}$ is equal to $\mathcal{C}_{%
	\mathbf{A}^{(N)}}+\mathcal{C}_{\mathbf{f}^{(N)}}+2D_{L}D_{N}$ flops; 2) the
cost $\mathcal{C}_{\mathbf{C}_{3,k,j}^{(N)}}$ is equal to $%
2D_{L}^{2}D_{N}+2D_{L}D_{N}^{2}-D_{L}D_{N}$ flops (the cost for computing $%
\mathcal{C}_{\mathbf{A}^{(N)}}$ and $\mathcal{C}_{\mathbf{f}^{(N)}}$ has
been already accounted for at point 1)); 3) the cost $\mathcal{C}%
_{D_{3,k,j}^{(n)}}$ is equal to $D_{N}^{3}/3+D_{N}^{2}+5D_{N}/3+2$ flops; 4)
the cost $\mathcal{C}_{Z_{3,k,j}^{(n)}}$ is equal to $2D_{N}^{2}+2D_{N}-1$
flops.

The fifth term appearing in the RHS of Eq. (\ref{C_2}) is evaluated as 
\begin{eqnarray}
	\mathcal{C}_{\mathrm{pm}^{(2)}} &=&N_{p}\left( \mathcal{C}_{\check{\eta}%
		_{z,k,j}^{(n)}}+\mathcal{C}_{\mathbf{\check{C}}_{z,k,j}^{(n)}}+\mathcal{C}_{%
		\mathbf{\check{W}}_{2,k,j}^{(n)}}+\mathcal{C}_{\mathbf{\check{w}}%
		_{2,k,j}^{(n)}}+\mathcal{C}_{D_{2,k,j}^{(n)}}+\mathcal{C}_{Z_{2,k,j}^{(n)}}\right) ,
	\label{C_PM2}
\end{eqnarray}%
where (see Eqs. (\ref{Z_pm})-(\ref{w_pm_x_N_l_j})): 1) the cost $\mathcal{C}%
_{\check{\eta}_{z,k,j}^{(n)}}$ is equal to $\mathcal{C}_{\mathbf{A}%
	^{(L)}}+2D_{L}^{2}$ flops; 2) the cost $\mathcal{C}_{\mathbf{\check{C}}%
	_{z,k,j}^{(n)}}$ is equal to $4D_{L}^{3}-D_{L}^{2}$ flops (the cost for
computing $\mathcal{C}_{\mathbf{A}^{(L)}}$ has been already accounted for at
point 1)); 3) the cost $\mathcal{C}_{\mathbf{\check{W}}_{2,k,j}^{(n)}}$ is
equal to $2D_{L}^{3}/3+5D_{L}^{2}/2+5D_{L}/6$ flops; 4) the cost $\mathcal{C}%
_{\mathbf{\check{w}}_{2,k,j}^{(n)}}$ is equal to $\mathcal{C}_{\mathbf{f}%
	^{(L)}}+4D_{L}^{2}-D_{L}$ flops; 5) the cost $\mathcal{C}_{D_{2,k,j}^{(n)}}$
is equal to $D_{L}^{3}/3+2D_{L}^{2}+5D_{L}/3+2$ flops; 6) the cost $\mathcal{%
	C}_{Z_{2,k,j}^{(n)}}$ is equal to $6D_{L}^{2}+3D_{L}-1$ flops.

The sixth term appearing in the RHS of Eq. (\ref{C_2}) is computed as%
\begin{equation}
	\mathcal{C}_{\mathrm{ms}^{(2)}}=N_{p}\left( \mathcal{C}_{\bar{\eta}%
		_{5,k,j}^{(n)}}+\mathcal{C}_{\mathbf{\bar{C}}_{5,k,j}^{(n)}}+\mathcal{C}%
	_{D_{5,k,j}^{(n)}}+\mathcal{C}_{Z_{5,k,j}^{(n)}}\right) ,  \label{C_MS2}
\end{equation}%
where (see Eqs. (\ref{eq:ms_fe_PF})-(\ref{Z_5})): 1) the cost $\mathcal{C}_{%
	\bar{\eta}_{5,k,j}^{(n)}}$ is equal to $\mathcal{C}_{\mathbf{B}}+\mathcal{C}%
_{\mathbf{g}}+2PD_{L}$ flops; 2) the cost $\mathcal{C}_{\mathbf{\bar{C}}%
	_{5,k,j}^{(n)}}$ is equal to $2PD_{L}^{2}+2P^{2}D_{L}-PD_{L}$ flops (the
cost for computing $\mathcal{C}_{\mathbf{B}}$ has been already accounted for
at point 1)); 3) the cost $\mathcal{C}_{D_{5,k,j}^{(n)}}$ is equal to $%
D_{L}^{3}/3+D_{L}^{2}+5D_{L}/3+2$ flops; 4) the cost $\mathcal{C}%
_{Z_{5,k,j}^{(n)}}$ is equal to $2P^{3}/3+7P^{2}/2+17P/6-1$ flops. It is
important to note that, if the forward weights $\{w_{\mathrm{fe},k,j}\}$ are
reused, the cost $\mathcal{C}_{\mathrm{ms}^{(2)}}$ appearing in Eq. (\ref%
{C_MS2}) is equal to zero.

The seventh term appearing in the RHS of Eq. (\ref{C_2}) is given by 
\begin{equation}
	\mathcal{C}_{\mathrm{be}2^{(2)}}=N_{p}\left( \mathcal{C}_{D_{6,k,j}^{(n)}}+%
	\mathcal{C}_{Z_{6,k,j}^{(n)}}+\mathcal{C}_{w_{6,k,j}^{(n)}}\right) ,
	\label{C_BE22}
\end{equation}%
where the costs $\mathcal{C}_{D_{6,k,j}^{(n)}}$ and $\mathcal{C}%
_{Z_{6,k,j}^{(n)}}$ are equal to $2$ flops, and the cost $\mathcal{C}%
_{w_{6,k,j}^{(n)}}$ is equal to $3$ flops (see Eqs. (\ref{weight_6b})-(\ref%
{Z_tot})). If the forward weights $\{w_{\mathrm{fe},k,j}\}$ are reused, the
costs $\mathcal{C}_{D_{6,k,j}^{(n)}}$ and $\mathcal{C}_{Z_{6,k,j}^{(n)}}$
are equal to $1$ flops, whereas the cost $\mathcal{C}_{w_{6,k,j}^{(n)}}$
remains unchanged.

The last term appearing in the RHS of Eq. (\ref{C_2}) is evaluated as 
\begin{equation}
	\mathcal{C}_{\mathrm{sm}^{(2)}}=\mathcal{C}_{w_{1,k,j}^{(n)}}+\mathcal{C}%
	_{W_{1,k,j}^{(n)}},  \label{C_SM2}
\end{equation}%
where the costs $\mathcal{C}_{w_{1,k,j}^{(n)}}$ and $\mathcal{C}%
_{W_{1,k,j}^{(n)}}$ are equal to $N_{p}$ and $2N_{p}-1$ flops, respectively
(see Eqs. (\ref{w_sm_j_Na})-(\ref{W_fe_2_x_N_l})).

\textbf{Phase III} - The overall computational cost of this task is
evaluated as 
\begin{equation}
	\mathcal{C}_{3}=\mathcal{C}_{\mathrm{be}^{(2)}}+\mathcal{C}_{\mathrm{pm}%
		^{(1)}}+\mathcal{C}_{\mathrm{be}1^{(1)}}+\mathcal{C}_{\mathrm{be}^{(1)}}.
	\label{C_3}
\end{equation}%
Here, the cost $\mathcal{C}_{\mathrm{be}^{(2)}}$ is equal to $\mathcal{C}%
_{S}(N_{p})$, that represents the total cost of a sampling step that
involves a particle set of size $N_{p}$; moreover, the costs $\mathcal{C}_{%
	\mathrm{pm}^{(1)}}$ and $\mathcal{C}_{\mathrm{be}1^{(1)}}$ are the same as
those appearing in the RHS of Eq. (\ref{C_2}), and $\mathcal{C}_{\mathrm{be}%
	^{(1)}}$ is computed as (see Eqs. (\ref{W_ms_x})-(\ref{w_be2_x})) 
\begin{eqnarray}
	\mathcal{C}_{\mathrm{be}^{(1)}} &=&\mathcal{C}_{\mathbf{W}_{\mathrm{ms},k}}+%
	\mathcal{C}_{\mathbf{w}_{\mathrm{ms},k}}+\mathcal{C}_{\mathbf{W}_{\mathrm{be}%
			2,k}}+ \\
	&&\mathcal{C}_{\mathbf{w}_{\mathrm{be}2,k}}+\mathcal{C}_{\mathbf{C}_{\mathrm{%
				be}}}+\mathcal{C}_{\eta _{\mathrm{be}}}.  \notag  \label{C_be1}
\end{eqnarray}%
Moreover, we have that: 1) the cost $\mathcal{C}_{\mathbf{W}_{\mathrm{ms},k}}
$ is equal to $\mathcal{C}_{\mathbf{H}}+2P^{2}D+2PD^{2}-D^{2}-PD$ flops; 2)
the cost $\mathcal{C}_{\mathbf{w}_{\mathrm{ms},k}}$ is equal to $\mathcal{C}%
_{\mathbf{B}}+\mathcal{C}_{\mathbf{g}}+2P^{2}D+3PD+2PD_{L}-P-D$ flops (the
cost for computing $\mathcal{C}_{\mathbf{H}}$ has been already accounted for
at point 1)); 3) the cost $\mathcal{C}_{\mathbf{W}_{\mathrm{be}2,k}}$ is
equal to $D^{2}$ flops; 4) the cost $\mathcal{C}_{\mathbf{w}_{\mathrm{be}%
		2,k}}$ is equal to $D$ flops; 5) the cost $\mathcal{C}_{\mathbf{C}_{\mathrm{%
			be}}}$ is equal to $2D^{3}/3+3D^{2}/2+5D/6$ flops; 6) the cost $\mathcal{C}%
_{\eta _{\mathrm{be}}}$ is equal to $2D^{2}-D$ flops. The expressions listed
at points 1)-6) can be exploited for the DDBSA too; in the last case,
however, $D_{N}=0$ and $D=D_{L}$ must be assumed. Note that the costs $%
\mathcal{C}_{\mathbf{W}_{\mathrm{ms},k}}$ and $\mathcal{C}_{\mathbf{w}_{%
		\mathrm{ms},k}}$ (see points 1) and 2)) are ignored if the precision matrix $%
\mathbf{W}_{\mathrm{ms},k}$ and the transformed mean vector $\mathbf{w}_{%
	\mathrm{ms},k}$ are stored in the forward pass (so that they do not need to
be recomputed in the backward pass). Moreover, if the SDBSA or the SDDBSA is
used, the cost $\mathcal{C}_{\mathrm{be}^{(2)}}$ in the RHS of Eq. (\ref{C_3}%
) becomes $D_{N}(2N_{p}-1)$ flops.

Finally, it is worth stressing that, if the DBSA or the DDBSA (the SDBSA or
the SDDBSA) is employed, the overall computational complexity is obtained by
multiplying the computational cost assessed for a single recursion by $M\,T$
(by $T$), where $M$ and $T$ denote the overall number of accomplished
backward passes and the duration of the observation interval, respectively.


\begin{thebibliography}{99}
	
	\bibitem{Anderson_1979} B. Anderson and J. Moore, \textbf{Optimal Filtering}%
	, Englewood Cliffs, NJ, Prentice-Hall, 1979.
	
	\bibitem{Sar_2013} S. S\"{a}rkk\"{a}, \textbf{Bayesian Filtering and
		Smoothing}. Cambridge, U.K.: Cambridge Univ. Press, 2013.
	
	\bibitem{Doucet_2000} A. Doucet, S. Godsill and C. Andrieu,
	\textquotedblleft On Sequential Monte Carlo Sampling Methods for Bayesian
	Filtering\textquotedblright, \emph{Statist. Comput.}, vol. 10, no. 3, pp.
	197-208, 2000.
	
	\bibitem{Kitagawa_1987} G. Kitagawa, \textquotedblleft Non-Gaussian
	state-space modeling of nonstationary time series\textquotedblright, \emph{%
		Journal of the American Statistical Association}, vol. 82, pp. 1032-1063,
	1987.
	
	\bibitem{Kitagawa_1994} G. Kitagawa, \textquotedblleft The two-filter
	formula for smoothing and an implementation of the Gaussian-sum
	smoother\textquotedblright, \emph{Annals of the Institute of Statistical
		Mathematics}, vol. 46, pp. 605-623, 1994.
	
	\bibitem{Bresler_1986} Y. Bresler, \textquotedblleft Two-filter formula for
	discrete-time non-linear Bayesian smoothing\textquotedblright, \emph{Int.
		Journal of Control}, vol. 43, no. 2, pp. 629-641, 1986.
	
	\bibitem{Vo_2012} B. N. Vo, B. T. Vo and R. P. S. Mahler, \textquotedblleft
	Closed-Form Solutions to Forward--Backward Smoothing\textquotedblright, 
	\emph{IEEE Trans. Sig. Proc.}, vol. 60, no. 1, pp. 2-17, Jan. 2012.
	
	\bibitem{Sar_2010} S. S\"{a}rkk\"{a} and J. Hartikainen, \textquotedblleft
	On Gaussian optimal smoothing of non-linear state space
	models\textquotedblright, \emph{IEEE Trans. Autom. Control}, vol. 55, no. 8,
	pp. 1938--1941, Aug. 2010.
	
	\bibitem{Kok_2016} J. Kokkala, A. Solin, and S. S\"{a}rkk\"{a},
	\textquotedblleft Sigma-point filtering and smoothing-based parameter
	estimation in nonlinear dynamic systems\textquotedblright, \emph{J. Adv.
		Inf. Fusion}, vol. 11, no. 1, pp. 15--30, 2016.
	
	\bibitem{Gar_2017} A. F. Garc\'{\i}a-Fern\'{a}ndez, L. Svensson and S. S\"{a}%
	rkk\"{a}, \textquotedblleft Iterated posterior linearisation
	smoother\textquotedblright, \emph{IEEE Trans. Autom. Control}, vol. 62, no.
	4, pp. 2056--2063, Apr. 2017.
	
	\bibitem{Douc_2011} R. Doucet, A. Garivier, E. Moulines and J. Olsson,
	\textquotedblleft Sequential Monte Carlo smoothing for general state space
	hidden Markov models\textquotedblright, \emph{Ann. Appl. Probab.}, vol. 21,
	no. 6, pp. 2109--2145, 2011.
	
	\bibitem{Kitagawa_1996} G. Kitagawa, \textquotedblleft Monte Carlo filter
	and smoother for non-Gaussian nonlinear state space
	models\textquotedblright, \emph{J. Comput. Graph. Statist.}, vol. 5, no. 1,
	pp. 1--25, 1996.
	
	\bibitem{Godsill_2004} S. J. Godsill, A. Doucet, and M. West,
	\textquotedblleft Monte Carlo smoothing for nonlinear time
	series\textquotedblright, \emph{J. Amer. Statist. Assoc.}, vol. 99, no. 465,
	pp. 156--168, Mar. 2004.
	
	\bibitem{Lindsten_2013} F. Lindsten and T. B. Sch\"{o}n, \textquotedblleft
	Backward simulation methods for Monte Carlo statistical
	inference\textquotedblright, \emph{Foundat. Trends Mach. Learn.}, vol. 6,
	no. 1, pp. 1--143, 2013.
	
	\bibitem{Hot_2019} R. Hostettler and S. S\"{a}rkk\"{a}, \textquotedblleft
	Rao--Blackwellized Gaussian Smoothing\textquotedblright, \emph{IEEE Trans.
		Autom. Control}, vol. 64, no. 1, pp. 305-312, Jan. 2019.
	
	\bibitem{Briers_2010} M. Briers, A. Doucet and S. Maskell, \textquotedblleft
	Smoothing algorithms for state-space models\textquotedblright, \emph{Ann.
		Inst. Statist. Math.}, vol. 62, no. 1, pp. 61--89, Feb. 2010.
	
	\bibitem{Fong_2002} W. Fong, S. J. Godsill, A. Doucet and M. West,
	\textquotedblleft Monte Carlo smoothing with application to audio signal
	enhancement\textquotedblright, \emph{IEEE Trans. Signal Process.}, vol. 50,
	no. 2, pp. 438--449, Feb. 2002.
	
	\bibitem{Lindsten_2016} F. Lindsten, P. Bunch, S. S\"{a}rkk\"{a}, T. B. Sch%
	\"{o}n and S. J. Godsill, \textquotedblleft Rao-Blackwellized Particle
	Smoothers for Conditionally Linear Gaussian Models\textquotedblright, \emph{%
		IEEE J. Sel. Topics in Sig. Proc.}, vol. 10, no. 2, pp. 353-365, March 2016.
	
	\bibitem{Vitetta_2018} G. M. Vitetta, E. Sirignano and F. Montorsi,
	\textquotedblleft Particle Smoothing for Conditionally Linear Gaussian
	Models as Message Passing over Factor Graphs\textquotedblright, \emph{IEEE
		Trans. Sig. Proc.}, vol. 66, no. 14, pp. 3633-3648, July 2018.
	
	\bibitem{Vitetta_DiViesti_2019} G. M. Vitetta, P. Di Viesti and E.
	Sirignano, \textquotedblleft Multiple Bayesian Filtering as Message
	Passing\textquotedblright, submitted to the \emph{IEEE Trans. Sig. Proc.},
	February 2019 (available on arXiv at https://arxiv.org/abs/1907.01358)
	
	\bibitem{Loeliger_2007} H.-A. Loeliger, J. Dauwels, Junli Hu, S. Korl, Li
	Ping, F. R. Kschischang, \textquotedblleft The Factor Graph Approach to
	Model-Based Signal Processing\textquotedblright, \emph{IEEE Proc.}, vol. 95,
	no. 6, pp. 1295-1322, June 2007.
	
	\bibitem{Kschischang_2001} F. R. Kschischang, B. Frey, and H. Loeliger,
	\textquotedblleft Factor Graphs and the Sum-Product
	Algorithm\textquotedblright, \emph{IEEE Trans. Inf. Theory}, vol. 41, no. 2,
	pp. 498-519, Feb. 2001.
	
	\bibitem{Vitetta_2019} G. M. Vitetta, E. Sirignano, P. Di Viesti and F.
	Montorsi, \textquotedblleft Marginalized Particle Filtering and Related
	Techniques as Message Passing\textquotedblright, \emph{IEEE Trans. Sig. Proc.%
	}, vol. 67, no. 6, pp. 1522-1535, Mar. 2019.
	
	\bibitem{Schon_2005} T. Sch\"{o}n, F. Gustafsson, P.-J. Nordlund,
	\textquotedblleft Marginalized Particle Filters for Mixed Linear/Nonlinear
	State-Space Models\textquotedblright, \emph{IEEE Trans. Sig. Proc.}, vol.
	53, no. 7, pp. 2279-2289, July 2005.
	
	\bibitem{Arulampalam_2002} M. S. Arulampalam, S. Maskell, N. Gordon and T.
	Clapp, \textquotedblleft A Tutorial on Particle Filters for Online
	Nonlinear/Non-Gaussian Bayesian Tracking\textquotedblright, \emph{IEEE
		Trans. Sig. Proc.}, vol. 50, no. 2, pp. 174-188, Feb. 2002.
	
	\bibitem{Hoteit_2016} B. Ait-El-Fquih and I. Hoteit, \textquotedblleft A
	variational Bayesian multiple particle filtering scheme for
	large-dimensional systems\textquotedblright, \emph{IEEE Trans. Sig. Proc.},
	vol. 64, no. 20, pp. 5409--5422, Oct. 2016.

\end{thebibliography}
\end{document}